\documentclass[11pt]{report}

\usepackage{epigraph}
\usepackage[dvips]{graphicx}
\usepackage{pstricks}
\usepackage[margin=1.5in]{geometry}
\usepackage{amsmath}
\usepackage{amsfonts}
\usepackage{amssymb}
\usepackage{dsfont} 
\usepackage{xcolor}
\usepackage{mathtools}
\usepackage{hyperref} 
\usepackage[T1]{fontenc}

\usepackage[francais]{babel}

\makeatletter
\renewcommand{\@seccntformat}[1]{ {\csname
the#1\endcsname}.\hspace{0.2em}}
\makeatother

\def \P {{\mathbb P}}
\def \C {{\mathbb C}}
\def\E{\mathbb E}
\def \R {{\mathbb R}}
\def \N {{\mathbb N}}
\def \Z {{\mathbb Z}}

\def \Ad {\mbox{Ad}}
\def \ad {\mbox{ad}}
\def \SU {\mbox{SU}}
\def \rad {\mbox{rad}}

\DeclareMathOperator{\Pf}{Pf} 

\newcommand{\PP}{\mathcal P}

\begin{document}
 \begin{titlepage}
\begin{center}
{  \Large{\sc  Universit\'e de Paris}}\\\vspace{2.5cm}
{\LARGE\bf {\Large Habilitation \`a diriger des recherches}}\\\vspace{1cm}{\Large{Sp\'ecialit\'e : Math\'ematiques}}\\\vspace{1.5cm}
\vspace{2cm}

{\LARGE{\sc Mouvement brownien et  alg\`ebres de \\  \vspace{0.4cm} Kac--Moody affines  }}\\\vspace{0.2cm}
 
 \vspace{1cm}
 
{\large \sc Manon Defosseux}\\\vspace{1.5cm} 
 
\begin{tabular}{llll}
Rapporteurs~: & Philippe &  
BIANE &  Universit\'e Paris-Est Marne-la-Vall\'ee  \\ & Persi & DIACONIS & Stanford University
\\ & Thierry & 
L\'EVY & Sorbonne Universit\'e   
\end{tabular}
\\
 \end{center}
 
\vspace{0.5cm}
\noindent{\, Soutenue   le  6 novembre 2020  devant le jury compos\'e de :}
{
\vspace{0.1cm}
\begin{center}
\begin{tabular}{lll} \vspace{0.2cm}
Philippe &  
BIANE &  Universit\'e Paris-Est  Marne-la-Vall\'ee \\  \vspace{0.2cm}
Philippe  &  
BOUGEROL & Sorbonne Universit\'e\\ \vspace{0.2cm}
Persi & DIACONIS & Stanford University   \\ \vspace{0.2cm}
Catherine  & 
DONATI-MARTIN & Universit\'e de Versailles Saint-Quentin-en-Yvelines   \\ \vspace{0.2cm}
Nathalie  & 
EISENBAUM & Universit\'e de  Paris   \\ \vspace{0.2cm}
Alice  & 
GUIONNET &  \'Ecole normale sup\'erieure de Lyon   \\ \vspace{0.2cm}
 Jean-Fran\c cois   & LE GALL  & Universit\'e Paris-Saclay    \\ \vspace{0.2cm}
Thierry & 
L\'EVY & Sorbonne Universit\'e   
\end{tabular}
\end{center}}
  \end{titlepage}

 \newtheorem{theo}{Th\'eor\`eme}[section]
\newtheorem{lemm}[theo]{Lemme} 
\newtheorem{defn}[theo]{D\'efinition} 
\newtheorem{prop}[theo]{Proposition}
\newtheorem{coro}[theo]{Corollaire} 
\newtheorem{rema}[theo]{Remarque}

\definecolor{gris25}{gray}{0.75}
 
\newcounter{encart}
\renewcommand{\theencart}{{}}
\newenvironment{encart}[1]
{
\refstepcounter{encart}
\small \sf \medskip
\noindent\colorbox{gris25}{
\makebox[\textwidth][c]{{\sf --\, P{\scriptsize{UBLICATIONS}}  
\theencart} -- \bfseries #1}}\medskip 
\addcontentsline{toc}{subsection}{\sf \small \hspace{0.85 cm} Encart \theencart. #1}
}
{\smallskip
\noindent\colorbox{gris25}{\makebox[\textwidth][c]{\hspace{1cm}}}
\rm \normalsize \medskip}

\newenvironment{theobis}[1]
  {\renewcommand{\thetheo}{\ref{#1}$\ bis $}%
   \addtocounter{theo}{-1}%
   \begin{theo}}
  {\end{theo}}
  \newpage
  \thispagestyle{empty}    
           $ $
     \newpage                 \thispagestyle{empty}           \quad \quad\quad\quad \quad  \quad\quad\quad\quad \quad  \quad     \quad  \quad \quad  \quad \quad {\it Pour Philippe et Philippe,  Joseph et Cl\'eo,}
                                    
                                    \quad\quad\quad\quad   \quad \quad\quad   \quad\quad\quad  \quad     \quad  \quad \quad  \quad \quad {\it  Agathe, Florent et Martin Gardner  } 
                                                                                  
                         \newpage      
                         $ $ \thispagestyle{empty}   
                         \newpage                           
  \setcounter{tocdepth}{1}
  \setcounter{page}{1}
\tableofcontents
\newpage
\enlargethispage{1.4cm} 
 {\bf{\Large Remerciements}} 
  \thispagestyle{empty} 
 \\

 De l'autre cot\'e du miroir, une fillette vient de lire un texte bien myst\'erieux, dont le sens se d\'evoile et se d\'erobe tour \`a tour, et voil\`a ce que nous dit Lewis Carroll :
 \og \c Ca a l'air tr\`es joli, dit Alice, quand elle eut fini de lire, mais c'est assez difficile \`a comprendre !    [...]   \c Ca me remplit la t\^ete de toutes sortes d'id\'ees, mais\dots\, mais je ne sais pas exactement quelles sont ces id\'ees ! En tout cas, ce qu'il y a de clair c'est que quelqu'un a tu\'e quelque chose\footnote{Serait-ce un brownien qu'on a tu\'e  ?  }\dots \fg. Au moment d'\'ecrire les remerciements, je pense \'evidemment d'abord \`a Philippe Bougerol et \`a nos innombrables conversations. Elles m'ont appris  quel engrais puissant peuvent devenir, pour celui qui les re\c coit, certaines  intuitions livr\'ees  dans un flou   g\'en\'ereux et stimulant\footnote{flou pour celui qui \'ecoute bien s\^ur, et  bien souvent accompagn\'e de r\'ef\'erences pr\'ecises pour se mettre au travail.}. Je remercie Philippe  Bougerol de n'avoir  cess\'e depuis que nous nous connaissons de me donner mati\`ere \`a penser, c'est peu dire que les math\'ematiques que je fais lui doivent \'enorm\'ement. 
 \\
 
 J'ai toujours eu une immense admiration pour les travaux et la  pratique math\'ematique de Philippe Biane, Persi Diaconis et Thierry L\'evy.  Cela a \'et\'e pour moi une joie immense  qu'ils acceptent tous les trois de rapporter mon habilitation, et je les en remercie.  Je remercie \'egalement Catherine Donati-Martin, Nathalie Eisenbaum, Alice Guionnet et Jean-Fran\c cois Le Gall  d'avoir accept\'e de faire partie de mon Jury. Je serai heureuse et  honor\'ee de   leur pr\'esenter mon travail.  \\
 
Charme, beaut\'e et cr\'eativit\'e, merci \`a C\'eline et Maya pour les r\'eunions joyeuses et caf\'ein\'ees qui ont accompagn\'e le tout d\'ebut de nos r\'edactions.  Merci \`a Reda, \`a Philippe encore, pour les am\'eliorations que  leurs relectures attentives du manuscrit ont permises,  \`a Thierry \`a nouveau, pour ses nombreuses suggestions  et pour la  d\'elicatesse de leur formulation. Et  merci  \`a tous ceux, en particulier Nathael, Nathalie et C\'eline encore, qui ont permis que les derni\`eres \'etapes   soient franchies avec succ\`es.
  \\ 
    
  Nul besoin de caravelles pour les  d\'ecouvertes math\'ematiques grandes ou modestes, et les voyages  et les temp\^etes sont int\'erieurs dans le bureau du math\'ematicien. Voil\`a quelques ann\'ees  que Rapha\"el et moi voguons ensemble, cinq jours par semaine, \`a bord du vaisseau 732-E. Son calme, sa mesure (de Lebesgue of course), sa gentillesse et sa bonne humeur   auront rendu  les escales plus douces et  les intemp\'eries moins p\'enibles, et je l'en remercie. Merci \'evidemment   aux autres passagers du MAP5 qui nous accompagnent et qui  auront donn\'e \`a la travers\'ee  eux aussi   une couleur  bien   joyeuse.  
\\
 
 Avant de finir, merci au papillon du Br\'esil (il faut n'oublier personne) et au l\'epidopt\'eriste qui ne l'a pas \'epingl\'e : sans eux assur\'ement, les \'ev\'enements auraient pris une tout autre tournure.
\\

Merci enfin \`a Philippe qui m'accompagne de son chant depuis de nombreuses ann\'ees.  \'El\'egies et ritournelles font  un autre Jabberwocky.

 \newpage
  \thispagestyle{empty}   
 $ $
 \newpage 

 \chapter*{Introduction}  
\begin{figure}
     \begin{center}
            \includegraphics[scale=0.9]{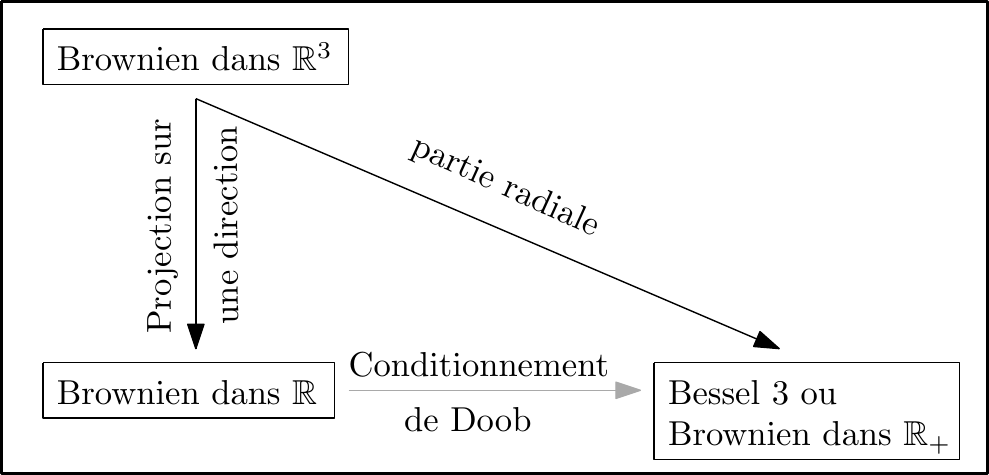}
            \caption{Transformation de Doob et Bessel $3$}
            \label{DR3}
                    \end{center} 
 \end{figure}

\begin{figure}
     \begin{center}
            \includegraphics[scale=0.9]{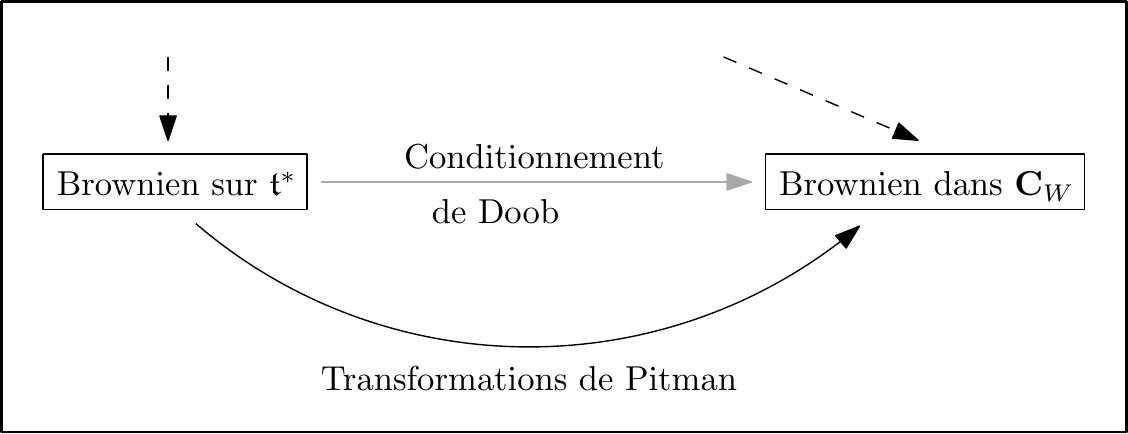}
            \caption{Transformation de Doob et transformations de Pitman - Le cas compact}
            \label{DCHP} 
                    \end{center}
 \end{figure}

\addcontentsline{toc}{chapter}{Introduction} 
  En 2014, j'ai commenc\'e \`a r\'efl\'echir aux liens entre mouvement brownien et alg\`ebre de Kac--Moody affine. Ce sont ces liens que je me propose d'exposer dans ce m\'emoire, tels qu'ils me sont apparus depuis. On trouvera dans cette introduction non pas l'\'etat de l'art dans mon domaine de recherche mais la description de mes contributions dans un champ que je m'efforcerai de pr\'eciser dans le corps du m\'emoire. Ces contributions ont toutes un analogue dans un cadre souvent qualifi\'e de cadre compact ou semi-simple, tandis qu'elles interviennent dans un cadre que l'on   qualifiera de cadre affine.
 Je commence  par rappeler les r\'esultats du cadre compact dont mes contributions sont un analogue. Davantage que sur les r\'esultats   eux m\^emes   c'est sur leur structure que je souhaite insister, telle qu'elle perdure dans un cadre affine.    
 \paragraph{\Large{Le cadre compact.}} Les r\'esultats que nous rappelons pour le cadre compact sont tous bien connus.
 \paragraph{Norme d'un brownien de dimension trois.}  On appelle partie radiale d'un vecteur de $\R^3$ la norme euclidienne de ce vecteur. Les orbites pour l'action du groupe $\mbox{SO(3)}$ sur $\R^3$ sont les sph\`eres de $\R^3$ et la partie radiale d'un vecteur d\'etermine l'orbite  \`a laquelle il appartient. Si on consid\`ere un mouvement brownien standard sur $\R^3$ alors le processus de sa partie radiale  est un processus de Bessel de dimension trois. C'est aussi une transformation de Doob d'un brownien  r\'eel tu\'e en $0$, via la fonction harmonique $h$ d\'efinie sur $\R_+$ par $h(x)=x$, $x\ge 0$.  La projection dans une direction d'un brownien de $\R^3$ est un brownien r\'eel et il y a ainsi deux fa\c cons, illustr\'ees \`a la figure \ref{DR3}, d'obtenir un Bessel $3$ \`a partir d'un brownien de dimension trois.  Sur le diagramme - c'est une convention que nous utiliserons toujours -  les fl\`eches noires figurent une transformation d\'eterministe, tandis que la grise d\'esigne une op\'eration sur des mesures. Pour finir rappelons que la loi en un temps fix\'e du brownien  dans $\R^3$ conditionnellement \`a sa partie radiale  est la mesure de probabilit\'e uniforme sur l'orbite  correspondante.  On a pour la transform\'ee de   Fourier de la probabilit\'e uniforme $\mu_r$ sur la sph\`ere $S^2_r$ de rayon $r$   la formule suivante
 \begin{align}\label{archi} \int_{S^2_r}e^{i \lambda(u,z)}\mu_r(dz)=\frac{\sin(r\lambda)}{ r\lambda},\quad \lambda\in \R,
\end{align}
pour tout vecteur $u\in \R^3$ de norme $1$. Autrement dit la coordonn\'ee dans une direction quelconque donn\'ee d'un point choisi selon la mesure uniforme sur une sph\`ere de rayon $r$ est uniform\'ement distribu\'ee sur l'intervalle $[-r,r]$\footnotemark\footnotetext{ \samepage Ce r\'esultat est connu sous le nom de th\'eor\`eme d'Archim\`ede. C'est l'autre pouss\'ee d'Archim\`ede : pouss\'ee vers l'avant de la mesure uniforme sur une sph\`ere par la projection sur une direction.}. Tout ceci est bien connu. Nous allons maintenant changer  de point de vue et d\'ecrire ces propri\'et\'es dans un cadre un peu diff\'erent, qui pr\'esente l'avantage de mieux mettre en \'evidence les propri\'et\'es structurelles qui nous int\'eressent. 
 \paragraph{Brownien sur le dual de $\mathfrak{su}(2)$.}
 Consid\'erons l'ensemble de matrices $\mathfrak{su}(2)$ d\'efini par
 $$\mathfrak{su}(2)=\{M\in \mathcal M_2(\C): M+M^*=0, \mbox{tr}(M)=0\},$$
 une base de $\mathfrak{su}(2)$
 $$x=\begin{pmatrix} 
0& 1 \\
-1& 0
\end{pmatrix},\, y=\begin{pmatrix} 
0& i \\
i& 0
\end{pmatrix}, \, z=\begin{pmatrix} 
i& 0 \\
0& -i
\end{pmatrix}, $$
et sa base duale $\{e_x,e_y,e_z\}$.  On d\'efinit une action coadjointe\footnotemark\footnotetext{C'est \`a dessein que nous travaillons ici sur le dual de $\mathfrak{su}(2)$ et non sur  $\mathfrak{su}(2)$. En effet, si en dimension finie, il est indiff\'erent de consid\'erer un espace vectoriel ou son dual, il sera essentiel dans le cadre infini-dimensionnel qui est le n\^otre de consid\'erer, conform\'ement \`a l'esprit de la m\'ethode des orbites de Kirillov, non pas les orbites adjointes mais les orbites coadjointes.} $\Ad^*$ de $\SU(2)$ sur le dual $\mathfrak{su}(2)^*$ de $\mathfrak{su}(2)$ en posant pour $\varphi\in \mathfrak{su}(2)^*$, $u\in \SU(2)$, $M \in \mathfrak{su}(2)$
$$(\Ad^*(u)\varphi)(M)=\varphi(u^{-1}Mu).$$
Pour toute forme $\varphi$ dans $\mathfrak{su}(2)^*$ il existe un unique r\'eel positif $r$ et un \'el\'ement $u\in \SU(2)$ tels que 
$$\Ad^*(u)(\varphi)=re_z.$$
 Autrement dit  l'espace quotient $\mathfrak{su}(2)^*/\Ad^*(\SU(2))$ s'identifie \`a $\R_+e_z$. On appelle $re_z$ la partie radiale de $\varphi$. En fait, si 
$$\varphi=ae_x+be_y+ce_z,$$
pour $a,b,c\in \R$,  alors $r=\sqrt{a^2+b^2+c^2}$.  Munissons maintenant $\mathfrak{su}(2)^*$ du produit scalaire usuel et consid\'erons un brownien standard $\{b_t: t\ge 0\}$ sur $\mathfrak{su}(2)^*$, i.e.
$$b_t=x_te_x+y_te_y+z_te_z,\quad t\ge 0,$$
o\`u $\{(x_t,y_t,z_t): t\ge 0\}$ est un brownien standard sur $\R^3$. Dans ce contexte, la projection sur une direction d'un vecteur de $\R^3$ devient  la projection  sur le dual d'une sous-alg\`ebre de Cartan de $\mathfrak{su}(2)$ et les propri\'et\'es du brownien de dimension trois rappel\'ees plus haut s'\'enoncent de la fa\c con suivante.  D'une part, le processus de la partie radiale de $\{b_t: t\ge 0\}$ et celui de sa projection sur $\R e_z$ tu\'ee en $0$ et  conditionn\'ee au sens de Doob \`a rester dans $\R_+e_z$ ont m\^eme loi. D'autre part, en un temps fix\'e $t>0$, la loi de $b_t$ conditionnellement \`a son $\Ad^*(\SU(2))$-orbite est la mesure uniforme sur cette orbite. 
 L'avantage d'une telle pr\'esentation est qu'on peut imm\'ediatement en donner une traduction valable lorsque  $\SU(2)$ est remplac\'e par un groupe de Lie compact connexe semi-simple quelconque.  Pour indiquer comment  consid\'erons un tel groupe    $G_0$ et son alg\`ebre de Lie   $\mathfrak{g}_0$. On suppose  sans 
 perte de g\'en\'eralit\'e que $G_0$ est un groupe de matrices.
 \begin{figure}
     \begin{center}
            \includegraphics[scale=0.9]{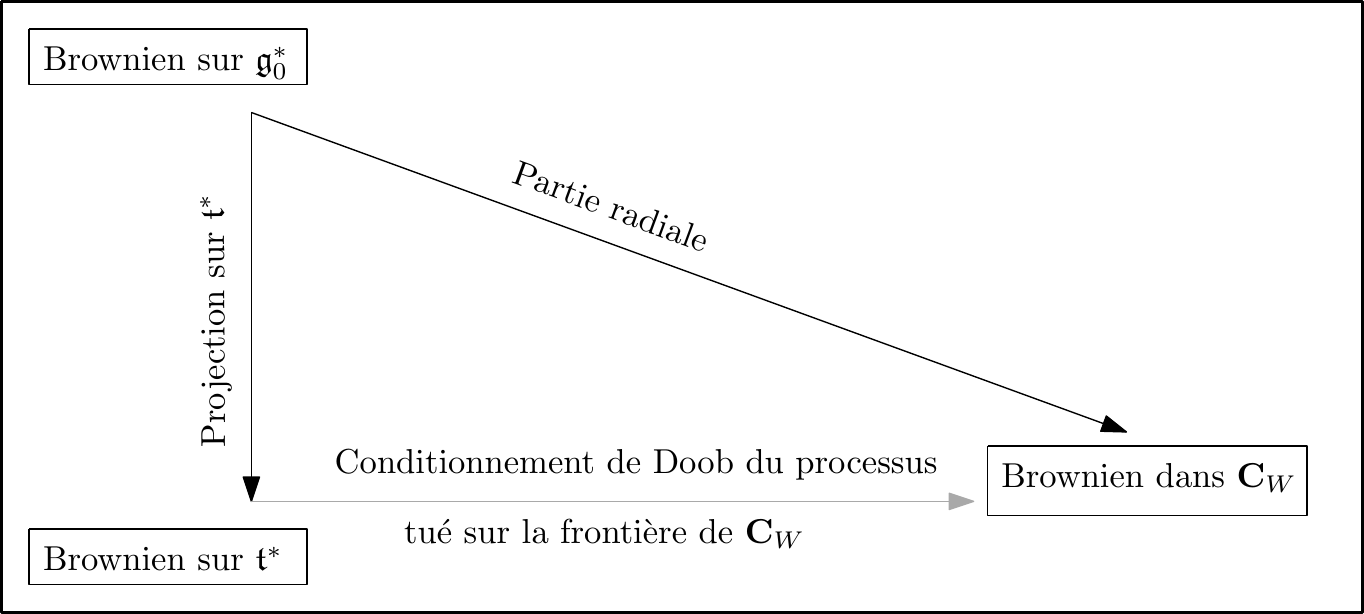}
            \caption{Partie radiale et transformation de Doob - le cas compact}
            \label{DCH} 
                    \end{center}
 \end{figure}
\paragraph{Brownien sur $\mathfrak{g}_0^*$.}  L'action coadjointe not\'ee $\Ad^*$ de $G_0$ sur $\mathfrak g_0^*$ est d\'efinie de m\^eme que celle de $\SU(2)$ sur $\mathfrak{su}(2)^*$.   On consid\`ere un tore maximal $T$ de $G_0$ et son alg\`ebre de Lie $\mathfrak t$ qui joue le r\^ole de $\R z$,  et on choisit une chambre de Weyl ${\bf C}_W$ dans $\mathfrak t^*$ qui joue le r\^ole de $\R_+e_z$. Il existe pour chaque classe d'\'equivalence dans  $\mathfrak{g}_0^*/\Ad^*(G_0)$ un unique repr\'esentant dans ${\bf C}_W$. Par exemple, lorsque $G_0$ est le groupe sp\'ecial unitaire $\SU(n)$, on peut choisir pour $T$ l'ensemble des matrices diagonales de $\SU(n)$ et pour ${\bf C}_W$ l'ensemble
$$\{\sum_{k=1}^{n-1}\lambda_ke_{k}^*: \lambda_1\ge \dots\ge \lambda_{n-1}\ge 0 \},$$ 
o\`u $e_{k}^*(x)=x_k$, si $x$ est une matrice     de taille $n\times n$ dont les \'el\'ements diagonaux sont $ix_1,\dots,ix_n$.    On appelle partie radiale d'une forme  sur $\mathfrak g_0$ le repr\'esentant dans ${\bf C}_W$ de l'orbite \`a laquelle elle appartient.   On munit $\mathfrak g^*_0$ d'un produit scalaire $\Ad^*(G_0)$-invariant. Alors le processus de la partie radiale d'un mouvement brownien standard sur $\mathfrak{g}^*_0$ et celui de sa projection  sur $\mathfrak t^*$ tu\'ee sur le bord de ${\bf C}_W$ et conditionn\'ee au sens de Doob \`a rester dans ${\bf C}_W$ ont  m\^eme loi. Le diagramme commutatif repr\'esent\'e en figure \ref{DCH}  illustre  cette observation. 
  En outre, la loi en un temps fix\'e du brownien sur $\mathfrak{g}_0^*$ conditionnellement \`a son $\Ad^*(G_0)$-orbite est la mesure de probabilit\'e uniforme sur cette orbite. Dans ce contexte, la formule   (\ref{archi}) devient une formule de Harish-Chandra, \'equivalente dans ce cadre compact \`a une formule des caract\`eres de Kirillov, ou \`a la formule de Itzykson--Zuber lorsque $G_0=\SU(n)$.  Elle donne la transform\'ee de Fourier de la mesure de Duistermaat--Heckman normalis\'ee associ\'ee \`a l'action   du tore $T$ sur une orbite coadjointe, c'est-\`a-dire la mesure image par la projection canonique sur $\mathfrak t^*$ de la mesure de probabilit\'e uniforme sur cette orbite. 
   \paragraph{Le th\'eor\`eme de Pitman.}
 La transformation de Pitman $\mathcal P$  op\`ere sur les chemins \`a valeurs r\'eelles, c'est-\`a-dire les  fonctions continues $f:\R_+\to \R$ telles que $f(0)=0$. Elle est  d\'efinie par
$$\mathcal Pf(t)=f(t)-2\inf_{0\le s\le t}f(s), \, \, t\ge 0.$$ 
Le th\'eor\`eme de Pitman s'\'enonce  ainsi : si $\{b(t),t\ge 0\}$ est un mouvement brownien r\'eel standard, alors $\{\mathcal Pb(t),t\ge 0\}$ est un processus de Bessel de dimension trois.  Ce th\'eor\`eme \'etablit donc une relation trajectorielle entre le brownien sur $\R$ et la transformation de Doob sur $\R_+$. Philippe Biane, Philippe Bougerol et Neil O'Connell ont montr\'e dans \cite{bbo} qu'une telle relation existait entre le mouvement brownien sur $\mathfrak t^*$ et celui dans  le c\^one ${\bf C}_W$.  Elle s'obtient en appliquant au brownien sur $\mathfrak t^*$ des transformations de type Pitman associ\'ees aux sym\'etries orthogonales par rapport aux hyperplans perpendiculaires aux racines simples de $\mathfrak g_0$.   On peut donc ajouter une fl\`eche au diagramme de la figure \ref{DCH} comme indiqu\'e \`a la figure \ref{DCHP}. Par ailleurs, la mesure normalis\'ee de Duistermaat--Heckman s'obtient en consid\'erant en un temps fix\'e la loi du brownien sur $\mathfrak t^*$ conditionnellement \`a son image par les transformations de Pitman.

 \paragraph{\Large{Le cadre affine}. }
 
 Tr\`es rapidement, un enjeu important de ma recherche   a \'et\'e de comprendre ce qu'un   diagramme tel que celui repr\'esent\'e en figure \ref{DCH} pouvait devenir dans un cadre affine. Dans un second temps, s'est pos\'ee la question de l'existence d'un th\'eor\`eme de type Pitman dans ce cadre. Dans la suite $G_0$ est un groupe de Lie suppos\'e compact, connexe, simple et simplement connexe et on consid\`ere   un produit scalaire   $(\cdot \vert \cdot)$ $\Ad(G_0)$-invariant sur  $\mathfrak g_0$.
 \paragraph{Alg\`ebres affines, alg\`ebres de lacets.} 
 Les alg\`ebres de Lie affines font partie d'une classe d'alg\`ebres de Lie appel\'ees alg\`ebres de Kac--Moody qui contient, outre les alg\`ebres de Lie semi-simples complexes, des alg\`ebres de Lie de dimension infinie. Bien que de dimension infinie les alg\`ebres affines  partagent avec les alg\`ebres de Lie semi-simples  de nombreuses propri\'et\'es.   Nous les pr\'esentons dans le chapitre \ref{chap-Algebre-affine}.     Dans ce m\'emoire nous consid\'erons une classe particuli\`ere d'alg\`ebres affines : les alg\`ebres de lacets \`a valeurs dans une alg\`ebre de Lie  simple complexe avec extension centrale. Dans la perspective d'un diagramme commutatif, on doit travailler sur la partie compacte d'une telle alg\`ebre. Nous consid\'erons donc l'alg\`ebre de lacets $L(\mathfrak g_0)$ \`a valeurs dans $\mathfrak{g}_0$ - les lacets sont index\'es par le cercle $S^1$ identifi\'e \`a $\R/\Z$ - et 
son extension centrale\footnotemark \footnotetext{Ici comme dans la suite nous ne  pr\'ecisons pas le degr\'e de r\'egularit\'e des lacets et supposons qu'ils sont toujours aussi r\'eguliers qu'on peut le souhaiter.}
 $$\widetilde{L}(\mathfrak{g}_0)= {L}(\mathfrak{g}_0)\oplus \R c,$$
 munie d'un crochet de Lie d\'efini par
 \begin{align}
[\xi+\lambda c,\eta+\mu c]=[\xi,\eta]_{\mathfrak g_0}+ \Big(\int_0^1(\xi'(s)\vert \eta(s))\, ds\,\Big)c,
\end{align}  
pour $\xi,\eta\in L(\mathfrak g_0)$, $\lambda,\mu\in \R$,  o\`u  $[\cdot ,\cdot]_{\mathfrak g_0}$ est le crochet de Lie sur ${\mathfrak g_0}$ et $[\xi,\eta]_{\mathfrak g_0}$ est  d\'efini point par point.  Le  crochet d\'efinit une action adjointe de $ {\widetilde{L}}(\mathfrak{g}_0)$ sur elle-m\^eme. En consid\'erant l'exponentielle de cette action   on d\'efinit l'action adjointe  d'un groupe de lacets  $L(G_0)$ \`a valeurs dans $G_0$  sur  $ { \widetilde{L}}(\mathfrak{g}_0)$. On consid\`ere l'action coadjointe $\Ad^*$ de $L(G_0)$ sur  $\widetilde{ {L}}(\mathfrak{g}_0)^*$ qui en d\'ecoule. On \'ecrit 
$$\widetilde{ {L}}(\mathfrak{g}_0)^*= L(\mathfrak{g}_0)^*\oplus \R\Lambda_0$$
o\`u $\Lambda_0$ est un poids dit fondamental d\'efini par 
$$\Lambda_0(c)=1, \quad \Lambda_0(L(\mathfrak{g}_0))=0.$$ 
L'alg\`ebre $\mathfrak t\oplus \R c$ est une sous-alg\`ebre ab\'elienne maximale de $\widetilde L(\mathfrak g_0)$ et son dual est $\mathfrak t^*\oplus \R\Lambda_0$.
Il nous faut maintenant d\'ecrire ce qui remplace dans ce cadre les objets du diagramme de la figure \ref{DCH}. 
 
 \paragraph{Browniens espace-temps.} Le premier travail effectu\'e dans cette voie est expos\'e dans le chapitre \ref{chap-BETaff}, qui reprend les r\'esultats de \cite{defo2} en les modifiant l\'eg\`erement, de sorte qu'ils r\'epondent \`a une exigence de coh\'erence. Il  a permis d'obtenir des candidats pour les processus   de la base du diagramme. Pour cela, je me suis inspir\'ee de ce que nous savions du cas compact. Dans ce cas en effet, comme on peut le voir dans les travaux de Philippe Biane \cite{biane1}, le diagramme peut se comprendre comme une d\'eg\'en\'erescence commutative d'un diagramme analogue valable pour des variables al\'eatoires non commutatives. Ces variables al\'eatoires sont par nature li\'ees \`a la th\'eorie des repr\'esentations de $\mathfrak g_0$ et dans un   contexte non commutatif le brownien sur $\mathfrak t^*$ est remplac\'e par une marche sur le r\'eseau des     poids   de $\mathfrak g_0$ et le brownien dans ${\bf C}_W$ par un processus de Markov sur celui de ses poids dominants. Les pas de la marche sont distribu\'es selon une mesure de probabilit\'e uniforme sur l'ensemble des poids d'une repr\'esentation  complexe  de dimension finie de $\mathfrak g_0$ et le noyau du processus de Markov s'exprime en fonction de la dimension des composantes isotypiques d'un produit tensoriel de repr\'esentations. Un th\'eor\`eme central limite relie ces processus \`a temps discret au brownien sur $\mathfrak t^*$ d'une part  et au mouvement brownien dans ${\bf C}_W$ d'autre part.   La th\'eorie des repr\'esentations des alg\`ebres de Kac--Moody affines poss\`ede de nombreux points communs avec celle des repr\'esentations des alg\`ebres de Lie semi-simples complexes et moyennant quelques ajustements on peut mimer la construction de ces processus dans un cadre affine\footnotemark \footnotetext{La construction est en fait valable pour une alg\`ebre de lacets avec extension centrale augment\'ee d'une d\'erivation. Il est cependant inutile d'ajouter la d\'erivation pour cet expos\'e liminaire, car si elle joue un grand r\^ole dans la construction des processus, il n'en reste pas trace dans l'\'enonc\'e des r\'esultats.}.  On obtient alors une marche sur le r\'eseau des poids d'une alg\`ebre affine et un processus de Markov sur celui de ses poids dominants.  Le cas affine diff\`ere cependant   du cas compact sur plusieurs points essentiels.  Le plus important \`a mentionner ici  est sans doute la pr\'esence pour les processus discrets introduits et leurs limites en temps long d'une coordonn\'ee temporelle d\'eterministe le long du poids   fondamental $\Lambda_0$.  \begin{figure}
     \begin{center}
            \includegraphics[scale=0.9]{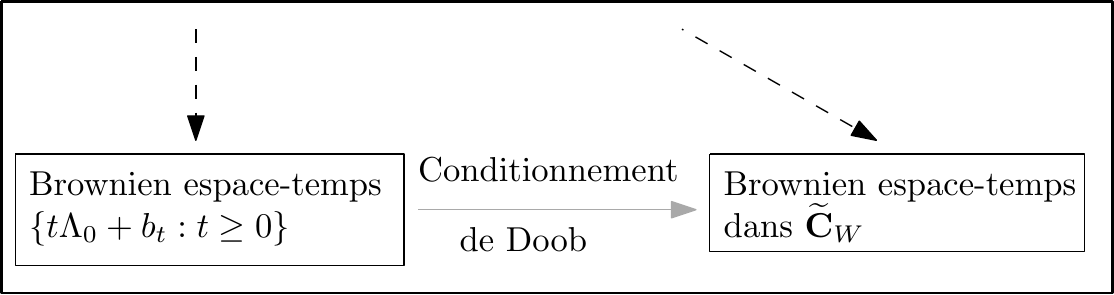}
            \caption{Processus et processus conditionn\'es - Le cas affine}
            \label{DCFlg1} 
        \end{center}
 \end{figure} C'est le niveau des poids des repr\'esentations consid\'er\'ees qui appara\^it ainsi.    Finalement, l'\'etude des limites en temps long des processus discrets nous permet de placer sur la base du diagramme, \`a gauche un processus espace-temps $$\{t\Lambda_0+b_t: t\ge 0\}$$ o\`u $\{b_t: t\ge 0\}$ est un brownien standard sur $\mathfrak t^*$, \`a droite  un processus espace-temps conditionn\'e au sens de Doob \`a rester dans une chambre de Weyl fondamentale ${\bf \widetilde {C}}_W$ dans    $\R_+\Lambda_0+ \mathfrak t^*$.   \begin{figure}[h]
     \begin{center}
            \includegraphics[scale=0.9]{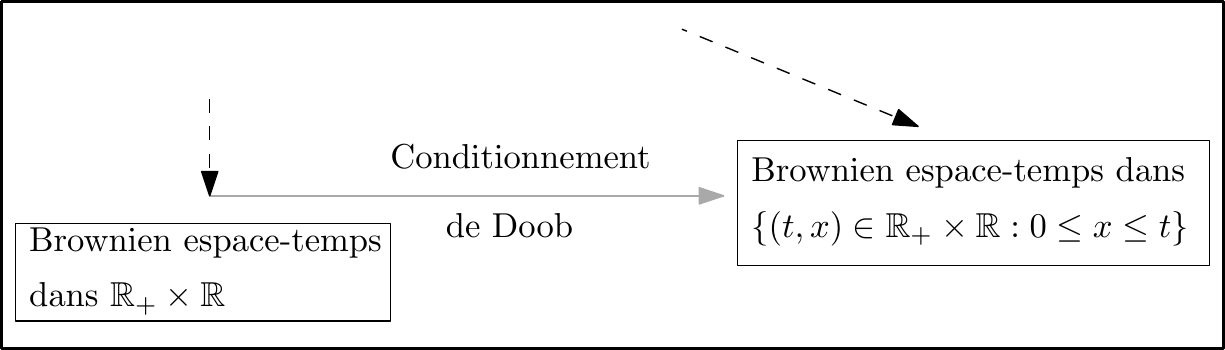}
            \caption{Processus et processus conditionn\'e  - Le cas $A_1^{(1)}$}
            \label{DCFlg1A} 
                    \end{center}
 \end{figure} Nous obtenons ainsi la base du diagramme  repr\'esent\'ee \`a la figure \ref{DCFlg1}. 
 Dans le cas o\`u $\mathfrak g_0=\mathfrak{su}(2)$, l'alg\`ebre affine est de type $A_1^{(1)}$ et  les ensembles  $\R_+\Lambda_0+\mathfrak t^*$ et  ${\bf \widetilde {C}}_W$ s'identifient respectivement \`a $\R_+\times \R$ et $$\{(t,x)\in \R_+\times \R: 0\le x\le t\}.$$ On repr\'esente \`a  la figure  \ref{DCFlg1A}  la base du diagramme dans ce cas.

\begin{figure}
     \begin{center}
            \includegraphics[scale=0.9]{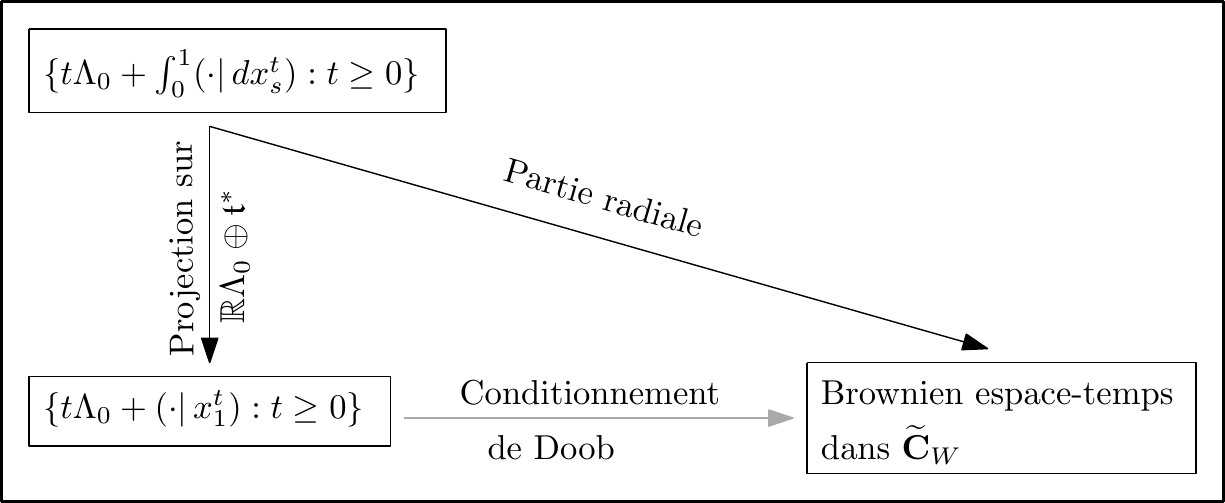}
            \caption{Partie radiale et conditionnement de Doob - Le cas affine}
            \label{DCA} 
               \end{center}
 \end{figure}

\paragraph{Drap brownien.} Une fois obtenue une premi\`ere ligne du diagramme commutatif, il restait \`a d\'eterminer deux processus, l'un devant jouer le r\^ole du brownien sur $\mathfrak g^*_0$ et l'autre celui du processus de sa partie radiale dans ${\bf C}_W$. C'est ce que j'ai   fait dans \cite{defo3}.  Dans \cite{Frenkel}, Igor  Frenkel d\'efinit une notion de partie radiale pour les formes de $\R\Lambda_0\oplus L(\mathfrak g_0)^*$ s'\'ecrivrant 
$$\Phi=t\Lambda_0+\int_0^1(\,\cdot\, \vert \,\dot x(s))\, ds,$$
avec $t$ un r\'eel strictement positif et $\{x(s) :s\in [0,1]\}$  un chemin r\'egulier \`a valeurs dans $\mathfrak g_0$, c'est-\`a-dire $$\Phi(c)=t\textrm{ et }\Phi(y)=\int_0^1(y(s)\vert \,\dot x(s)) \, ds,$$ pour $y\in L(\mathfrak g_0)$. Dans ce cas, il existe un unique \'el\'ement dans l'intersection de   la chambre de Weyl fondamentale ${\bf \widetilde {C}}_W$  et de l'orbite $\Ad^*(L(G_0))\{\Phi\}$. C'est cet \'el\'ement qu'on appelle  partie radiale de $\Phi$. Il s'obtient en r\'esolvant l'\'equation diff\'erentielle  (en $s$)
$$t \, dX(s)=X(s) \, dx(s),$$
avec $X(0)=I$, o\`u $I$ est la matrice identit\'e dans $G_0$ et en consid\'erant l'orbite de $X(1)$ dans $G_0$ pour l'action par conjugaison de $G_0$ sur lui-m\^eme.  C'est cette orbite dans $G_0$ qui d\'etermine l'orbite de $\Phi$ dans $\widetilde L(\mathfrak g_0)^*$. Nous l'avons dit, dans le cas compact, la mesure uniforme sur une orbite peut s'obtenir en consid\'erant  la loi d'une gaussienne sur $\mathfrak g_0^*$ conditionnellement \`a sa partie radiale. Dans \cite{Frenkel} Igor  Frenkel a l'id\'ee, pour d\'efinir une mesure sur une orbite coadjointe dans un cadre affine, de remplacer la trajectoire r\'eguli\`ere  par un mouvement brownien $\{x_s: s\in [0,1]\} $ sur $\mathfrak g_0$, de consid\'erer le processus $\{X_s: s\in[0,1]\}$  issu de $I$,  solution de l'\'equation diff\'erentielle stochastique
$$t \, d X  =X\circ dx,$$
o\`u $\circ$ d\'esigne l'int\'egrale de Stratonovitch\footnotemark\footnotetext{Une telle solution est un processus \`a valeurs dans $G_0$.} et de d\'efinir une mesure sur une orbite coadjointe \`a partir de la loi de $\{x_s: s\in[0,1]\}$ conditionnellement \`a l'orbite de $X_1$ dans $G_0$. C'est pour une telle mesure que Frenkel \'etablit une formule des caract\`eres de Kirillov dans un cadre affine. Cette formule  donne la transform\'ee de Fourier de la mesure image  par la projection sur $\R\Lambda_0\oplus \mathfrak t^*$ de la mesure de Frenkel sur une orbite.   Ainsi dans ce cadre c'est un couple form\'e d'un r\'eel positif $t$ et d'un brownien $\{x_s: s\in[0,1]\}$ qui joue le r\^ole de la gaussienne sur $\mathfrak g_0^*$ et  ce couple doit \^etre pens\'e comme une forme al\'eatoire\footnotemark\footnotetext{Dans notre travail, tous les calculs sont faits en utilisant le mouvement brownien ou le drap brownien sur $\mathfrak g_0$ qui sont bien d\'efinis. L'\'ecriture sous forme de forme al\'eatoire est utilis\'ee dans l'introduction comme dans le corps du m\'emoire     pour rappeler que les lois des trajectoires sur $\mathfrak g_0$ jouent dans le cadre affine le r\^ole de mesures sur des orbites coadjointes.}
\begin{align}\label{mesorbaff}
t\Lambda_0+\int_0^1(\, \cdot \,\vert\, dx_s).
\end{align}
Nous sommes d\'esormais en mesure de compl\'eter le diagramme de la figure \ref{DCH} dans un contexte affine. Au sommet se trouve un processus de L\'evy \`a valeur dans $\widetilde L(\mathfrak g_0)^*$ dont la projection sur $\R\Lambda_0\oplus\mathfrak t^*$ est un brownien espace-temps $\{t\Lambda_0+b_t:t\ge 0\}$ et qui en chaque temps fix\'e $t$, a m\^eme loi que (\ref{mesorbaff}) pour un brownien bien choisi. Un tel processus s'obtient en consid\'erant un drap brownien $\{x_s^t: s, t\ge 0\}$ sur $\mathfrak g_0$ et en lui associant le processus 
$$\{t\Lambda_0 +\int_0^1 (\, \cdot\, \vert \, dx_s^t) :t\ge 0\}.$$
Le processus de la partie radiale s'obtient en r\'esolvant pour chaque $t$ l'\'equation diff\'erentielle stochastique (en $s$)
$$t\, dX_s=X_s \circ dx_s^t,$$ avec $X_0=I$ et en consid\'erant l'orbite dans $G_0$ de $X_1$. On obtient alors le diagramme commutatif repr\'esent\'e en figure \ref{DCA}.

 \paragraph{Th\'eor\`eme de localisation et formule de Kirillov--Frenkel.} Faisons une parenth\`ese dans la pr\'esentation de nos contributions. Nous l'avons dit, Igor Frenkel \'etablit dans \cite{Frenkel} une formule des caract\`eres de type Kirillov dans le cadre des alg\`ebres affines. Dans un cadre compact, la formule des caract\`eres de Kirillov est au coeur du diagramme commutatif de la figure \ref{DCH}. Il en va de m\^eme pour notre diagramme, repr\'esent\'e en figure \ref{DCA}, et  la formule   des caract\`eres de Frenkel. Celle-ci permet en effet d'obtenir une relation d'entrelacement d'op\'erateurs \`a partir de laquelle on d\'emontre \`a la fa\c con de  Chris Rogers et Jim Pitman \cite{RP} l'identit\'e en loi entre le transform\'e de Doob et le processus de la partie radiale.   Comprendre le travail de Frenkel a donc  \'et\'e  pour mon propre travail une \'etape essentielle. Dans le cadre compact la formule de Kirillov peut se voir comme une formule de localisation de Duistermaat et de Heckman. Il m'a sembl\'e int\'eressant de pr\'esenter  celle de  Frenkel comme une formule de localisation en dimension infinie. C'est ce que nous faisons dans le chapitre \ref{chap-loc}. Dans le contexte symplectique, la mesure de Wiener, qui joue le r\^ole de mesure sur une orbite coadjointe dans la formule des caract\`eres, apparait naturellement lorsque l'on consid\`ere l'action hamiltonnienne d'un cercle sur une orbite. Cette id\'ee n'est pas nouvelle, le chapitre \ref{chap-loc} en propose n\'eanmoins une r\'e\'ecriture, qui me para\^it un peu originale, sous une forme adapt\'ee \`a notre contexte.  

\begin{figure}
     \begin{center}
            \includegraphics[scale=0.9]{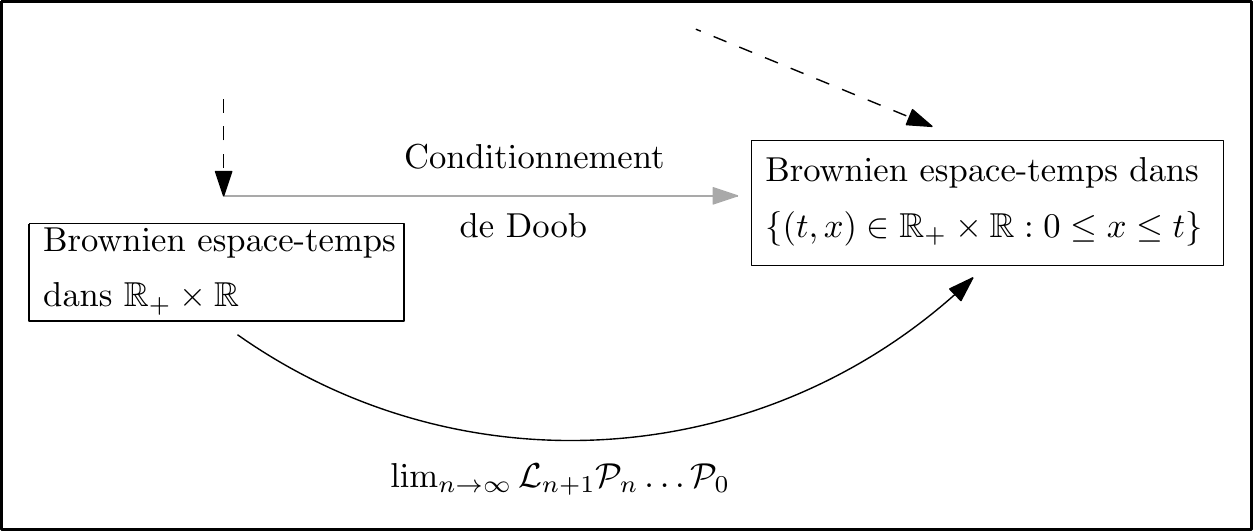}
            \caption{Transformation de Doob et transformation de Pitman - Le cas $A_1^{(1)}$}
            \label{DCFlg1P} 
                    \end{center}
 \end{figure}

 \paragraph{Un th\'eor\`eme de Pitman.}   Nous l'avons rappel\'e, Philippe Biane, Philippe Bougerol et Neil O'Connell ont \'etabli  dans \cite{bbo} un th\'eor\`eme de repr\'esentation de type Pitman pour le brownien dans la chambre de Weyl ${\bf C}_W$ associ\'ee \`a l'alg\`ebre de Lie $\mathfrak g_0$. Le mouvement brownien dans ${\bf C}_W$ s'obtient en appliquant successivement \`a un brownien dans $\mathfrak t^*$ des transformations de Pitman associ\'ees aux  r\'eflexions simples engendrant le groupe de Weyl associ\'e \`a $\mathfrak g_0$. Comme cela est indiqu\'e par les auteurs, les transformations de  Pitman jouent un r\^ole important dans les mod\`eles de chemins  de Littelmann qui sont des mod\`eles combinatoires pour les repr\'esentations d'alg\`ebres de Lie. Ces mod\`eles sont valables pour les alg\`ebres affines et il est donc naturel de chercher une repr\'esentation de type Pitman pour le brownien espace-temps apparu pr\'ec\'edemment dans la   chambre  fondamentale
${\bf \widetilde {C}}_W$ associ\'ee \`a une alg\`ebre affine.  Dans le cas compact, le nombre de transformations \`a  appliquer au processus non conditionn\'e est la longueur de l'\'el\'ement le plus long du groupe de Weyl.   Le groupe de Weyl associ\'e \`a une alg\`ebre affine est de cardinal infini et il n'existe pas de plus long \'el\'ement. On ne peut donc esp\'erer qu'un r\'esultat asymptotique. Philippe Bougerol et moi-m\^eme avons trait\'e le cas de l'alg\`ebre affine $A_1^{(1)}$ dans \cite{boubou-defo}. Il est expos\'e dans le chapitre \ref{chap-pit}. Dans ce cas le groupe de Weyl est engendr\'e par deux r\'eflexions auxquelles correspondent deux transformations de Pitman. Un r\'esultat inattendu\footnotemark\footnotetext{Nous ne l'attendions pas en tous cas !} de \cite{boubou-defo}  est que le brownien espace-temps conditionn\'e de la figure \ref{DCFlg1A} ne s'obtient pas exactement en appliquant successivement et alternativement ces transformations \`a un brownien espace-temps. Un th\'eor\`eme de repr\'esentation existe cependant qu'on obtient en   ajoutant une petite correction  \`a la suite de ces transformations successives,   correction provenant  du manque de r\'egularit\'e des trajectoires browniennes. D\'efinissons les transformations de Pitman $\mathcal P_0$ et $\mathcal P_1$ correspondant aux deux r\'eflexions ainsi que des versions modifi\'ees $\mathcal L_0$ et $\mathcal L_1$. Elles  agissent sur un chemin (espace-temps) $\eta(t)=(t,f(t)), t\in \R_+$, o\`u $f(t)\in \R$ et $f(0)=0$, de la fa\c con suivante. Pour $t\ge 0$ on a
\begin{align*}
\mathcal P_0\eta(t)&=(t, f(t)+2\inf_{s\le t} (s-f(s))), \quad  \mathcal P_1\eta(t)=(t, f(t)-2\inf_{s\le t} f(s))\\ 
 \mathcal L_0\eta(t)&=(t, f(t)+\inf_{s\le t} (s-f(s))), \quad \mathcal L_1\eta(t)=(t, f(t)-\inf_{s\le t} f(s)).
\end{align*}  
On pose $\mathcal P_{2n}=\mathcal P_0$, $\mathcal L_{2n}=\mathcal L_0$, $\mathcal P_{2n+1}=\mathcal P_1$ et $\mathcal L_{2n+1}=\mathcal L_1$, $n\ge 0$.
Nous avons montr\'e que si $\{B(t)=(t,b_t): t\ge 0\}$ est un brownien espace-temps et $\{A(t)=(t,a_t): t\ge 0\}$ le brownien espace-temps de la figure \ref{DCFlg1A} conditionn\'e \`a rester dans 
$$\{(t,x)\in \R_+\times \R: 0\le x\le t\},$$  alors les suites de processus 
$$\{\mathcal L_{n+1}\mathcal P_n \dots \mathcal P_0 B(t): t\ge 0\}\textrm{ 
et }\{\mathcal L_{n+1}\mathcal P_n \dots \mathcal P_1 B(t): t\ge 0\},\quad n\ge 0,$$
convergent en loi vers $\{A(t): t\ge 0\}$ quand $n$ tend vers l'infini. Ce r\'esultat  nous permet d'ajouter une fl\`eche au diagramme de la figure \ref{DCFlg1A} comme indiqu\'e \`a la figure \ref{DCFlg1P}. Notons que notre motivation premi\`ere pour ce travail \'etait l'obtention d'un th\'eor\`eme de repr\'esentation pour le brownien dans l'intervalle, qu'on d\'eduit du premier en appliquant une inversion temporelle au processus espace-temps conditionn\'e.  Disons enfin que comme dans le cas compact   la loi  en un temps fix\'e du brownien  conditionnellement \`a son image par nos transformations successives est une mesure de Duistermaat--Heckman, ici   la mesure image par la projection sur le dual d'un tore maximal de $\SU(2)$ de la mesure de Frenkel sur une orbite coadjointe de $\widetilde{L}(\mathfrak{su}(2))^*$.  
 
 \paragraph{Produit de fusion.} Parlons enfin d'un dernier point qui occupe une place un peu \`a part  dans le m\'emoire dont les diagrammes des figures \ref{DCA} et \ref{DCFlg1P} composent la charpente,  puisqu'il porte sur des r\'esultats  n'intervenant  pas directement dans ces diagrammes. Il concerne l'hypergoupe de la fusion obtenu \`a partir de ce qu'on appelle le produit de fusion, tel qu'il est d\'efini dans \cite{Kac}. Pour un entier $k$  fix\'e, l'hypergoupe de la fusion est l'ensemble des poids dominants entiers  de niveau $k$ d'une alg\`ebre affine $\widetilde{L}(\mathfrak g_0)$\footnote{C'est-\`a-dire l'ensemble des poids ayant une coordonn\'ee le long de $\Lambda_0$ \'egale \`a $k$.} muni d'une loi de composition donn\'ee par le produit de fusion. Il est tr\`es important de remarquer qu'ici le niveau des poids est fix\'e, contrairement \`a ce qu'on observe lorsque l'on consid\`ere un produit tensoriel de repr\'esentations d'une alg\`ebre de Lie affine.  Nous pr\'esentons dans le chapitre \ref{chap-fusion}  les r\'esultats de \cite{defo0} qui donnent une interpr\'etation probabiliste \`a ce produit. Bien qu'un peu \`a part,  comme nous l'avons dit,  au sein du m\'emoire, ces r\'esultats ne sont cependant pas sans rapport avec les questions que nous avons jusqu'ici  abord\'ees.
 
 On reconna\^itra  d'une part dans un probl\`eme de Horn multiplicatif l'esprit de l'approximation semi-classique de la mesure Duistermaat--Heckman ou encore de la m\'ethode des orbites de Kirillov. Rappelons en effet informellement la cinqui\`eme recommandation du \og User's guide \fg\,  de \cite{Kirillov}: si vous souhaitez d\'ecrire la d\'ecomposition d'un produit tensoriel  de repr\'esentations, vous devez consid\'erer la somme des orbites correspondantes que vous devez ensuite d\'ecomposer en orbites coadjointes. J'ai \'etabli dans \cite{defo0} qu'une relation similaire existait entre le produit de fusion et le produit de convolution sur un groupe de Lie compact, les orbites \'etant cette fois les orbites pour l'action par conjugaison du groupe sur lui-m\^eme.  Ce r\'esultat r\'esout une conjecture de \cite{Shaffaf}.
 
 On reconna\^itra d'autre part dans le lien \'etabli entre le produit de fusion et certaines marches dans des alc\^oves la logique de construction de cha\^ines de Markov \`a valeurs dans des chambres de Weyl \`a partir de produits tensoriels de repr\'esentations.  Je montre en effet dans \cite{defo0} que les coefficients de fusion   jouent pour une large classe de marches al\'eatoires dans des alc\^oves le m\^eme r\^ole que les coefficients de Littelwood-Richardson et leurs g\'en\'eralisations pour les marches dans une chambre de Weyl associ\'ee \`a une alg\`ebre de Lie semi-simple. Nous avons ainsi r\'epondu   positivement \`a une question pos\'ee par David Grabiner dans \cite{grabiner}  quant \`a l'existence de liens \'eventuels  entre certaines marches al\'eatoires dans des alc\^oves et   la th\'eorie des repr\'esentations.  
 
  \begin{figure}
     \begin{center}
            \includegraphics[scale=0.9]{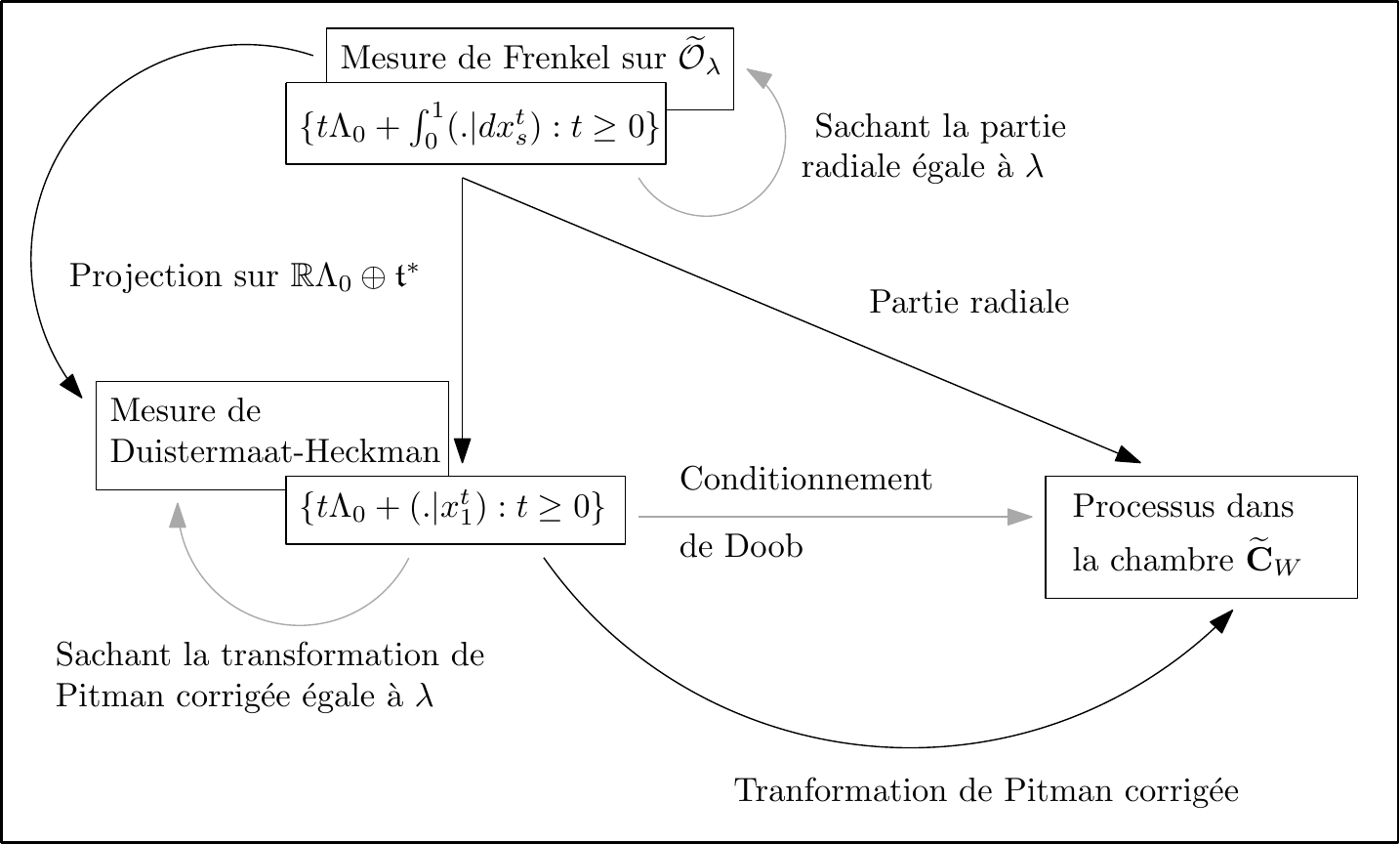}
            \caption{Sch\'ema de synth\`ese}
            \label{tout} 
                    \end{center}
 \end{figure}
 
 \paragraph{Pour finir.} Le sch\'ema de la figure \ref{tout} est une tentative de synth\`ese des r\'esultats  expos\'es dans cette introdution,  ceux concernant le produit de fusion  en \'etant cependant exclus. Pour $\lambda\in {\bf \widetilde C}_W$, $\widetilde{\mathcal O}_\lambda$ y est l'orbite de $\lambda$ sous l'action coadjointe  de  $L(G_0)$. Les fl\`eches impliquant la transformation de Pitman corrig\'ee ne peuvent pour l'instant \^etre trac\'ees  que dans le cas  d'une alg\`ebre de Lie affine  de type $A_1^{(1)}$.  Le chapitre \ref{chap-non-com} dont nous n'avons pas encore parl\'e porte sur une r\'eflexion en cours dont   l'objectif est d'\'etablir un tel diagramme dans un contexte de probabilit\'es non commutatives. Disons enfin que les  lemmes,   propositions ou th\'eor\`emes \'enonc\'es dans ce m\'emoire, \`a l'exception du th\'eor\`eme \ref{DH} de Duistermaat et de Heckman, et du theor\`eme \ref{theo-Frenkel}  dont la paternit\'e revient \`a Igor  Frenkel, sont issus de travaux que j'ai \'ecrits entre 2014 et aujourd'hui. 
 
 \newpage   
 \begin{encart}{}
\begin{enumerate}
\item Avec Philippe Bougerol, \textit{Pitman transforms and Brownian motion in the interval viewed as an affine alcove}, \`a para\^itre aux Annales Scientifiques de l'\'Ecole normale sup\'erieure  
\item \textit{Kirillov--Frenkel character formula for loop groups, radial part and Brownian sheet}, Annals of probability  (2019), 20 pp.
\item  \textit{Affine Lie algebras  and conditioned space-time Brownian motions in   affine Weyl chambers},  Probability Theory and Related Fields (2016), 17 pp.
\item \textit{Fusion coefficients and random walks in alcoves}, Ann. Inst. Henri Poincar\'e Probab. Stat. (2016), 20 pp.
\item  \textit{An interacting particle model and a Pieri-type formula for the orthogonal group\/},
Journal of Theoretical Probability (2012), 21 pp.
\item \textit{Interacting particle models and the Pieri-type formulas: the symplectic case with non equal weights\/},
Electron. Commun. Probab. 17 (2012), no. 32, 12 pp.
\item Avec  Fran\c cois Chapon, 
\textit{Quantum random walks and minors of Hermitian Brownian motion\/},
Canad. J. Math. 64 (2012), 16 pp.
\item \textit{Generalized Laguerre unitary ensembles and an interacting particles model with a wall\/},
Electron. Commun. Probab. 16 (2011), 11 pp

\item \textit{Orbit measures, random matrix theory and interlaced determinantal processes},
Ann. Inst. Henri Poincar\'e Probab. Stat. 46 (2010),  41 pp.
\item \textit{Orbit measures  and interlaced determinantal processes}, C. R. Acad. Sci. Paris, Ser. I 346 (2008).
\end{enumerate}

\end{encart}
 \chapter{Les alg\`ebres affines et leurs repr\'esentations\\ \vspace{0.5cm}
\small{\it o\`u l'on dit deux ou trois choses que l'on doit savoir d'$\mathcal{L}(\mathfrak{g})$}}\label{chap-Algebre-affine}

 Les alg\`ebres de Kac--Moody  forment une classe d'alg\`ebres de Lie sur $\C$ que l'on peut d\'efinir par   g\'en\'erateurs  et relations \`a partir d'une matrice de Cartan g\'en\'eralis\'ee. Cette classe contient, outre les alg\`ebres de Lie semi-simples complexes, des alg\`ebres de Lie de dimension infinie, et parmi elles,  des alg\`ebres de Kac--Moody de type affine non tordu. Ces derni\`eres pr\'esentent l'immense avantage de poss\'eder une r\'ealisation simple \`a partir d'une alg\`ebre   de dimension finie.   Certains de nos r\'esultats sont valables pour une alg\`ebre affine quelconque. Cependant, ce sont celles de type non tordu que nous avons choisi de consid\'erer exclusivement dans ce m\'emoire, pour la simplicit\'e de leur pr\'esentation, mais surtout parce que leur r\'ealisation \`a  partir d'une alg\`ebre de lacets sera dans les chapitres \ref{chap-loc} et \ref{chap-drap} plus ajust\'ee  \`a nos questionnements. Les d\'efinitions et propri\'et\'es   \'enonc\'ees sont toutes issues du livre de Victor G. Kac \cite{Kac}.
\section{Alg\`ebres de Kac--Moody  affines}  
Nous d\'efinissons dans cette section les alg\`ebres de lacets $\widehat{\mathcal L}(\mathfrak g)$  \`a valeurs dans $\mathfrak g$, avec  extension centrale et d\'erivation, o\`u $\mathfrak g$ est une alg\`ebre de Lie complexe simple. Remarquons que la d\'erivation sera fondamentale dans notre travail.  En effet,  les caract\`eres de certaines repr\'esentations de $\widehat{\mathcal L}(\mathfrak g)$ y jouent un r\^ole essentiel  et ils ne sont pas d\'efinis si on restreint ces repr\'esentations \`a l'alg\`ebre sans sa d\'erivation.  Par ailleurs, d'un point de vue plus analytique, la d\'erivation provient d'une action hamiltonienne du cercle $S^1$ qui     apparait de mani\`ere fondamentale dans le chapitre \ref{chap-loc} portant sur les th\'eor\`emes de localisation en dimension infinie. Comme nous le verrons cependant, m\^eme si la d\'erivation est cruciale pour nos constructions, nos r\'esultats ne portent \`a la fin que sur l'alg\`ebre sans d\'erivation\footnotemark\footnotetext{C'est-\`a-dire  l'alg\`ebre affine dite d\'eriv\'ee\dots}.

 \subsection{Extensions d'alg\`ebres de lacets}
Consid\'erons une alg\`ebre de Lie simple  complexe $\mathfrak{g}$ de dimension finie munie d'un crochet  de Lie $[\cdot\, ,\cdot]_\mathfrak{g}$ et   l'alg\`ebre de lacets 
$$\mathcal L(\mathfrak{g})\coloneqq\mathcal L\otimes_\C \mathfrak g,$$ 
o\`u $\mathcal L$ est l'alg\`ebre $\C[z^{-1},z]$ des polyn\^omes de Laurent en $z$.  Autrement dit $\mathcal L(\mathfrak{g})$ est l'ensemble des polyn\^omes de Laurent \`a coefficients dans  $\mathfrak{g}$. Poser $z=e^{2i\pi \theta}$, $\theta\in [0,1]$, identifie $\mathcal L(\mathfrak g)$ \`a une alg\`ebre de lacets. C'est une alg\`ebre de Lie pour le crochet $[\cdot  ,\cdot]$ d\'efini par 
$$[P\otimes x ,Q\otimes y]=PQ\otimes [x,y]_\mathfrak{g},\quad  x,y\in \mathfrak{g}, \quad P,Q\in \mathcal L.$$ On munit $\mathfrak{g}$ d'une forme  bilin\'eaire      sym\'etrique $(\cdot\vert \cdot)$  invariante et non d\'eg\'en\'er\'ee puis on consid\`ere l'extension centrale avec d\'erivation   de l'alg\`ebre de lacets 
$$\widehat{\mathcal L}(\mathfrak g)=\mathcal L(\mathfrak{g})\oplus \C c\oplus \C d,$$
munie d'un crochet, qu'on note toujours $[\cdot ,\cdot]$,  d\'efini par
\begin{align}\label{crochet-lie}
[z^{n_1}\otimes  x_1 \oplus & a_1c\oplus b_1d,z^{n_2}\otimes x_2 \oplus a_2c\oplus b_2d] \nonumber \\
&= z^{n_1+n_2}\otimes [x_1,x_2]_{\mathfrak{g}}+b_1n_2z^{n_2}\otimes x_2-b_2n_1z^{n_1}\otimes x_1\nonumber\\
&\quad \quad\quad\quad\quad\quad\quad\quad\quad\quad\quad\quad\quad\quad\quad\quad+n_1\delta_{n_1+n_2}(x_1\vert x_2)c,
\end{align}
 pour $ x_i\in \mathfrak{g}$, $a_i,b_i\in \C$, $n_i\in \Z$, $i\in\{1;2\}$. L'alg\`ebre de Lie de dimension infinie ainsi obtenue est une r\'ealisation d'alg\`ebre de Kac--Moody de type affine non tordu. Plus pr\'ecis\'ement si $\mathfrak{g}$ est une alg\`ebre de Lie de rang $n$ de type $X_n$ alors $\widehat{\mathcal{L}}(\mathfrak{g})$ est une alg\`ebre affine de type $X_n^{(1)}$.

\subsection{Racines et coracines de $\widehat{\mathcal{L}}(\mathfrak{g})$}\label{section-roots-coroots-A}
On choisit une sous-alg\`ebre de Cartan $\mathfrak{h}$ de $\mathfrak{g}$ suppos\'ee de rang $n$. On consid\`ere la d\'ecomposition en sous-espaces radiciels
\begin{align}\label{radiC}
\mathfrak{g}=\mathfrak{h}\oplus\bigoplus_{\alpha\in \Phi} \mathfrak{g}_\alpha,
\end{align}
o\`u  $\Phi$ l'ensemble des racines de $\mathfrak{g}$ et $\mathfrak{g}_\alpha=\{x\in \mathfrak{g}: \forall h\in \mathfrak{h}, [h,x]=\alpha(h)x\}$.
On choisit un sous-ensemble de racines simples  $\Pi=\{\alpha_1,\dots,\alpha_n\}$
et on note  $\Phi_+$ l'ensemble des racines positives.    On consid\`ere   une sous-alg\`ebre ab\'elienne maximale de $\widehat{\mathcal{L}}(\mathfrak{g})$ $$\widehat{\mathfrak{h}}=\mathfrak{h}\oplus \C c\oplus \C d$$ et on   \'etend les racines de $\Phi$ \`a $\widehat{\mathfrak h}$ en posant $\alpha(c)=\alpha(d)=0$, $\alpha\in \Phi$. Alors l'ensemble des racines de l'alg\`ebre de Lie $\widehat{\mathcal L}(\mathfrak{g})$ est $$\widehat{\Phi}=\{\alpha+k\delta:  \alpha\in \Phi, k\in \Z\}\cup\{k\delta:  k\in \Z^*\},$$
o\`u $\delta$ est la racine nulle, qui est nulle sur $\mathfrak{h}\oplus \C c$ et v\'erifie $\delta(d)=1$. On a en fait la  d\'ecomposition en sous-espaces  radiciels par rapport \`a $\widehat{\mathfrak{h}}$   
\begin{align}\label{radiA}
\widehat{\mathcal{L}}(\mathfrak{g})=\widehat{\mathfrak{h}}\oplus \bigoplus_{\alpha\in \Phi,\, k\in \Z}\widehat{\mathcal{L}}(\mathfrak{g})_{\alpha+k\delta}\oplus \bigoplus_{k\in\Z^*}\widehat{\mathcal{L}}(\mathfrak{g})_{ k\delta},
\end{align}
avec  
$\widehat{\mathcal{L}}(\mathfrak{g})_{\alpha+k\delta}=z^k\otimes  \mathfrak{g}_\alpha,\, \alpha\in \Phi, k\in \Z,$ et $\widehat{\mathcal{L}}(\mathfrak{g})_{k\delta}=z^k\otimes\mathfrak{h}, \, k\in \Z^*. $
 Pour $\alpha\in \Pi$, les sous-espaces radiciels $\mathfrak{g}_\alpha$ sont de dimension $1$ et  on choisit un triplet  $(h_\alpha,e_\alpha,f_\alpha)\in \mathfrak{h}\times \mathfrak{g}_\alpha\times\mathfrak{g}_{-\alpha}$ tel que 
$$[e_\alpha ,f_\alpha]=h_\alpha,\quad [h_\alpha ,e_\alpha]=2e_\alpha,\quad [h_\alpha ,f_\alpha]=-2f_\alpha.$$ On identifie $\mathfrak h$ et $\mathfrak h^*$ via 
\begin{align}\label{iso}
\nu :h\in \mathfrak h\to (h\vert \cdot)\in \mathfrak h^*,
\end{align} et on munit $\mathfrak h^*$ du produit scalaire induit par cet isomorphisme.  On choisit une normalisation telle que $(\theta\vert \theta)=2$, o\`u $\theta$ est la racine la plus haute de $\mathfrak g$. La forme $(\cdot\vert\cdot)$ est ainsi la forme bilin\'eaire invariante dite standard normalis\'ee.  On pose $\alpha_0=\delta-\theta$. On choisit un triplet $(h_{\alpha_0},e_{\alpha_0},f_{\alpha_0})\in \widehat{\mathfrak{h}}\times \widehat{\mathcal L}(\mathfrak{g})_{\alpha_0}\times \widehat{\mathcal L}(\mathfrak{g})_{-\alpha_0}$ tel que 
$$[e_{\alpha_0} ,f_{\alpha_0}]=h_{\alpha_0},\quad [h_{\alpha_0} ,e_{\alpha_0}]=2e_{\alpha_0},\quad [h_{\alpha_0} ,f_{\alpha_0}]=-2f_{\alpha_0}.$$ Alors $h_{\alpha_0}=c- \theta^\vee$   o\`u $\theta^\vee$ est la coracine la plus haute, et    $$\widehat\Pi=\{\alpha_0=\delta-\theta,\alpha_1,\dots,\alpha_n\}$$ est un ensemble de racines simples de $\widehat{\mathcal L}(\mathfrak g)$ et 
 $$\widehat\Pi^\vee=\{\alpha^\vee_0=h_{\alpha_0}, \alpha_1^\vee= h_{\alpha_1},\dots,\alpha_n^\vee= h_{\alpha_n}\}$$ est l'ensemble des coracines simples. 
 
 Remarquons que les d\'ecompositions en sous-espaces radiciels (\ref{radiC}) et (\ref{radiA}) permettent d'identifier respectivement   $\mathfrak h^*$ et $\widehat{\mathfrak{h}}^*$ aux sous-espaces des formes de  $\mathfrak g^*$ et de $\widehat{\mathcal{L}}(\mathfrak g)^*$ nulles sur les sous-espaces  radiciels autres que la sous-alg\`ebre de Cartan.
 
\subsection{Groupe et chambre de Weyl affine}
On consid\`ere l'\'el\'ement $\Lambda_0$ de $\widehat{\mathfrak{h}}^*$ d\'efini par 
$$ \Lambda_0(\alpha_0^\vee)=1,\quad \Lambda_0(d)=0,\quad \Lambda_0(\alpha_i^\vee)=0,\quad  i\in\{1,\dots,n\}.$$
Alors $\{\Lambda_0,\alpha_0,\dots,\alpha_n\}$ est une base de $\widehat{\mathfrak{h}}^*$. On consid\`ere la forme bilin\'eaire  $(\cdot \vert\cdot)$ que l'on restreint \`a $\mathfrak{h}$ et \'etend \`a $\widehat{\mathfrak{h}}$ en posant pour $x\in \mathfrak{h}$,
$$(c\vert x)=(d\vert x)=(c\vert c)=(d\vert d)=0, \quad \textrm{et } \quad (c\vert d)=1.$$ 
L'isomorphisme lin\'eaire
$\nu:\, h\in \widehat{\mathfrak{h}}\mapsto (h\vert \cdot) \in\widehat{\mathfrak{h}}^* $, identifie $\widehat{\mathfrak{h}}$ et $\widehat{\mathfrak{h}}^*$. La forme  bilin\'eaire $(\cdot\vert \cdot)$ sur $\widehat{\mathfrak{h}}^*$ induite par l'identification v\'erifie
\begin{align*}
&(\delta\vert \alpha_i)=0,\quad  i=0,\dots,n,\quad (\delta\vert\delta)=0,\quad (\delta\vert\Lambda_0)=1.
\end{align*}    
Le groupe de Weyl   $\widehat{W}$ de $\widehat{\mathcal L}(\mathfrak{g})$ est un sous-groupe  de $\mbox{GL}(\widehat{\mathfrak{h}}^*)$ engendr\'e par  les r\'eflexions   fondamentales $s_\alpha$, $\alpha\in \widehat\Pi$, d\'efinies par $$s_\alpha(\beta)=\beta-  \beta(\alpha^\vee) \alpha,\ \quad \beta\in \widehat{\mathfrak{h}}^*.$$   
 La forme $(\cdot \vert \cdot)$ est $\widehat W$-invariante. Notons $Q^\vee$  le r\'eseau des coracines de $\mathfrak g$ et $W$ le groupe de Weyl engendr\'e par les r\'eflexions $s_\alpha$ pour $\alpha\in \Pi$.  Le groupe de Weyl affine $\widehat W$  est le produit semi-direct  $W\ltimes  {\Gamma}$  o\`u  $\Gamma$ est  le groupe des  transformations $t_{\gamma}$, $\gamma\in \nu( {Q}^\vee)$, d\'efinies par 
$$t_\gamma(\lambda)=\lambda+\lambda(c)\gamma-((\lambda\vert \gamma)+\frac{1}{2}(\gamma\vert\gamma)\lambda(c))\delta, \quad \lambda\in \widehat{\mathfrak{h}}^*.$$ 
On  note $\widehat{\mathfrak{h}}_\R= \bigoplus_{i=1}^n\R\alpha_i^\vee\oplus \R d\oplus \R c$ et $\widehat{\mathfrak{h}}^*_\R$ l'ensemble des formes lin\'eaires \`a valeurs r\'eelles sur $\widehat{\mathfrak{h}}_\R$, i.e. $\widehat{\mathfrak{h}}^*_\R=\bigoplus_{i=1}^n\R\alpha_i \oplus \R\Lambda_0\oplus \R\delta$. L'ensemble 
$$\widehat {  C}_W=\{\lambda\in\widehat{ \mathfrak{h}}_\R^*:   \lambda(\alpha^\vee)\ge 0, \, \alpha\in \widehat\Pi\}$$ est appel\' e chambre fondamentale ou chambre de Weyl. C'est un domaine  fondamental pour l'action de $\widehat W$ sur un c\^one de $\widehat{\mathfrak{h}}_\R^*$ appel\'e c\^one de Tits. 
\section{Repr\'esentations, repr\'esentations int\'egrables}\label{section-rep-aff}
  On note    $\widehat P$ (resp. $\widehat P_+$) l'ensemble des poids entiers   (resp. dominants)   d\'efini par  
$$\widehat P=\{\lambda\in \widehat{\mathfrak{h}}^*:   \lambda(\alpha^\vee_i )\in \Z, \, i=0,\dots, n\},$$
 $$(\textrm{resp. } \widehat P_+=\{\lambda\in\widehat P:  \lambda(\alpha^\vee_i )\ge 0, \, i=0,\dots,n\}).$$
Le niveau d'un poids   $\lambda\in \widehat P$, est l'entier $(\delta\vert\lambda)$, i.e. la coordonn\'ee le long de $\Lambda_0$. Pour $k\in \N$, on note     $\widehat P^k$ (resp. $\widehat P^k_+)$ l'ensemble des poids entiers (resp. dominants) de niveau $k$, i.e.
$$\widehat P^k=\{\lambda\in \widehat P: (\delta\vert\lambda)=k\}.$$
$$(\textrm{resp. }\widehat P^k_+=\{\lambda\in \widehat P_+: (\delta\vert\lambda)=k\}.)$$
Un $\widehat{\mathcal{L}}(\mathfrak{g})$-module $V$ est  dit $\widehat{\mathfrak{h}}$-diagonalisable s'il admet une d\'ecomposition en espaces de poids   $V=\bigoplus_{\lambda\in \widehat{\mathfrak{h}}^*}V_\lambda$ o\`u  $V_\lambda$ est d\'efini par 
$$V_\lambda=\{v\in V: \forall h\in \widehat{\mathfrak{h}},\, h.v=\lambda(h)v\}.$$ 
La cat\'egorie mono\"idale  $\mathcal O$ est d\'efinie comme l'ensemble des $\widehat{\mathcal{L}}(\mathfrak{g})$-modules $V$  $\widehat{\mathfrak{h}}$-diagonalisables dont les espaces de poids sont de dimension finie et tels qu'il existe un nombre fini d'\'el\'ements    $\lambda_1,\dots,\lambda_s\in\widehat{\mathfrak{h}}^*$ tels que 
\begin{align}\label{O} 
P(V)\subset \cup_{i=1}^s\{\mu\in \widehat{\mathfrak{h}}^*: \lambda_i-\mu\in \N \widehat{\Phi}_+\},
\end{align}
o\`u $P(V)=\{\lambda\in \widehat{\mathfrak{h}}^*: V_\lambda\ne \{0\}\}$. Les modules de cette cat\'egorie, bien que  de dimension infinie, partagent de nombreuses propri\'et\'es avec les modules d'alg\`ebres de Lie semi-simples complexes de dimension finie. Tout d'abord, le fait que les espaces de poids d'un module de $\mathcal O$ soient de dimension finie et la condition (\ref{O}) assurent  qu'on puisse d\'efinir pour  un module $V$ de $\mathcal O$ un    caract\`ere formel $\mbox{ch}(V)$     en posant
$$\mbox{ch}(V)=\sum_{\mu\in P(V)}\dim(V_\mu)e^{\mu}.$$  
Nous verrons dans le chapitre \ref{chap-BETaff} une proc\'edure de construction de cha\^ines de Markov \`a partir de d\'ecompositions en composantes  irr\'eductibles de produits tensoriels de repr\'esentations. La cat\'egorie $\mathcal O$ poss\`ede une sous-cat\'egorie  dont les modules sont compl\`etement r\'eductibles. C'est la cat\'egorie des modules int\'egrables, sur lesquels  les g\'en\'erateurs de Chevalley $e_i,f_i$, $i\in\{0,\dots,n\}$ agissent de mani\`ere localement nilpotente. Elle est not\'ee $\mathcal O_{int}$. Les modules irr\'eductibles de  $\mathcal O_{int}$ sont les modules ayant un plus haut poids  dans $\widehat P_+$. Pour $\lambda\in \widehat P_+$ on  dispose pour un module $V(\lambda)$ de plus haut poids $\lambda$,  d'une formule de 
  Weyl des caract\`eres     qui s'\'ecrit  
\begin{align}\label{Weyl}
\mbox{ch}(V(\lambda))=\frac{\sum_{w\in\widehat W}\det(w)e^{w(\lambda+\hat\rho)-\hat\rho}}{\prod_{\alpha\in\widehat \Phi_+}(1-e^{-\alpha})^{\mbox{mult}(\alpha)}},
\end{align}
o\`u $\mbox{mult}(\alpha)$ est la dimension de l'espace radiciel $\widehat{\mathcal{L}}(\mathfrak{g})_\alpha$ et   $\hat\rho$ est un vecteur de  $\widehat{\mathfrak{h}}^*$ appel\'e vecteur de Weyl, choisi tel que $\hat\rho(\alpha_i^\vee)=1$, pour tout $i\in\{0,\dots,n\}$. Notons  $\langle\cdot,\cdot\rangle$  l'appariement dual canonique entre  $\widehat{\mathfrak h}$ et $\widehat{\mathfrak h}^*$.
Pour un module  $V$ de plus haut poids $\lambda\in P_+$ de $\mathcal{O}_{int}$ et $h\in \widehat{\mathfrak{h}}$ la s\'erie $$\sum_{\mu\in P(V)}\dim(V_\mu)e^{\langle\mu,h\rangle}$$   est absolument convergente si et seulement si $\mbox{Re}\langle \delta, h\rangle>0$.   Ce r\'esultat de convergence est tr\`es important pour nous, puisque les caract\`eres \'evalu\'es en des vecteurs bien choisis nous permettront de d\'efinir les processus de Markov introduits au chapitre \ref{chap-BETaff}.

\section{L'alg\`ebre de Kac--Moody affine $\widehat{\mathcal{L}}(\mathfrak{sl}_2(\C))$.}\label{section-affA11}
Nous pr\'ecisons  dans cette section  ce que sont les objets d\'efinis ci-dessus dans le cas  o\`u   $\mathfrak{g}$ est l'alg\`ebre de Lie $\mathfrak{sl}_2(\C)$ d\'efinie par $$\mathfrak{sl}_2(\C)=\{M\in \mathcal{M}_2(\C): \mbox{tr}(M)=0\}.$$
Les  matrices
$$h=\begin{pmatrix} 
1& 0 \\
0& -1
\end{pmatrix}, \, e=\begin{pmatrix} 
0& 1 \\
0& 0
\end{pmatrix}, \, f=\begin{pmatrix} 
0& 0 \\
1& 0
\end{pmatrix},$$
forment une base de $\mathfrak{sl}_2(\C)$ v\'erifiant les relations de commutation
$$[h,e]=2e,\, [h,f]=-2f,\, [e,f]=h.$$
La sous-alg\`ebre des matrices diagonales $\C h$ est une  sous-alg\`ebre de Cartan et  
$\mathfrak{sl}_2(\C)$ poss\`ede une seule racine positive $\alpha_1$  d\'efinie par $\alpha_1(h)=2$. On munit $\mathfrak{sl}_2(\C)$ d'une forme sesquilin\'eaire $(\cdot \vert \cdot)$ d\'efinie par $(x\vert y)=\mbox{tr}(xy^*)$, $x,y\in \mathfrak{sl}_2(\C)$.
On consid\`ere l'alg\`ebre de Lie $$\widehat{\mathcal{L}}(\mathfrak{sl}_2(\C))=\mathcal L( \mathfrak{sl}_2(\C))\oplus \C c\oplus \C d,$$
le crochet de Lie \'etant d\'efini par (\ref{crochet-lie}). C'est une alg\`ebre de Lie de type $A_1^{(1)}$.
\paragraph{Racines et sous-espaces radiciels.} On pose $\widehat{\mathfrak{h}}=\C h\oplus \C c\oplus \C d$ et on   \'etend $\alpha_1$ \`a $\widehat{\mathfrak{h}}$  en posant $\alpha_1(c)=\alpha_1(d)=0$. On d\'efinit la racine nulle $\delta \in \widehat{\mathfrak{h}}^*$ en posant   $$\delta(c)=\delta(h)=0\, \textrm{ et } \, \delta(d)=1.$$ On pose $\alpha_0=\delta-\alpha_1$.  On v\'erifie facilement que l'on a pour $n\in \Z$,
\begin{alignat*}{5}
&[h,z^n\otimes e]&=&2z^n\otimes e,\,\,
&[d,z^n\otimes e]&=&nz^n\otimes e,\\
&[h,z^n\otimes f]&=&-2z^n\otimes f,\,\, 
 &[d,z^n\otimes f]&=&nz^n\otimes f,\\
 &[d,z^n\otimes h]&=&nz^n\otimes h, & &
 \end{alignat*}
 ce qui montre que les sous-espaces radiciels sont
 \begin{align*}
&\widehat{\mathcal{L}}(\mathfrak{sl}_2(\C))_{\alpha_1+n\delta}=\C z^n\otimes e, \, \,
\widehat{\mathcal{L}}(\mathfrak{sl}_2(\C))_{-\alpha_1+n\delta}=\C z^n\otimes f, \,\, n\in \Z,\\
&\widehat{\mathcal{L}}(\mathfrak{sl}_2(\C))_{n\delta}=\C z^n\otimes h,\, \, n\in \Z^*,
\end{align*}
et les g\'en\'erateurs de Chevalley 
$$\alpha_1^\vee=h_{\alpha_1}=1\otimes h, \quad e_{\alpha_1}=1\otimes e,\quad  f_{\alpha_1}=1\otimes f,$$
$$\alpha_0^\vee=h_{\alpha_0}=c-h_{\alpha_1}, \quad e_{\alpha_0}=z\otimes f, \quad f_{\alpha_0}=z^{-1}\otimes e.$$
 \paragraph{Ensemble des poids et ensemble des poids dominants.}
 Soit le poids fondamental $\Lambda_0$ d\'efini par $\Lambda_0(h)=\Lambda_0(d)=0$ et $\Lambda_0(c)=1$. 
 L'ensemble des poids de $\widehat{\mathcal L}(\mathfrak{sl}_2(\C))$ est $$\widehat P=\{x\Lambda_0+y\frac{\alpha_1}{2}+z\delta :  x,y\in \Z, z\in \C\},$$ et celui des poids dominants   $$\widehat P_+=\{x\Lambda_0+y\frac{\alpha_1}{2}+z\delta :  x,y\in \Z, z\in \C, 0\le y\le x\}.$$
 
\begin{figure}[!t]
\centering
\includegraphics[scale=0.7]{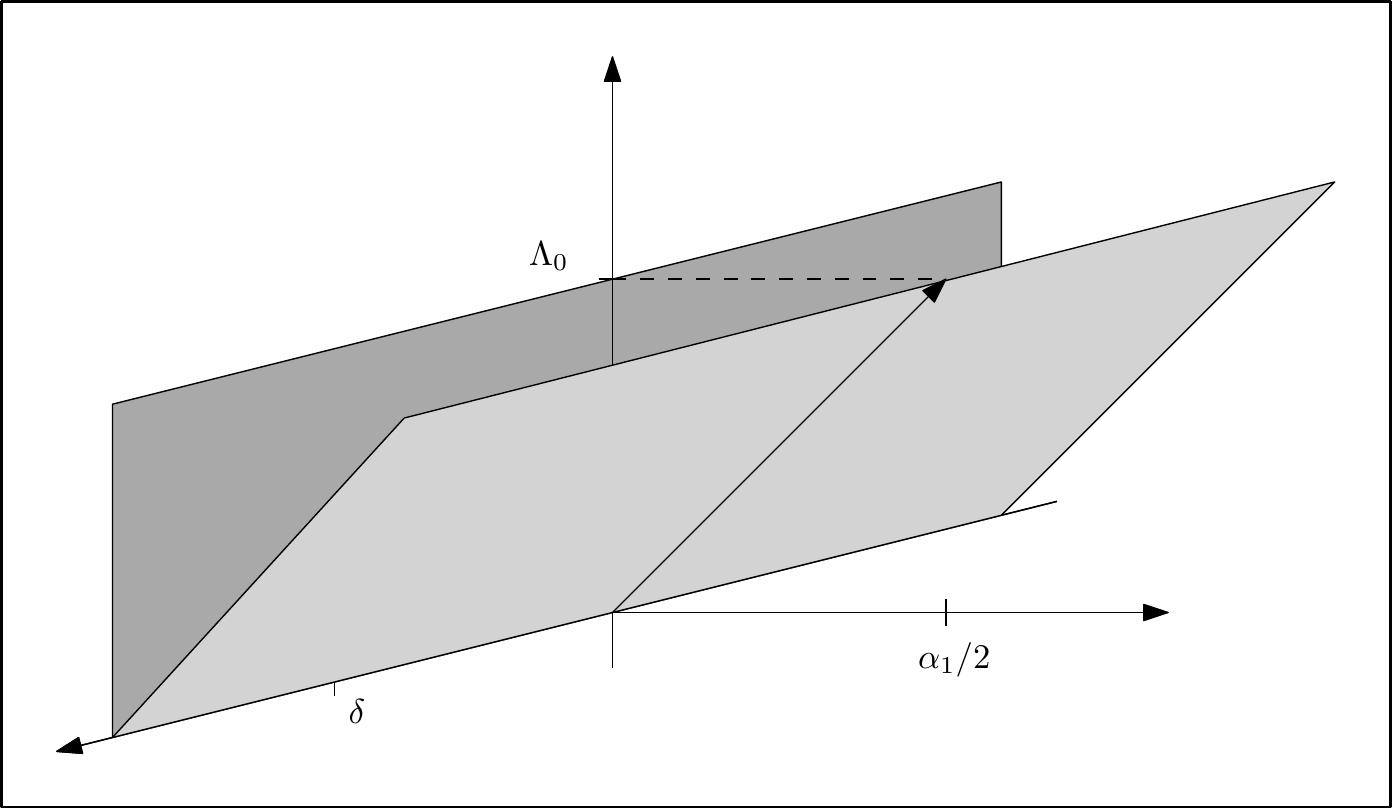}
\caption{Chambre de Weyl pour l'alg\`ebre $\widehat{\mathcal{L}}(\mathfrak{sl}_2(\C))$} 
\label{CAA}
\end{figure}

 \paragraph{Groupe de Weyl et domaine fondamental.} Le groupe de Weyl associ\'e \`a $\mathfrak{sl}_2(\C)$ est le groupe des transformations engendr\'e par la r\'eflexion $s_{\alpha_1}$ d\'efinie par 
 $$s_{\alpha_1}(x\Lambda_0+y\alpha_1/2)=x\Lambda_0-y\alpha_1/2, \quad x,y\in \R,$$ 
 Le r\'eseau des coracines est $Q^\vee=\Z\alpha^\vee_1$ et ici $\nu(Q^\vee)=\Z\alpha_1$. Ainsi  le groupe de Weyl affine $\widehat W$  est le produit semi-direct  $\{s_{\alpha_1},\mbox{Id}_{\mathfrak h^*}\} \ltimes  {\Gamma}$  o\`u  $\Gamma$ est  le groupe des  transformations $t_{k\alpha_1}$, $k\in \Z$, d\'efinies par 
$$t_{k\alpha_1}(x\Lambda_0+y\alpha_1/2)=x\Lambda_0+(y+2kx)\alpha_1/2-(ky+k^2x)\delta, \quad x,y\in \R, k\in \Z.$$
La chambre de Weyl est  ici
$$\widehat C_W=\{x\Lambda_0+y\alpha_1/2+z\delta: 0\le y\le x,\,  x,y\in\R^*_+, z\in \R\}.$$
Elle est repr\'esent\'ee en figure \ref{CAA}. C'est le domaine compris entre les deux hyperplans gris. 
Remarquons que si l'on consid\`ere l'action de $\widehat W$   sur l'espace   quotient $\Lambda_0+\R\alpha_1+\R\delta\mod \delta$, alors le groupe de Weyl  s'identifie au groupe de transformations de $\R$  engendr\'e par les  r\'eflexions par rapport aux entiers. Vu dans l'espace quotient\'e par $\R\delta$,
$$\{\Lambda_0+x\alpha_1/2+y\delta: 0\le x\le 1,\,  y\in \R\}$$ s'identifie  a l'intervalle $[0,1]$ qui est un  domaine fondamental de ce groupe de transformations non-lin\'eaires. Il est repr\'esent\'e sur la figure \ref{CAA} par la ligne en pointill\'es.

 \paragraph{Formule des caract\`eres de Weyl.} Le vecteur $\hat \rho=2\Lambda_0+\frac{1}{2}\alpha_1$ est un vecteur de Weyl et la formule des caract\`eres de Weyl   s'\'ecrit pour $\lambda=x\Lambda_0+y\alpha_1/2$ et $z=ad+b\alpha^\vee_1$, avec $a>0$,
 \begin{align}\label{char-a1}
 \mbox{ch}(V(\lambda))(z)&=\frac{\sum_{k\in\Z}\sinh(b(y+1)+2bk(x+2))e^{-a(k(y+1)+k^2(x+2))}}{\sum_{k\in\Z}\sinh(b+4bk)e^{-a(k+2k^2)}}\\
 &=\frac{\sum_{k\in\Z}\sinh((y+1)(b+ak))e^{-(x+2)(2kb+k^2a)}}{\sum_{k\in \Z}\sinh(b+ak)e^{-2(k2b+k^2a)}}.\end{align}

 \chapter{Browniens espace-temps  et repr\'esentations d'alg\`ebres affines\\ \vspace{0.5cm}
\small{\it  o\`u l'on met de  petits poids sur les bas poids}}\label{chap-BETaff}

Les formules de Clebsch--Gordan donnent en particulier la d\'ecomposition en composantes irr\'eductibles du produit tensoriel d'une repr\'esentation de dimension finie de $\SU(2)$ et de la repr\'esentation standard de ce groupe.  D\`es 1975, Bernard Roynette   \cite{Roynette} remarque que ces formules fournissent les probabilit\'es de transition d'un Bessel  $3$ discret.  Au  d\'ebut des ann\'ees 90, Philippe Biane   fait appara\^itre ce Bessel  dans un contexte non commutatif \cite{biane1}, donnant ainsi un nouvel \'eclairage \`a cette co\"incidence. Outre la joie des correspondances, que  nous apporte un tel rapprochement ?  Le processus de Bessel $3$ discret est un processus de Markov \`a valeurs dans l'ensemble des entiers naturels, c'est-\`a-dire le r\'eseau des poids dominants de $\SU(2)$.  Sa loi est  la transformation de Doob de celle de la marche simple sur $\Z$ tu\'ee en $-1$.  Les pas de la marche simple  sont eux distribu\'es selon la mesure de probabilit\'e uniforme  sur l'ensemble des poids de la repr\'esentation standard de $\SU(2)$. Si on remplace ce groupe  sp\'ecial unitaire par un groupe de Lie compact connexe semi-simple quelconque et la repr\'esentation standard par une repr\'esentation complexe   de dimension finie du groupe, on peut de la m\^eme mani\`ere associer \`a celle-ci deux processus, une marche al\'eatoire et un processus de Markov. Les pas de la marche sont distribu\'es selon la mesure de probabilit\'e uniforme sur les poids de la repr\'esentation compt\'es avec leur multiplicit\'e et  le noyau de transition du processus de Markov s'exprime \`a partir de la dimension des composantes isotypiques d'un produit tensoriel de repr\'esentations. Les repr\'esentations irr\'eductibles du groupe sont   d\'etermin\'ees par leur plus haut poids et on obtient ainsi une marche \`a valeurs dans  le r\'eseau des poids du groupe et un processus de Markov \`a valeurs  dans l'ensemble de ses poids dominants. Ces processus  fournissent d'une part une approximation d'un brownien  dans le dual d'une alg\`ebre de Cartan de l'alg\`ebre de Lie du groupe, d'autre part une approximation  du brownien  conditionn\'e au sens de Doob \`a vivre dans une chambre de Weyl    associ\'ee au syst\`eme de racines du groupe. La combinatoire des repr\'esentations de ce groupe \'eclaire alors certaines propri\'et\'es de ces browniens. C'est une approche possible par exemple du th\'eor\`eme de Pitman et de ses g\'en\'eralisations, qui se trouvent ainsi li\'es au mod\`ele combinatoire des chemins de Littelmann \cite{bbo}. Pour nous, une telle approche dans un cadre affine aura fourni les processus de la base du diagramme repr\'esent\'e \`a la figure  \ref{tout} de l'introduction.

Dans \cite{defo1} et \cite{defo2} nous consid\'erons  en effet  des processus \`a valeurs dans le r\'eseau des poids   d'une alg\`ebre de Lie affine ou dans celui de ses poids dominants construits   selon un canevas analogue \`a celui que nous venons de d\'ecrire bri\`evement. Ces processus sont pour la premi\`ere fois introduits dans \cite{LLP}. Le cas affine diff\`ere   du cas compact sur plusieurs points essentiels. Tout d'abord,  les repr\'esentations irr\'eductibles consid\'er\'ees sont de dimension infinie et  il n'existe   pas de mesure de probabilit\'e uniforme sur l'ensemble des poids d'une telle repr\'esentation. Nous sommes donc amen\'es \`a affecter  \`a ces poids  un poids de  Boltzmann.  
D'autre part,  les processus discrets introduits  pr\'esentent une coordonn\'ee temporelle le long d'un poids fondamental tout comme en cons\'equence  les processus limites qui  apparaissent par ailleurs  \`a une \'echelle de temps de l'ordre de $n$, et une \'echelle d'espace de l'ordre de $\frac{1}{n}$ au lieu du facteur $\frac{1}{\sqrt{n}}$ habituel. 

  Nous d\'etaillons dans la section \ref{section-MA-su(2)} le cas de la marche simple sur $\Z$ avec un drift \'eventuel conditionn\'ee \`a rester positive. Nous expliquons dans quelle mesure cette marche est li\'ee aux  repr\'esentations complexes de dimension finie  de $\mathfrak{sl}_2(\C)$.    La vocation de cette section dans laquelle est d\'ecrite une situation tr\`es simple est de d\'egager des \'el\'ements de structure qui perdurent dans des situations plus complexes. Nous exposons   le cas affine dans la section \ref{section-affine}. Le th\'eor\`eme \ref{TCL-cond} est le principal r\'esultat du chapitre. Il fait apparaitre des processus limites qui sont pour notre travail des objets essentiels que nous reverrons aux chapitres  \ref{chap-drap} et  \ref{chap-pit}. Leur mise en \'evidence  aura \'et\'e la premi\`ere \'etape de notre exploration des liens entre mouvement brownien et alg\`ebre de Kac--Moody affine.   
 \section{Marches al\'eatoires et repr\'esentations de  $\mathfrak{sl}_2(\C)$}\label{section-MA-su(2)}
    \paragraph{Repr\'esentations de $\mathfrak{sl}_2(\C)$.}  D\'eterminer les repr\'esentations complexes de dimension finie du groupe compact $\SU(2)$ revient \`a d\'eterminer celles  de   l'alg\`ebre de Lie complexifi\'ee  de  l'alg\`ebre de Lie $\mathfrak{su}(2)$ de $\SU(2)$. Cette alg\`ebre est l'alg\`ebre $\mathfrak{sl}_2(\C)$  d\'efinie par  $$\mathfrak{sl}_2(\C)=\{M\in \mathcal{M}_2(\C): \mbox{tr}(M)=0\},$$
et les  matrices
$$h=\begin{pmatrix} 
1& 0 \\
0& -1
\end{pmatrix}, \, e=\begin{pmatrix} 
0& 1 \\
0& 0
\end{pmatrix} \,\textrm{ et } f=\begin{pmatrix} 
0& 0 \\
1& 0
\end{pmatrix}$$
forment une base de $\mathfrak{sl}_2(\C)$.
  Une repr\'esentation   $V$ de dimension finie de $\mathfrak{sl}_2(\C)$ se d\'ecompose en espaces de poids 
  $$V=\bigoplus_{z\in \Z} V_z,$$
  o\`u $V_z=\{v\in V : h.v=z v\}$. Si $V$ est irr\'eductible il existe un unique entier $x\in \N$ tel que $V_x\ne \{0\}$ et $V_{x+2}=\{0\}$. On dit que $x$ est le plus haut poids de $V$ et que $V$ est la repr\'esentation de plus haut poids $x$. On la note $V(x)$. L'ensemble  des entiers naturels est l'ensemble des poids dominants de $\mathfrak{sl}_2(\C)$. Pour $x\in \N$, il existe une repr\'esentation irr\'eductible $V(x)$ de $\mathfrak{sl}_2(\C)$ de plus haut poids $x$. L'action de $\mathfrak{sl}_2(\C)$ sur $V(x)$  est donn\'ee par 
  $$\left\{
    \begin{array}{l}
      h\cdot v_i=(x-2i)v_i \\
      e\cdot v_i=(x-i+1)v_{i-1}\\
      f\cdot v_i=(i+1)v_{i+1}, 
          \end{array}
\right.
$$
  pour $i\in\{0,\dots,x\},$
  o\`u $v_0\in V(x)_x$, $v_{-1}=0$ et $v_i=(1/i!)f^{i} \cdot v_0$.  L'ensemble des poids de $V(x)$ est $\{x
  -2i: i\in \{0,\dots,x\}\}$, et chaque poids appara\^it avec une multiplicit\'e $1$.
  Le caract\`ere $\mbox{ch}(V)$ d'une repr\'esentation de dimension finie $V$ de $\mathfrak{sl}_2(\C)$ est d\'efini par 
  $$\mbox{ch}(V)(q)=\sum_{z\in \Z}\dim(V_z)q^z, \quad q\in \R_+^*.$$
  En particulier pour $x\in \N$, 
    \begin{align*}
    \mbox{ch}(V(x))(q)&=\sum_{i=0}^{x} q^{x-2i}
    =\frac{q^{x+1}-q^{-(x+1)}}{q-q^{-1}},
    \end{align*} 
    la derni\`ere identit\'e n'\'etant rien d'autre que la formule des caract\`eres de Weyl.

\begin{figure}[!t]
\centering
\includegraphics[scale=0.8]{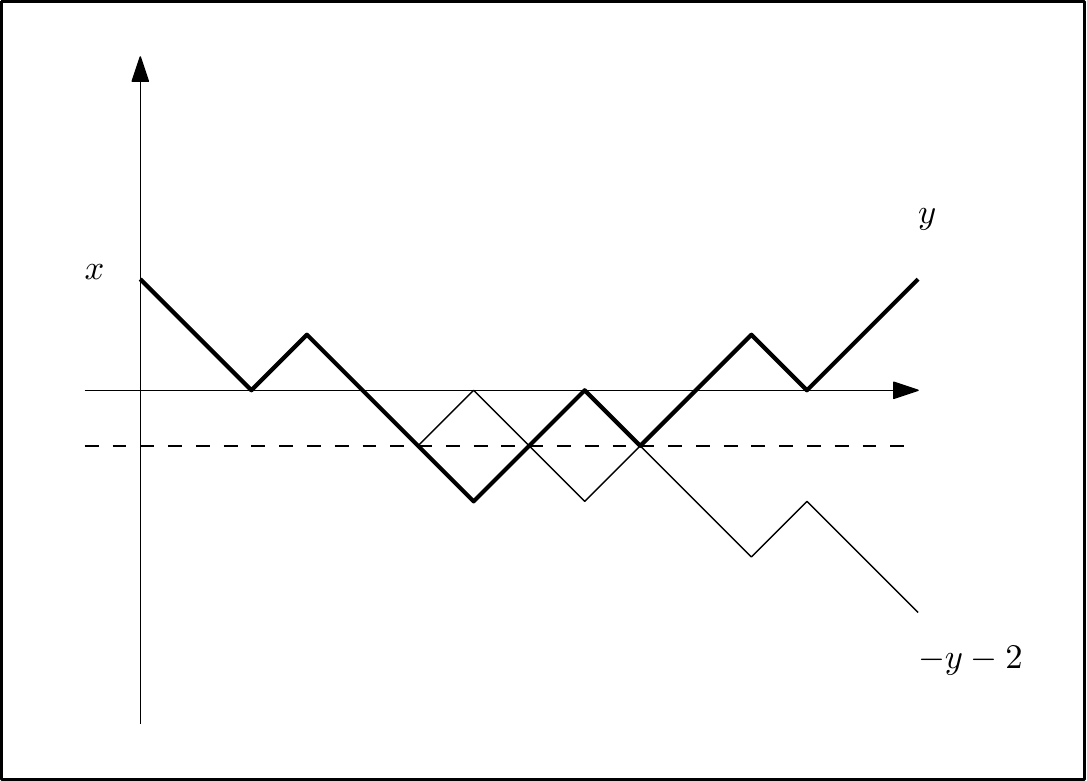}
\caption{Principe de r\'eflexion} 
\label{refl}
\end{figure}

   \paragraph{Marches al\'eatoires,  principe de r\'eflexion et transformation de Doob.}     Fixons un r\'eel strictement positif      $q$. On consid\`ere   une marche al\'eatoire simple  $\{X(n), n\ge 0\}$  avec drift  de  noyau  de transition
   $$K(x,y)=\frac{q^{y-x}}{q+q^{-1}} \mathds{1}_{\vert x-y\vert=1}, \quad x,y\in \Z.$$ 
Ici nous consid\'erons le cas d'une marche avec drift car il montre davantage de similitudes que le cas centr\'e avec ce que nous ferons dans le cas affine. On introduit le temps d'atteinte de $-1$ qu'on note $T$. Une \'etude rapide de suite r\'ecurrente lin\'eaire d'ordre $2$ montre qu'\`a un coefficient multiplicatif pr\`es, il existe pour la marche tu\'ee $\{X(n\wedge T), n\ge 0\}$ une unique fonction harmonique sur $\N\cup\{-1\}$, positive sur $\N$ et nulle en $-1$. C'est la fonction $x\in \N\cup\{-1\} \mapsto q^{-x}s_x(q)$, o\`u $s_x(q)$ est la fonction de Schur  d\'efinie par 
  \begin{align}\label{schur}
  s_x(q)=\frac{q^{x+1}-q^{-(x+1)}}{q-q^{-1}}, \quad x\in \N\cup\{-1\}.
  \end{align}
   On d\'efinit une transformation de Doob de la marche tu\'ee via cette fonction harmonique. C'est le processus de Markov \`a valeurs dans  $\N$ de probabilit\'e  de transition 
   \begin{align}\label{MASP}
   L(x,y)=\frac{q^{-y}s_y(q)}{q^{-x}s_x(q)}K(x,y), \quad x,y\in \N.
   \end{align}
  Un principe de r\'eflexion pour $q=1$ illustr\'e en figure \ref{refl}, suivi d'un th\'eor\`eme de Girsanov, montre que pour $x,y\ge -1$, on a 
  $$\P_x(X(n)=y,\, T\ge n+1)=\P_x(X(n)=y)-q^{2y+2}\P_x(X(n)=-y-2),$$
  ou de mani\`ere \'equivalente,
  \begin{align}\label{refl-su}
  L^n(x,y)&=\frac{q^{-y}s_y(q)}{q^{-x}s_x(q)}\big(K^n(x,y)-q^{2y+2}K^n(x,-y-2)\big)
  \end{align}
 Nous allons voir que ce dernier processus fait partie d'une classe de processus obtenus \`a partir de produits tensoriels de repr\'esentations, en donnant   une proc\'edure un peu syst\'ematique de construction de marches al\'eatoires et de processus de Markov, \`a valeurs respectivement dans le r\'eseau des poids $\Z$ de $\mathfrak{sl}_2(\C)$ et dans celui de ses poids dominants $\N$, proc\'edure que nous appliquerons ensuite dans le cadre des alg\`ebres affines.

   \paragraph{Marches al\'eatoires,  repr\'esentations et transformation de Doob.}    
 Choisissons un r\'eel strictement positif $q$ et un entier naturel non nul $\omega$  et consid\'erons    la repr\'esentation irr\'eductible $V(\omega)$ de $\mathfrak{sl}_2(\C)$ de plus haut poids $\omega$. On associe alors \`a la repr\'esentation $V(\omega)$ une marche al\'eatoire  $\{X_\omega(n),n\ge 0\}$ sur $\Z$, dont les accroissements sont distribu\'es selon une loi $\mu_\omega$ d\'efinie par 
\begin{align}\label{acc-su} 
\mu_\omega(y)=\frac{q^{y}}{ \mbox{ch}(V(\omega))(q)},\quad y\in \{-\omega,-\omega+2,\dots,\omega+2,\omega\}.
\end{align}
En particulier pour $\omega=1$ et $q=1$, $\{X_\omega(n),n\ge 0\}$ est la marche simple sur $\Z$. Pour $x\in \N$, on consid\`ere la d\'ecomposition  en composantes irr\'eductibles
$$V(x)\otimes V(\omega) \simeq \bigoplus_{z\in \N} M_{x,\omega}^z\, V(z)$$ 
ou  de mani\`ere \'equivalente, 
\begin{align}\label{prod-car-su} \mbox{ch}(V(x))\times  \mbox{ch}(V(\omega))=\sum_{z\in \N} M_{x,\omega}^z\, \mbox{ch}(V(z)),\end{align}
o\`u $M_{x,\omega}^z$ est la multiplicit\'e du module $V(z)$ dans le produit tensoriel $V(x)\otimes V(\omega)$. On d\'efinit  un noyau markovien $L_\omega$ sur $\N$ en posant 
\begin{align} \label{MarkovSU2} 
L_\omega(x,z)=\frac{\mbox{ch}(V(z))(q)}{\mbox{ch}(V(x))(q)\mbox{ch}(V(\omega))(q)}\, M_{x,\omega}^z, \quad x, z\in \N.
\end{align}
Autrement dit
\begin{align} 
L_\omega(x,z)=\frac{s_z(q)}{s_x(q)s_\omega(q)}\, M_{x,\omega}^z, \quad  x,z\in \N.
\end{align} 
  Quand   $\omega=1$, les coefficients  $M_{x,\omega}^z$, $z\in \N$,  valent $1$ si   $z\in \{x-1,x+1\}$ et    $0$ sinon.  La d\'ecomposition (\ref{prod-car-su}) \'evalu\'ee en $q$ dit dans ce cas que la fonction $x\mapsto q^{-x}s_x(q)$ est harmonique   pour la marche $X$ tu\'ee en $-1$  et le noyau $L_\omega$ est alors la probabilit\'e de transition $L$ d\'efinie plus haut.  
 
 \paragraph{Principe de r\'eflexion et r\`egle de Brauer et de Klimyk.}  Pour une alg\`ebre de Lie semi-simple complexe,   la r\`egle de Brauer et de Klimyk \cite{Stembridge} fournit une expression des multiplicit\'es des composantes irr\'eductibles dans un produit tensoriel de repr\'esentations montrant une somme altern\'ee de dimensions d'espaces de poids portant sur le groupe de Weyl associ\'e au syst\`eme de racines de l'alg\`ebre de Lie. Dans un cadre probabiliste ou combinatoire de telles formules sont connues sous le nom de formules de  Lindstr\"om-Gessel-Viennot ou Gessel--Zeilberger qui,  sous le prisme de la th\'eorie des repr\'esentations,  se comprennent  donc   comme des formules de multiplicit\'e\footnote{On pourrait aussi citer Karlin--McGregor, pour le cas des chemins en temps continu.} .   La r\`egle de Brauer et de Klimyk   pour l'alg\`ebre de Lie  $\mathfrak{gl}_n(\C)$ par exemple est \'equivalente \`a  une formule de Gessel et de Zeilberger portant sur le nombre de trajectoires possibles d'un point \`a un autre  d'une marche simple sur $\Z^n$  dont les coordonn\'ees  sont assujetties \`a ne pas se croiser.  Plus g\'en\'eralement, on peut faire la m\^eme remarque pour toute marche al\'eatoire dont les pas sont  \`a valeurs dans une r\'eunion d'ensembles de poids de repr\'esentations minuscules d'une alg\`ebre de Lie semi-simple complexe fix\'ee.      Nous d\'etaillons ici les r\`egles de Brauer et de Klimyk et leur lien avec le principe de r\'eflexion dans le  cas de $\mathfrak{sl}_2(\C)$. Dans la suite nous notons 
 $$s_x(q)=(q^{x+1}-q^{-(x+1)})/(q-q^{-1}),$$ pour tout $x\in \Z$, $q\in \R^*$.  Pour une repr\'esentation  $V$ de dimension finie de $\mathfrak{sl}_2(\C)$,  si, pour $y\in \Z$, $m_V(y)$ est la multiplicit\'e du poids $y$ dans $V$, on a imm\'ediatement,  pour $z\in \N$,   l'identit\'e 
 $$s_{z}(q)\mbox{ch}(V)(q)=\sum_{y\in \Z}m_V(y)s_{z+y}(q).$$
 Cette identit\'e est connue sous le nom de  r\`egle de Brauer--Klimyk \cite{Stembridge}. Elle d\'ecoule de la formule de Weyl. Si $M_{V\otimes V(x)}^z$ est la multiplicit\'e de $V(z)$ dans $V\otimes V(x)$ , la   r\`egle de Brauer--Klimyk implique que     
 \begin{align}\label{Klimyk-su}
 M_{V\otimes V_x}^z=m_V(z-x)-m_V(-(z+1)-(x+1)).
 \end{align} En prenant $V=V(\omega)^{\otimes n}$, on obtient 
 \begin{align*}
 L_\omega^n(x,z)=\frac{s_z(q)}{s_x(q)s_1(q)^n}(m_{V(\omega)^{\otimes n}}(z-x)- m_{V(\omega)^{\otimes n}}(-(z+1)-(x+1))),
 \end{align*} qui est  l'identit\'e (\ref{refl-su})  lorsque  $\omega=1$. 
 
 \paragraph{Convergence vers le brownien.} Soit $\gamma$ un r\'eel positif. La fonction $\phi_\gamma$ d\'efinie sur $\R_+$ par  
$$\phi_\gamma(x)= 1-e^{-2\gamma x}, \quad x\in \R_+,$$
est harmonique pour le brownien standard sur $\R$ avec drift $\gamma$ tu\'e en $0$,  positive sur $\R_+$ et nulle en $0$.  On consid\`ere le transform\'e   de Doob $\{A_t^\gamma: t\ge 0\}$  via la fonction  
$\phi_\gamma$ du brownien avec drift tu\'e en $0$. Ses  densit\'es   de transition  $\{q_t: t\ge 0 \}$ s'expriment en fonction du noyau de la chaleur  $\{p_t: t\ge 0\}$ sur $\R$. On a 
 \begin{align*}
q_t(x,y)&=\frac{\phi_\gamma(y)}{\phi_\gamma(x)}(p_t(x,y)-e^{-2\gamma x}p_t(-x,y)),\quad x,y,t>0.
\end{align*}  
Consid\'erons $\{Z^{(n)}(k): k\ge 0\}$ la transformation de Doob issue de $0$ de la marche simple tu\'ee en $-1$ dont le noyau est donn\'e par (\ref{MASP}) avec  $q=e^{ {\gamma}/{\sqrt{n}}}$. On a alors la convergence en loi,
\begin{align} \label{conv-su}
\{\frac{{Z^{(n)}}\big([nt]\big)}{\sqrt{n}},t\ge 0\}\underset{n\to\infty}\longrightarrow \{A_t^\gamma,t\ge 0\},
\end{align}
qui s'obtient   ais\'ement en utilisant  la convergence de la marche simple vers le brownien et en consid\'erant $L^n(\sqrt n x,\sqrt n y)$, o\`u  $L^n$   est donn\'ee en (\ref{refl-su}) avec  $q=e^{ \gamma/\sqrt{n}}$. Une telle convergence reste valable dans le cas d'un processus de Markov de noyau de transition d\'efini par  (\ref{MarkovSU2}).

\section{Marches al\'eatoires et repr\'esentations d'alg\`ebres affines}\label{section-affine}
Nous rempla\c cons maintenant l'alg\`ebre de Lie $\mathfrak{sl}_2(\C)$   par une alg\`ebre de Kac--Moody affine. Contrairement \`a celles des alg\`ebres de Lie semi-simples complexes, les repr\'esentations de plus haut poids d'une  telle  alg\`ebre   ne sont pas de dimension finie. Dans la situation d\'ecrite dans la section pr\'ec\'edente, affecter la valeur $1$ au param\`etre $q$ revient \`a munir l'ensemble des poids d'une repr\'esentation d'une mesure de probabilit\'e uniforme. Cela  n'a pas de sens dans un cadre affine, les repr\'esentations que nous consid\'erons ayant une dimension infinie. Cependant nous pouvons mettre en \oe uvre la proc\'edure d\'ecrite ci-dessus avec  un drift bien choisi pour lequel les quantit\'es consid\'er\'ees pour $\mathfrak{sl}_2(\C)$ gardent un sens. Nous obtenons alors une marche sur le r\'eseau des poids d'une alg\`ebre affine et un processus de Markov sur le r\'eseau de ses poids dominants.

On se place dans le cadre  du  chapitre \ref{chap-Algebre-affine}. On consid\`ere une alg\`ebre de Lie simple complexe   et l'alg\`ebre affine associ\'ee.    Nous l'avons vu plus haut, l'alg\`ebre de Lie $\mathfrak{sl}_2(\C)$  n'intervient  pas  directement dans la description des marches et des processus de Markov. Seules comptent  les d\'ecompositions en espaces de poids des repr\'esentations ou leur d\'ecomposition en composantes irr\'eductibles. Nous n'avons besoin de consid\'erer que la sous-alg\`ebre de Cartan 
$$\widehat{\mathfrak{h}}=\mathfrak{h}\oplus \C c \oplus \C d,$$ 
ainsi que sa partie r\'eelle 
$$\widehat{\mathfrak{h}}_\R=\mathfrak{h}_\R\oplus \R c \oplus \R d,$$
o\`u $\mathfrak{h}_\R=\mbox{vect}_\R\{\alpha_i^\vee :  i\in \{1,\dots, n\}\}$. Nous consid\'erons de m\^eme   le dual de la sous-alg\`ebre de Cartan  $\widehat{\mathfrak{h}}^*$ et sa partie r\'eelle  
$$\widehat{\mathfrak{h}}^*_\R=\mathfrak{h}^*_\R\oplus \R \Lambda_0 \oplus \R \delta.$$
o\`u $\mathfrak{h}^*_\R=\mbox{vect}_\R\{\alpha_i : i\in \{1,\dots, n\}\}$.  Nous consid\'erons   le r\'eseau des poids $\widehat P$ qui jouera le r\^ole de $\Z$ et   celui des   poids dominants $\widehat P_+$ qui jouera le r\^ole de $\N$.  Notons qu'une forme sur $\widehat{\mathfrak{h}}$ est d\'etermin\'ee par sa restriction \`a $\widehat{\mathfrak{h}}_\R$.   Dans ce chapitre, les   poids sont vus comme des \'el\'ements de $\widehat{\mathfrak h}_\R^*$ qu'on peut, quand c'est n\'ecessaire pour les d\'efinitions, \'etendre \`a $\widehat{\mathfrak h}$ par $\C$-lin\'earit\'e.  Nous l'avons rappel\'e, les repr\'esentations irr\'eductibles  de la cat\'egorie $\mathcal O_{int}$  de l'alg\`ebre affine sont les repr\'esentations de plus haut poids dans $\widehat{P}_+$.  Pour un poids dominant  $\lambda\in \widehat P_+$, on consid\`ere le caract\`ere  de la repr\'esentation $V(\lambda)$ de plus haut poids   $\lambda$. Il est d\'efini par une s\'erie formelle 
\begin{align*} 
 \sum_{\beta}\mbox{dim} (V(\lambda)_\beta)e^\beta , 
\end{align*}
o\`u  $V(\lambda)_\beta$ est l'espace de poids  $\beta$ de $V(\lambda)$. Nous avons rappel\'e que la s\'erie $$\sum_\beta  \mbox{dim} (V(\lambda)_\beta)e^{\langle \beta,h\rangle}$$ converge  absolument pour tout  $h$ dans   $Y$ avec $Y=\{h\in\widehat{ \mathfrak h}: \mbox{Re}(\langle\delta,h\rangle)>0\}$.  Pour $h\in Y$, on pose   $$\mbox{ch}_\lambda(h)=\sum_\beta  \mbox{dim} (V(\lambda)_\beta)e^{\langle \beta,h\rangle}.$$ 
\paragraph{Marche al\'eatoire sur le r\'eseau des poids.}  Fixons un poids dominant non nul $\omega\in \widehat P_+$.
 Pour $h\in Y_\R$ avec $Y_\R=Y\cap\widehat{\mathfrak{h}}_\R$,  on d\'efinit une mesure de probabilit\'e $\mu_\omega$ sur $\widehat P$ en posant
\begin{align}\label{paffine}
\mu_\omega(\beta)=\frac{  \mbox{dim} (V(\omega)_\beta)}{\mbox{ch}_\omega(h)}e^{\langle \beta,h\rangle}, \quad \beta \in \widehat P.
\end{align}
On consid\`ere une marche al\'eatoire issue de z\'ero $\{X_\omega(n),n\ge 0\}$ \`a valeurs dans $\widehat P$ dont les accroissements sont distribu\'es selon la loi $\mu_\omega$.\paragraph{Cha\^ine de Markov sur le r\'eseau des poids dominants.}
Pour un poids dominant  $\lambda$, on consid\`ere la d\'ecomposition suivante
\begin{align}\label{decomp-aff}
 V(\omega)\otimes V(\lambda) =\sum_{\beta\in \widehat P_+} M_{\omega,\lambda}^\beta V(\beta),
\end{align}
o\`u $M_{\omega,\lambda}^\beta$ est la multiplicit\'e du module de plus haut poids    $\beta$ dans la repr\'esentation $V(\omega)\otimes V(\lambda)$. Cela permet de  d\'efinir comme nous l'avons fait dans la section pr\'ec\'edente une probabilit\'e de transition $Q_\omega$ sur $\widehat P_+$ en posant, pour $\beta$ et $\lambda$ deux poids dominants, 
\begin{align}\label{qaffine}
Q_\omega(\lambda,\beta)=\frac{\mbox{ch}_\beta(h)}{\mbox{ch}_\lambda(h) \mbox{ch}_\omega(h)}M_{\omega,\lambda}^\beta.
\end{align} 
On consid\`ere une cha\^ine de Markov issue de z\'ero $\{Z_\omega(k): k\ge 0\}$ \`a valeurs dans $\widehat{P}_+$ de probabilit\'e de transition $Q_\omega$. 
La question qui nous a int\'eress\'ee est celle de la convergence, apr\`es mise \`a l'\'echelle, des suites de processus  $\{X_\omega([nt]), t \ge 0\}$ et $\{Z_\omega([nt]):t\ge 0\}$, $n\ge 0$.   Il est important de noter que les poids d'une repr\'esentation irr\'eductible   sont tous de m\^eme niveau.  Ainsi  $X_\omega(k)$ et $Z_\omega(k)$ sont, pour $k\in \N$, des poids de niveau $k\langle \delta,\omega\rangle$.  Autrement dit  la coordonn\'ee le long de $\Lambda_0$ de ces processus, qui est donc \`a penser comme une coordonn\'ee temporelle, vaut $k\langle \delta,\omega\rangle$. Une renormalisation standard en $1/\sqrt{n}$ fait    apparaitre  une coordonn\'ee   explosant le long de $\Lambda_0$.  Il semble donc plus naturel, afin de conserver la coordonn\'ee le long de $\Lambda_0$, de chercher une suite  $\{h_n,n\ge 1\}$ \`a valeurs dans  $Y_\R$, telle que  si $X^{(n)}_\omega$ et $Z_\omega^{(n)}$ sont respectivement la marche et la cha\^ine  d\'efinies par  (\ref{paffine}) et (\ref{qaffine}), avec $h=h_n$,  alors les suites de processus $$\{\frac{1}{n}X_\omega^{(n)}([nt]): t\ge 0\}\, \textrm{ et  }\,\{\frac{1}{n}Z_\omega^{(n)}([nt]): t\ge 0\}, \quad n\ge 1,$$  convergent.   En fait, pour une suite $\{h_n: n\ge 1\}$ telle que  les $h_n\sim \frac{h}{n}$, avec $h\in Y_\R$,  les projections  de ces processus sur $\R\Lambda_0\oplus \mathfrak{h}^*_\R$ convergent tandis que  leurs coordonn\'ees le long de $\delta$ divergent. C'est aux parties convergentes que nous nous int\'eressons.  La projection sur $\mathfrak{h}^*_\R$  est ce que nous appellons la partie spatiale. La coordonn\'ee le long de $\Lambda_0$ est la coordonn\'ee temporelle. Enfin, la coordonn\'ee le long de $\delta$ peut \^etre pens\'ee comme une version discr\`ete  de l'\'energie du brownien, qui appara\^it dans l'expression symbolique 
$$\int_{\mathcal C([0,T],  \mathfrak{h}_\R^*)}f(x)e^{-\frac{1}{2}\vert\vert \dot x\vert\vert^2}\, d x$$
de l'int\'egrale d'une fonction contre la mesure de Wiener sur l'ensemble des chemins \`a valeurs dans  $\mathfrak{h}_\R^*$ (voir chapitre \ref{chap-loc}, section \ref{subsectionorbite}). Pour un \'el\'ement $\xi \in \widehat{\mathfrak{h}}_\R^*$, nous notons $\tilde\xi$ la projection de $\xi$ sur $\R\Lambda_0\oplus\mathfrak{h}^*_\R$.   

\paragraph{Principe de r\'eflexion, r\`egle de de Brauer et de Klimyk.} Comme dans le cas de $\mathfrak{sl}_2(\C)$, le processus de Markov sur le r\'eseau des poids dominants  satisfait un principe de r\'eflexion. Il r\'esulte comme dans le cas semi-simple d'une  r\`egle  de Brauer et de Klimyk.  Pour   une repr\'esentation $V$ de $\mathcal O_{int}$ et $\lambda\in \widehat{P}_+$, cette r\`egle donne, pour $\beta\in \widehat P$, une expression de la multiplicit\'e $M_{V\otimes V(\lambda)}(\beta)$ de $V(\beta)$ dans $V\otimes V(\lambda)$ comparable \`a l'expression (\ref{Klimyk-su}), la  somme  portant dans ce cas sur le groupe de Weyl affine $\widehat{W}$ \cite{Stembridge}.  On a ainsi la formule  suivante qui relie les multiplicit\'es  $M_{V \otimes V(\lambda)}(\beta)$ aux dimensions des espaces de poids de la repr\'esentation~$V$
\begin{align} \label{BrauerKlimyk} 
M_{V\otimes V(\lambda) }(\beta)=\sum_{w\in \widehat W}\det(w)m_{V}(w(\beta+\hat\rho)-(\lambda+\hat\rho)),
\end{align} 
o\`u  $m_V(\beta)$ est la dimension de  l'espace de poids $V_\beta$.
Nous l'avons dit, c'est la projection $\widetilde{Z}_\omega$ de la cha\^ine   $Z_\omega$ sur  $\R\Lambda_0\oplus \mathfrak{h}_\R^*$ qui nous importe.   Si  $\lambda_1$ et $\lambda_2$ sont deux poids dominants tels que  $\lambda_1= \lambda_2\mod \delta$, alors les modules irr\'eductibles  $V(\lambda_1)$ et $V(\lambda_2)$ sont isomorphes. Le processus $(\widetilde{Z}_\omega(k),k\ge 0)$ reste donc markovien.  On note $\widetilde Q_\omega$ sa probabilit\'e de transition. La formule (\ref{BrauerKlimyk})  appliqu\'ee \`a la repr\'esentation $V=V(\omega)^{\otimes n}$ permet d'obtenir pour   $\widetilde Q_\omega$ l'expression suivante.

\begin{prop}\label{reflectiondiscrete}  Soient $\beta$ et $\lambda$ deux poids dominants de  $ \R\Lambda_0\oplus  {\mathfrak{h}}^*_\R$, et $n$  un entier naturel. La probabilit\'e de transition $\widetilde Q_\omega$ v\'erifie  
\begin{align*}
\widetilde Q_\omega^n(\lambda,\beta)& =\frac{\mbox{ch}_{\beta}(h)e^{-\langle\beta,h\rangle}}{\mbox{ch}_{\lambda}(h)e^{-\langle\lambda,h\rangle} }  \sum_{w\in \widehat W}\det(w)e^{\langle w(\lambda+\hat\rho)-(\lambda+\hat\rho),h\rangle}\widetilde P_\omega^n(\widetilde{w(\lambda+\hat\rho)}-\widetilde{\hat\rho},\beta)
\end{align*}
o\`u $\widetilde P_\omega$ est la probabilit\'e de transition de la projection $\widetilde X_\omega $ de $X_\omega$ sur $\R\Lambda_0\oplus \mathfrak{h}_\R^*$.
\end{prop}

\section{Brownien espace-temps et alg\`ebres affines}
Ici $ \mathfrak h$ est muni de la forme bilin\'eaire standard  normalis\'ee $(\cdot \vert\cdot)$ que l'on \'etend \`a $\widehat{\mathfrak{h}}$ comme dans le chapitre \ref{chap-Algebre-affine}. On rappelle que l'isomorphisme lin\'eaire $\nu$ est d\'efini par $$\nu: h\in \widehat{\mathfrak{h}}\to (h\vert \cdot)\in  \widehat{\mathfrak{h}}^*,$$ qu'il identifie $\widehat{\mathfrak h}$ et $\widehat{\mathfrak h}^*$ et qu'on note    toujours   $(\cdot\vert \cdot)$ la forme bilin\'eaire sur $\widehat{\mathfrak h}^*$ induite par l'identification. Dans \cite{defo2} nous avons consid\'er\'e une suite de drifts  $\{h_n: n\ge 1\}$ avec $h_n=\frac{1}{n} \nu^{-1}(\hat\rho)$, o\`u $\hat\rho$ est un vecteur de Weyl. Les processus limites poss\`edent alors eux-m\^eme un drift $\nu^{-1}(\hat \rho)$. Ici nous pr\'esentons les r\'esultats pour une suite $\{h_n:n\ge 1\}$ avec un drift nul sur $\mathfrak{h}$, afin d'obtenir des processus limites eux-m\^emes sans drift spatial, c'est-\`a-dire des browniens espace-temps standard. Les cas avec ou sans drift ne pr\'esentent pas de diff\'erences essentielles pour nous. Ils sont li\'es par un th\'eor\`eme de Girsanov.   

Soit $\{b_t: t\ge 0\}$ un mouvement brownien standard sur $\mathfrak h_\R^*$ et un processus $$\{B_t=t\Lambda_0+b_t: t\ge 0\}.$$  Comme pour les marches introduites dans la section pr\'ec\'edente, la coordonn\'ee temporelle est le long de $\Lambda_0$ tandis que la partie spatiale est la projection sur $\mathfrak{h}_\R^*$. Le processus $B$   peut ainsi \^etre vu comme un brownien espace-temps.   La chambre de Weyl affine,  vue dans l'espace quotient $  \widehat{\mathfrak{h}}_\R^* /\R\delta$ identifi\'e \`a $\R\Lambda_0\oplus  {\mathfrak{h}}_\R^* $, est 
$$\widetilde C_W=\{\lambda\in \R\Lambda_0\oplus  {\mathfrak{h}}_\R^*  : \lambda(\alpha^\vee_i)\ge 0,  i\in\{0,\dots,n\}\},$$
o\`u $n$ est le rang de $\mathfrak{h}$.
Si   $\{\lambda_k:k\ge 0\}$ est une suite de plus hauts poids telle que  $\lambda_k\sim k\lambda,$ avec $\lambda\in  \widetilde C_W$, et $\{h_k:k\ge 0\}$ une suite \`a valeurs dans $Y_\R$ telle que $h_k\sim h/k$, avec $h\in Y_\R$, alors le num\'erateur de $\mbox{ch}(V(\lambda_k))(h_k)$ dans la formule de Weyl converge vers 
 \begin{align} \label{num-char}
 \sum_{w\in \widehat W}\det(w)e^{\langle w\lambda,h\rangle}.\end{align}
  Pour $\lambda\in  \widetilde C_W$ et $h\in Y$, on pose   
   $$\widehat\varphi_h(\lambda)=\frac{1}{\pi( {h}/{\delta(h)})}\sum_{w\in \widehat W}\det(w)e^{\langle w\lambda,h\rangle}.$$
o\`u $\pi(\cdot)=\prod_{k=1}^n(e^{i\pi\alpha_k(\cdot)}-e^{-i\pi\alpha_k(\cdot)})$.  On peut montrer que  cette expression est bien d\'efinie par continuit\'e pour tout  vecteur $h$ dans $Y$.  
La proposition suivante permet de consid\'erer un conditionnement de Doob du brownien espace-temps tu\'e sur la fronti\`ere de   $\widetilde C_W$.  
\begin{prop} La fonction $\widehat \varphi_d$ est   harmonique pour le brownien espace-temps tu\'e sur  la fronti\`ere de $\widetilde C_W$. Elle est strictement positive \`a l'int\'erieur de  $\widetilde C_W$ et nulle sur sa fronti\`ere. 
\end{prop}
\begin{defn}\label{Brownien-esp-tps-cond}  Le processus $$\{A_t=t\Lambda_0+a_t:t\ge 0\}$$ est un processus de Markov issu de $0$, \`a valeurs dans l'int\'erieur de $\widetilde C_W$ pass\'e le temps initial. C'est le transform\'e de Doob issu de $0$ via $\widehat\varphi_d$ du brownien espace-temps $\{B_t: t\ge 0\}$ tu\'e sur la fronti\`ere de $\widetilde C_W$.  
\end{defn}
Il est notable que  le brownien tu\'e satisfait un principe de r\'eflexion analogue \`a celui \'enonc\'e dans  la proposition \ref{reflectiondiscrete}. C'est un point important pour les preuves. Notons $T$ le temps d'atteinte de la fronti\`ere de $\widetilde C_W $.
\begin{prop} Pour $r,t> 0$ et $\lambda$ un \'el\'ement de niveau $r$ dans l'int\'erieur de $\widetilde C_W$,  on a 
 $$\P(b_{(t+r)\wedge T}\in dy\, \vert\,  B_r=\lambda)=\sum_{w\in \widehat W}\det(w)e^{\langle w(\lambda)-\lambda ,d\rangle}p_t(\widetilde{w(\lambda)}-r\Lambda_0,y)dy,$$
o\`u $p_t$ est le noyau de la chaleur sur $\mathfrak h_\R^*$.
\end{prop}
Pour \'enoncer le th\'eor\`eme, nous prenons $\omega=\Lambda_0$. Un autre plus haut poids  non nul conviendrait, les parties temporelles et spatiales des processus limites ne d\'ependant du plus haut poids choisi qu'\`a travers deux constantes multiplicatives.  Pour un entier $n\ge 1$, consid\'erons la marche $\{X^{(n)}(k):k\ge 0\}$ issue de $0$, dont les accroissements sont distribu\'es selon la loi $\mu_{\Lambda_0}$ d\'efinie par (\ref{paffine}) avec $h=\frac{1}{n}d$, et la cha\^ine de Markov $\{Z^{(n)}(k):k\ge 0\}$ issue de $0$, dont la probabilit\'e de transition est le noyau $Q_{\Lambda_0}$ d\'efini en (\ref{qaffine}) avec $h=\frac{1}{n}d$.  On consid\`ere les projections $\widetilde X^{(n)}$ et $\widetilde Z^{(n)}$  de ces processus sur $\R\Lambda_0\oplus \mathfrak{h}_\R^*$. On a alors le th\'eor\`eme suivant. 
\begin{theo}[M.D, \cite{defo2}] \label{TCL-cond} Quand $n$ tend vers l'infini,
\begin{enumerate} 
\item la suite de processus $\{\frac{1}{n}\widetilde X^{(n)}([nt]): t\ge 0\}$, $n\ge 0$, converge en loi vers le brownien espace-temps $$\{B_t=t\Lambda_0+b_t:t\ge 0\},$$
\item la suite de processus $\{\frac{1}{n}\widetilde Z^{(n)}([nt]): t\ge 0\}$, $n\ge 0$, converge en loi vers le brownien espace-temps conditionn\'e $$\{A_t=t\Lambda_0+a_t:t\ge 0\}.$$
\end{enumerate}
\end{theo}
Ces processus espace-temps sont tr\`es importants pour nous. Ils r\'eapparaitront au chapitre \ref{chap-drap}, construits \`a partir d'une alg\`ebre affine quelconque et au chapitre \ref{chap-pit} en prenant $\mathfrak{g}$ \'egale \`a $\mathfrak{sl}_2(\C)$.

 \section{Le cas de l'alg\`ebre affine $\widehat{\mathcal L}(\mathfrak{sl}_2(\C))$}\label{section-sl2-affine}
\begin{figure}[!t]
\centering
\includegraphics[scale=0.9]{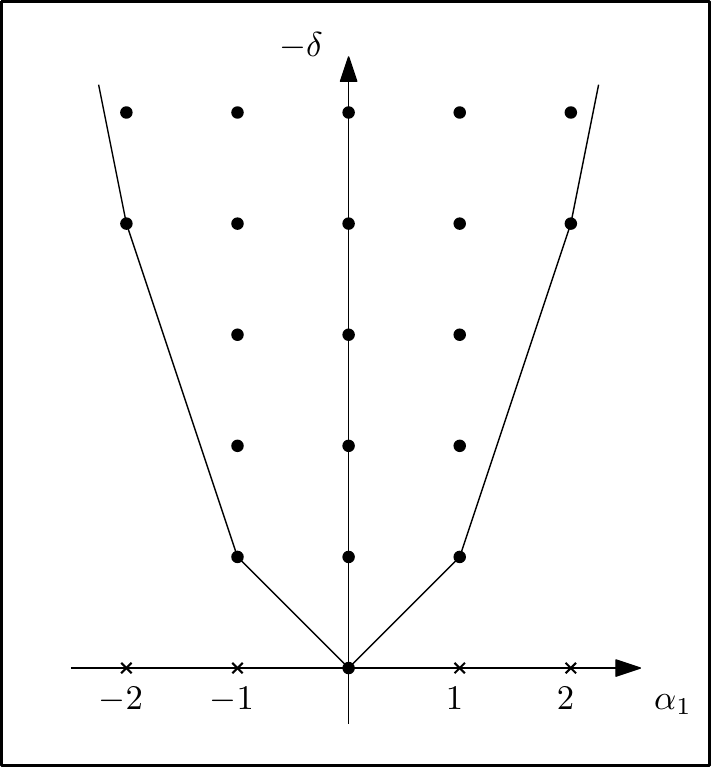}
\caption{Support de $\mu_{\Lambda_0}$} 
\label{RF}
\end{figure}

\begin{figure}[!t]
\centering
\includegraphics[scale=0.7]{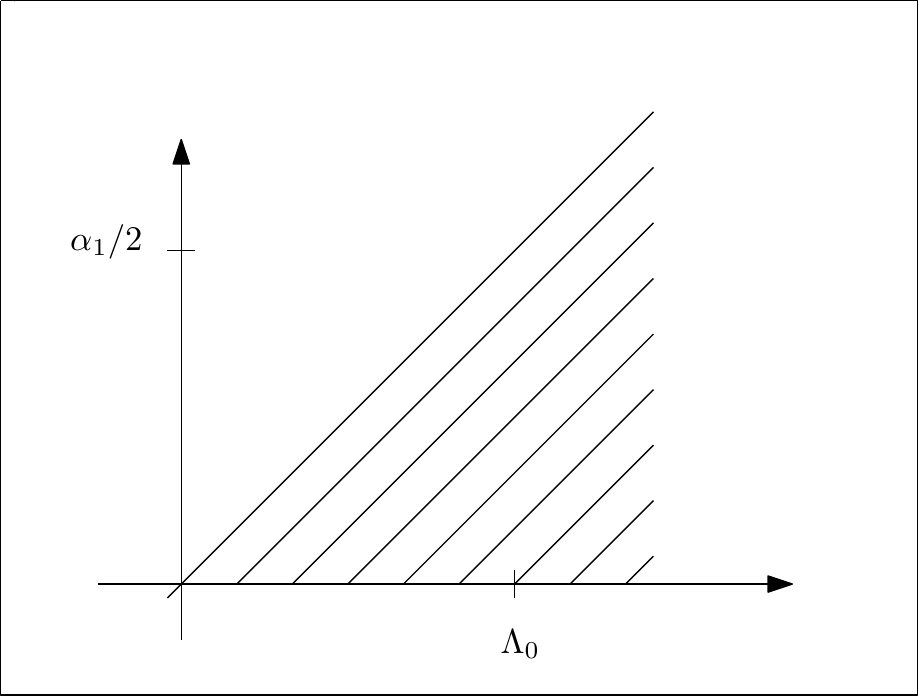}
\caption{Chambre de Weyl dans $\R\Lambda_0\oplus \R\alpha_1/2$} 
\label{CAC}
\end{figure}

Dans cette section nous   pr\'ecisons    ce que sont  les processus  pr\'ec\'edents lorsque l'alg\`ebre affine est  l'alg\`ebre affine associ\'ee \`a $\mathfrak{sl}_2(\C)$. Ces processus r\'eappara\^itront  dans le chapitre \ref{chap-pit} et nous pourrons ainsi nous y r\'ef\'erer plus facilement.  Nous reprenons les notations de  la section \ref{section-affA11} du chapitre \ref{chap-Algebre-affine}. Dans ce cas, $\widehat{\mathfrak h}_\R=\R c\oplus  \R\alpha_1^\vee \oplus\R d$ et $\widehat{\mathfrak h}^*_\R=\R \Lambda_0\oplus \R\alpha_1\oplus\R \delta$. L'espace $\R^3$  et son dual fournissent une r\'ealisation de ces espaces. On munit $\R^3$ du produit scalaire usuel. En identifiant $\R^3$ et son dual on peut poser
$$d=(0,0,1), \, c=(1,0,0), \,\alpha_1^\vee=(0,1,0)$$ et $$\Lambda_0=(1,0,0), \, \delta=(0,0,1), \, \alpha_1=(0,2,0).$$  Nous choisissons pour d\'efinir les mesures sur les poids et les plus hauts poids une suite $\{h_n:n\ge 1\}$, $h_n=\frac{2}{n}d$. Le facteur $2$ que nous ajoutons ici a pour simple vocation de faire apparaitre un brownien standard le long de $\alpha_1/2$ dans le th\'eor\`eme  limite. C'est la normalisation que nous consid\'erons dans le chapitre \ref{chap-pit}.
\paragraph{Marches et   processus de Markov.} Les accroissements de la marche al\'eatoire $\{X^{(n)}(k): k\ge 0\}$ issue de z\'ero sont distribu\'es selon la mesure $\mu_{\Lambda_0}$ d\'efinie en (\ref{paffine}), avec $\omega=\Lambda_0$, et $h=h_n$. Le support de $\mu_{\Lambda_0}$ est    l'intersection de l'enveloppe convexe  de $\widehat W\cdot\{\Lambda_0\}$ avec l'ensemble $\Lambda_0+\Z\alpha_1+\Z \delta$, c'est-\`a-dire l'ensemble des points \`a coordonn\'ees enti\`eres de  
 $$\{\Lambda_0-tx\alpha_1-(1-t)x^2\delta: x\in \Z, t\in[0,1]\}.$$
Il est  repr\'esent\'e \`a la figure \ref{RF},  la coordonn\'ee le long de $\Lambda_0$ \'etant omise.
La marche  $\{X^{(n)}(k): k\ge 0\}$ est \`a valeurs dans $$\{k\Lambda_0+x\frac{\alpha_1}{2}+y\delta, k\in \N,x,y\in \Z\}=\N\times \Z^2.$$ Le processus de Markov $\{Z^{(n)}(k): k\ge 0\}$ dont le noyau de transition est donn\'e par l'identit\'e (\ref{qaffine})  dans laquelle on aura pos\'e $\omega=\Lambda_0$, et $h=h_n$, est \`a valeurs dans  
$$\{k\Lambda_0+x\frac{\alpha_1}{2}+y\delta, 0\le x\le k, k,x\in \N,y\in \Z\}.$$
Par ailleurs, pour pour tout $k\in \N$, les coordonn\'ees de $X^{(n)}(k)$ et  de $Z^{(n)}(k)$ le long de $\Lambda_0$ sont toutes les deux $k$. 
\paragraph{Brownien espace-temps et brownien espace-temps conditionn\'e.} La projection $\widetilde Z^{(n)}(k)$ de  $Z^{(n)}(k)$ sur $\R\Lambda_0\oplus \mathfrak{h}_\R^*$ est dans $\widetilde C_W$ qui est ici d\'efini par
$$\widetilde C_W=\{t\Lambda_0+x\alpha_1/2: 0\le x\le t\}=\{(t,x)\in \R^2: 0\le x\le t\}.$$
On   repr\'esente $\widetilde C_W$ en figure \ref{CAC}. L'axe temporel $\R\Lambda_0$ repr\'esent\'e en figure \ref{CAA} par l'axe des  ordonn\'ees est ici repr\'esent\'e par celui des abscisses, comme c'est l'habitude  pour les repr\'esentations graphiques des processus.  
On d\'efinit une  fonction $\widehat\varphi_{2d}$ sur $\widetilde C_W$ en posant 
$$\widehat\varphi_{2d}(t,x)=\sum_{k\in \Z}(x+2kt)e^{-2(kx+k^2t)}, \quad (t,x)\in  \widetilde C_W.$$
Et on consid\`ere le processus $$\{t\Lambda_0+b_t\frac{\alpha_1}{2}: t\ge 0\}=\{(t,b_t): t\ge 0\},$$ o\`u $\{b_t: t\ge 0\}$ est un brownien r\'eel standard sur $\R$. La fonction $\widehat \varphi_{2d}$ est   harmonique  pour le processus $\{(t,b_t): t\ge\}$ tu\'e sur les bords de $\widetilde C_W$. Elle est nulle sur les bords de $\widetilde C_W$ et de signe constant sur $\widetilde C_W$. Le processus espace-temps   
\begin{align}\label{cond-sl2}
\{t\Lambda_0+a_t\alpha_1/2: t\ge 0\}=\{(t,a_t): t\ge 0\}
\end{align}
 est la transformation de Doob via $\widehat\varphi_{2d}$ du  brownien espace-temps $\{(t,b_t):t\ge 0\}$  tu\'e sur les bords de $\widetilde C_W$. C'est un processus \`a valeurs dans $\widetilde C_W$ pass\'e le temps initial. Le th\'eor\`eme \ref{TCL-cond} devient ici le th\'eor\`eme de convergence en loi finie dimensionnelle suivant.
\begin{theo}
Quand $n$ tend vers l'infini, les processus  $\{\frac{1}{n}\widetilde X^{(n)}([nt]): t\ge 0\}$ et $\{\frac{1}{n}\widetilde Z^{(n)}([nt]): t\ge 0\}$ convergent respectivement vers      $\{(t,b_t): t\ge 0\} $ et $ \{(t,a_t): t\ge 0\}.$
\end{theo}
\chapter{Formule de localisation en dimension infinie\\ \vspace{0.5cm}
\small{\it  o\`u l'on mise  sur les orbites}}\label{chap-loc}

 \bigskip
 \bigskip
 \bigskip
 
Les formules de localisation de Duistermaat et de Heckman \cite{DH} interviennent dans un cadre   symplectique, une forme symplectique sur une vari\'et\'e de dimension finie fournissant une mesure naturelle sur cette vari\'et\'e : la mesure de Liouville. Ces formules permettent d'exprimer la transform\'ee de Fourier de l'image de cette mesure par une application  moment associ\'ee \`a l'action hamiltonienne d'un tore sur la vari\'et\'e  comme une somme portant sur les points critiques du champ fondamental de cette action. Selon les mots de Mich\`ele Audin   \cite{Audin2}, les orbites coadjointes constituent pour les formules de Duistermaat et de Heckman \guillemotleft a machine producing examples\guillemotright. Par ailleurs, la m\'ethode dite des orbites de Kirillov, qui est davantage une philosophie qu'un th\'eor\`eme, parie   sur une correspondance entre   les orbites coadjointes pour l'action d'un groupe de Lie $G$ sur son alg\`ebre de Lie et les  repr\'esentations irr\'eductibles de $G$. Dans le cas  d'un groupe compact, la formule des caract\`eres de Kirillov et  la formule de localisation co\"incident. 

Igor Frenkel \cite{Frenkel} a le premier \'etabli \`a la mani\`ere de Kirillov  une   correspondance entre int\'egrale orbitale et caract\`ere de repr\'esentation pour un groupe $G$ \'egal \`a  un groupe de lacets  agissant sur une alg\`ebre de Kac--Moody affine, c'est-\`a-dire une alg\`ebre de lacets avec extension centrale et d\'erivation pouvant \^etre consid\'er\'ee comme l'alg\`ebre de Lie de $G$. Les formules de localisation de Duistermaat et de Heckman ne sont pas valables   de mani\`ere rigoureuse dans le cadre infini-dimensionnel consid\'er\'e par Igor Frenkel. Cependant,  beaucoup des objets que le th\'eor\`eme de localisation  met en jeu poss\`edent un analogue d\'efini dans le cadre des alg\`ebres de lacets. En les substituant aux objets intervenant dans la formule de localisation usuelle, on obtient au moins formellement une formule de type localisation pour les alg\`ebres de lacets qui est la formule de Kirillov--Frenkel\footnote{Jean-Michel  Bismut a \'egalement propos\'e une preuve de la formule de Kirillov-Frenkel. Son approche est expliqu\'ee par exemple dans \cite{Bismut}.}.

Il est selon moi int\'eressant d'introduire la formule de Kirillov--Frenkel en la pr\'esentant  comme une formule de localisation en dimension infinie et nous d\'eveloppons ce point  de vue dans ce chapitre qui n'a pas vocation \`a pr\'esenter des \'enonc\'es  compl\`etement rigoureux.  Nous consid\'erons en effet l'action sur une alg\`ebre de lacets d'un groupe de lacets   \`a valeurs dans un groupe compact. Espaces tangents,  diff\'erentielles etc, sont \og d\'efinis\fg \, par simple analogie avec le cas fini-dimensionnel. De m\^eme, nous appelons chemin toute application suffisamment r\'eguli\`ere, sans pr\'eciser le degr\'e de r\'egularit\'e.  Nous pr\'esenterons dans le chapitre suivant la formule de Kirillov--Frenkel telle qu'elle est pr\'esent\'ee dans \cite{Frenkel} et reprise dans \cite{defo3}.  Aucun th\'eor\`eme original n'est \'enonc\'e dans ce chapitre o\`u l'on fait simplement pour changer de point de vue un petit pas de c\^ot\'e.

Dans ce chapitre $G_0$ est un groupe de Lie compact connexe  simplement connexe d'alg\`ebre de Lie $\mathfrak g_0$.  On supposera sans perte de g\'en\'eralit\'e que $G_0$ est un groupe de matrices. L'action adjointe de $G_0$ sur son alg\`ebre de Lie not\'ee $\Ad$ est alors l'action par conjugaison. Nous notons $\Ad^*$ l'action codadjointe qui en d\'ecoule et nous munissons $\mathfrak g_0$  d'un produit scalaire $\Ad(G_0)$-invariant qu'on note~$(\cdot \vert \cdot)$.

  \section{Formule de localisation de Duistermaat et de Heckman}\label{sectionFormule-DH}
 Le cadre usuel est le suivant. On dispose d'une vari\'et\'e symplectique compacte  $(M,\omega)$ de dimension $2n$ et d'un tore $T$ agissant sur $M$ de fa\c con hamiltonienne, l'application moment associ\'ee \'etant ici not\'ee  $H$. On note $\mathfrak{t}$ l'alg\`ebre de Lie de $T$. On rappelle que $H$ est une application de $M$ dans $\mathfrak t^*$.  Pour $X\in \mathfrak t$,  $M_0(X)$ d\'esigne l'ensemble des points critiques en $X$, c'est-\`a-dire l'ensemble des points de $M$ o\`u le champ fondamental $\underline X$ associ\'e \`a $X$ s'annule.     Le th\'eor\`eme de Duistermaat et de Heckman tel qu'il est repris dans \cite{Audin} s'\'enonce ainsi. 
 
\begin{theo}\label{DH}  Soit $X\in \mathfrak t$ tel que $M_0(X)$ soit de cardinal fini. Alors  
$$\int_M e^{ H^X}\, \omega^{\wedge n}\propto  \sum_{\xi\in M_0(X)} \frac{e^{ H^X(\xi)}}{\Pf(\mathcal L_\xi(X))},$$ 
o\`u $\propto$ signifie proportionnelle \`a\footnote{la constante de proportionnalit\'e ne d\'ependant pas de $X$.}, $H^X(\xi)=H(\xi)(X)$ et $\mathcal L_\xi(X)\in \mathcal L(T_\xi M,T_\xi M)$, est l'action infinit\'esimale de $X$ sur $T_\xi M$ l'espace tangent  \`a $M$ en $\xi$, $\Pf(\mathcal L_\xi(X))$ d\'esignant son pfaffien, pour l'orientation associ\'ee \`a $\omega^{\wedge n}$. 
\end{theo}
 
 La mesure image par $H$ de la mesure $\omega^{\wedge n}$ est appel\'ee mesure de Duistermaat--Heckman associ\'ee \`a l'action hamiltonienne de $T$ sur $M$. Consid\'erons l'exemple classique d'une vari\'et\'e $M$ \'egale \`a une orbite  dans $\mathfrak g_0^*$ pour l'action coadjointe de $G_0$ sur $\mathfrak{g}_0^*$. Si $T$ est un tore maximal de $G_0$ d'alg\`ebre de Lie $\mathfrak t$, la formule de localisation de Duistermaat et de Heckman, pour une application moment \'egale \`a la projection canonique  sur $\mathfrak t^*$,  devient la formule de Harish-Chandra pour les groupes compacts, ou encore la formule des caract\`eres de Kirillov.  Dans ce contexte, la formule de localisation porte sur la transform\'ee de Fourier de  la mesure image par la projection canonique sur $\mathfrak t^*$ de la mesure uniforme sur une orbite coadjointe. Nous d\'etaillons ci-apr\`es le cas o\`u $G_0=\SU(n)$.
\paragraph{Formule de localisation pour $G_0=\SU(n)$.}  Notons $\mathfrak{su}(n)$ l'alg\`ebre de Lie de $\SU(n)$. Pour $\lambda\in \mathfrak{su}(n)^*$ consid\'erons    l'orbite  $\mathcal{O}_\lambda$ de $\lambda$ sous l'action de $\SU(n)$, i.e. $$\mathcal{O}_\lambda=\{\Ad^*(u)(\lambda): u\in \SU(n)\},$$
o\`u $\Ad^*$ d\'esigne l'action coajointe de $\SU(n)$ sur $\mathfrak{su}(n)^*$. Pour $X\in \mathfrak{su}(n)$, on note $\underline X$ le champ fondamental associ\'e \`a $X$ qui est d\'efini par  $$\underline X_\xi=\frac{d}{dt}\Ad^*(e^{tX})(\xi)_{\vert_{t=0}}, \quad \xi\in  \mathcal O_\lambda.$$ Ce champ fondamental est donc d\'efini par la relation 
$$\langle \underline X_\xi,Y\rangle = \langle \xi,[Y,X]\rangle,\quad \xi\in \mathcal O_\lambda,\,Y\in \mathfrak{su}(n),$$ o\`u $\langle.,.\rangle$ est l'appariement dual canonique. L'espace tangent  \`a $\mathcal O_\lambda$ en $\xi$ est $$T_\xi\mathcal O_\lambda=\{ \underline X_\xi: X\in \mathfrak{su}(n)\}.$$  On choisit $\lambda\in \mathfrak{su}(n)^*$ distinct du vecteur nul. L'orbite est ainsi non r\'eduite \`a un point. On d\'efinit une forme symplectique $\omega$ sur $\mathcal O_\lambda$ en posant 
$$\omega_\xi(\underline X_\xi,\underline Y_\xi)=\langle \xi,[X,Y]\rangle, \quad \xi \in \mathcal O_\lambda.$$ 
Alors $(\mathcal O_\lambda,\omega)$ est une vari\'et\'e symplectique et l'action coadjointe de $\SU(n)$ sur $\mathcal O_\lambda$ est symplectique. L'application identit\'e est une application moment pour l'action coadjointe de $\SU(n)$ sur $\mathcal O_\lambda$. 
Consid\'erons   un tore maximal  $T$ de $\SU(n)$, qu'on choisit   \'egal \`a l'ensemble des matrices diagonales de $\SU(n)$.  La restriction de l'application moment \`a l'alg\`ebre de Lie $\mathfrak t$ du tore reste une application moment.   On suppose que $\lambda\in \mathfrak t^*$.
  On choisit un \'el\'ement $X$ de  $\mathfrak{t}$ dont les valeurs propres   sont distinctes, ce qui assure que $M_0(X)$ soit de cardinal fini. On a 
   $$M_0(X)=\{\Ad^*(\boldsymbol \sigma) (\lambda): \sigma \in \mathfrak{S}_n \},$$
   o\`u $\mathfrak{S}_n$ est  le groupe des permutation de  $[\![1;n]\!]$ et $\boldsymbol \sigma$ est un \'el\'ement  de $\SU(n)$ dont l'action adjointe sur $\mathfrak t$ est l'action de permutation des \'el\'ements diagonaux par $\sigma$.  L'action infinit\'esimale de $X$ en   $\xi\in M_0(X)$   est donn\'ee par    
   \begin{align*}
   \mathcal L_\xi(X)(\underline Y_\xi)&=\frac{d}{dt}\big(\Ad^*(e^{tX})(\frac{d}{ds}\Ad^*(e^{sY})\xi_{\vert_{s=0}}\big)_{\vert_{t=0}}\\
   &= \underline{\mbox{ad}_X(Y)}_\xi.
   \end{align*}
Les valeurs propres de  $\mathcal L_\xi(X)$ sont $\{\pm i(x_l-x_k): 1\le l< k\le n\}$, o\`u les $ix_k$ sont les valeurs propres de $X$. La mesure de    Liouville \'etant  $\Ad$-invariante  le th\'eor\`eme \ref{DH} donne la formule suivante, qui est la formule  de Harish-Chandra pour le groupe sp\'ecial unitaire.  
\begin{align}\label{AC-SU}
\int_{\mbox{\small SU}(n)} e^{\langle \Ad^*(u)(\lambda),X\rangle }\, du\propto \frac{1}{\Delta(\lambda)\Delta(X)}\sum_{\sigma \in \mathfrak{S}_n}\det(\sigma)e^{\langle \mbox{\small Ad}^*(\boldsymbol\sigma)( \lambda),X\rangle},\end{align}
o\`u $\Delta(\lambda)=\prod_{i<j}\langle\lambda,E_{ii}-E_{jj}\rangle,$ $\Delta(X)=\prod_{i<j}(x_{i}-x_{j})$.
Si l'on interpr\`ete la somme dans le membre de droite  comme le num\'erateur d'un caract\`ere de repr\'esentation du groupe $\SU(n)$ tel qu'il appara\^it dans la formule des caract\`eres de Weyl, on reconna\^it une formule des caract\`eres de  Kirillov \cite{K}. Ici la mesure de Duistermaat--Heckman normalis\'ee est l'image par la projection canonique sur $\mathfrak t^*$ de la mesure de probabilit\'e uniforme sur $\mathcal O_\lambda$  et la formule (\ref{AC-SU}) en donne la transform\'ee de Fourier (en rempla\c cant $X$ par $iX$).

\section{Alg\`ebres de lacets, orbites coadjointes} \label{section-lacets-orbites}
 Nous identifions le cercle $S^1$ \`a $\R/ \Z$. Nous appelons  boucle ou lacet, tout chemin d\'efini sur $S^1$ suffisament r\'egulier et consid\'erons  l'ensemble   $L(G_0)$ (resp. $L(\mathfrak{g}_0)$) des boucles   \`a valeurs dans $G_0$  (resp. $\mathfrak{g}_0$). $L(\mathfrak{g}_0)$ est un $\R$-espace vectoriel que nous munissons d'un produit scalaire $\Ad(G_0)$-invariant toujours not\'e $(\cdot\vert \cdot)$ en posant 
$$(\eta\vert \xi)=\int_0^1(\eta(s)\vert\xi(s))\, ds, \quad \eta,\xi\in L(\mathfrak g_0).$$ Nous construisons l'extension centrale avec d\'erivation de l'alg\`ebre de lacets \`a valeurs dans  $\mathfrak{g}_0$  comme nous l'avons fait dans le chapitre \ref{chap-Algebre-affine}.  C'est cette fois une alg\`ebre de Lie sur $\R$.
\paragraph{Alg\`ebre de lacets avec extension centrale et d\'erivation.} Le crochet de Lie   sur $\mathfrak{g}_0$ est not\'e  $[\cdot,\cdot]_{\mathfrak{g}_0}$.  On consid\`ere l'extension centrale  avec d\'erivation   $$\widehat{L}(\mathfrak g_0)=L(\mathfrak g_0)\oplus \mathbb Rc\oplus \R d$$ qu'on munit d'un crochet de Lie $[\cdot,\cdot]$ en posant 
\begin{align}\label{Liehat} 
[\xi+\lambda c+\beta d,\eta+\mu c+\nu d]=[\xi,\eta]_{\mathfrak{g}_0}+ (\xi'\vert \eta)c+\beta \eta'-\nu\xi',
\end{align}
pour $\xi,\eta\in L(\mathfrak g_0)$, $\lambda,\beta,\mu,\nu\in \R$,  o\`u $[\xi,\eta]_{\mathfrak{g}_0}$ est  d\'efini point par point. Munie de ce crochet $\widehat{L}(\mathfrak g_0)$ est une $\R$-alg\`ebre de Lie. Une sous-alg\`ebre de Cartan de $\widehat{L}(\mathfrak g_0)$ est la sous-alg\`ebre    $\mathfrak t\oplus \R c\oplus \R d$, o\`u  $\mathfrak t$ est identifi\'e  \`a l'ensemble des boucles constantes \`a valeurs dans $\mathfrak t$.  Le crochet de Lie d\'efinit une action qu'on note $\mbox{ad}$ de   $L(\mathfrak g_0) \oplus \R c\oplus \R d$ sur elle m\^eme qui, restreinte \`a $L(\mathfrak g_0)\oplus \R d$, est donn\'ee par 
$$\mbox{ad}(\xi+\beta d)(\eta+\mu c+\nu d)=[\xi,\eta]_{\mathfrak{g}_0}+ (\xi'\vert \eta)c+\beta \eta'-\nu\xi',$$ 
pour tout  $\xi,\eta \in L(\mathfrak g_0),\beta,\mu,\nu\in \R$. 

\paragraph{Action  de $S^1\ltimes L(G_0)$.} On consid\`ere les actions de $S^1$ sur $L(G_0)$ et sur $L(\mathfrak{g}_0)$  d\'efinies par 
$$(\theta\cdot\gamma)(s)=\gamma(s+\theta),\quad (\theta\cdot\eta)(s)=\eta(s+\theta), \quad \theta,s\in S^1,$$ 
$\gamma\in L(G_0)$,  $\eta\in L(\mathfrak g_0)$. Nous consid\'erons le produit semi-direct $S^1\ltimes L(G_0)$, c'est-\`a-dire le groupe produit  $S^1\times L(G_0)$ muni d'une op\'eration de groupe donn\'ee par 
$$(\theta_1,\gamma_1)(\theta_2, \gamma_2)=(\theta_1+\theta_2,\gamma_1(\theta_1\cdot\gamma_2)),\quad \theta_1,\theta_2\in S^1, \gamma_1,\gamma_2\in L(G_0),$$
 et d\'efinissons une action $\Ad_{0}$ de $S^1\ltimes L(G_0)$ sur $L(\mathfrak g_0)$ en posant pour $\theta\in S^1$ et $\gamma\in L(G_0)$,
 $$\Ad_{0}(\theta,\gamma)\eta=\gamma(\theta\cdot\eta)\gamma^{-1}, \quad \eta\in L(\mathfrak g_0).$$
L'action coadjointe $\Ad^*_{0}$ de $S^1\ltimes L(G_0)$ sur $L(\mathfrak g_0)^*$ est donn\'ee par 
\begin{align*}
\Ad^*_{0}(\theta,\gamma)\phi&=\phi\circ\Ad_{0}((\theta,\gamma)^{-1})\\
&=\phi\circ\Ad_{0}(-\theta,(-\theta)\cdot\gamma^{-1}).
\end{align*}
On notera $\Ad^*_0(\gamma)$ au lieu de $\Ad^*_0(0,\gamma)$. On d\'efinit une action $\Ad$ de $ S^1\ltimes L(G_0)$ sur  $L(\mathfrak g_0) \oplus \R c\oplus \R d$ en posant 
\begin{align}\label{adjoint}
 \Ad(\theta,\gamma)(\eta+\mu c+\nu d)&=\mbox{Ad}_{0}(\theta,\gamma)(\eta)-\nu \gamma'\gamma^{-1}+\nu d \nonumber\\
 &       \quad +\big(\mu +(\gamma^{-1}\gamma'\vert\theta\cdot\eta)-\frac{\nu}{2}(\gamma'\gamma^{-1}\vert \gamma'\gamma^{-1})\big)c,
\end{align}
pour tout $\gamma\in L (G_0),(\eta,\lambda,\nu)\in L(\mathfrak g_0)\times \R\times \R$.  C'est l'action adjointe de   $S^1\ltimes L(G_0)$ sur $L(\mathfrak g_0) \oplus \R c\oplus \R d$. Elle v\'erifie\footnote{Nous surtout.}
$$\mbox{ad}(\theta d+\xi)=\frac{d}{dt}\Ad(t\theta,e^{t\xi})\vert_{t=0}.$$
 Consid\'erons  l'espace dual  $\widehat{L}(\mathfrak{g}_0)^*=L(\mathfrak{g}_0)^*\oplus \R\Lambda_0\oplus \R\delta$, o\`u  $$\Lambda_0(L(\mathfrak{g}_0))=\delta(L(\mathfrak{g}_0))=\delta(c)=\Lambda_0(d)=0,\,\Lambda_0(c)=\delta(d)=1.$$
 L'action coadjointe $\Ad^*$ de  $S^1\ltimes L(G_0)$ sur  $\widehat L(\mathfrak{g}_0)^*$ est donn\'ee par 
 \begin{align*}
 \Ad^*(\theta,\gamma)(\phi+\tau \Lambda_0+\nu \delta)&=(\phi+\tau \Lambda_0+\nu \delta)\circ\Ad((\theta,\gamma)^{-1})\\
 &= (\phi+\tau \Lambda_0+\nu \delta)\circ\Ad(-\theta,(-\theta)\cdot\gamma^{-1}).
 \end{align*}
 On obtient ainsi
 \begin{align}\label{coadjoint} 
\Ad^*(\theta,\gamma)(\phi+\tau \Lambda_0+\nu \delta)& =[\Ad_{0}^*(\theta,\gamma)\phi-\tau (\gamma'\gamma^{-1}\vert\cdot)]+\nu \delta+\tau\Lambda_0\nonumber \\ &\quad \quad +(\phi((-\theta)\cdot\gamma^{-1}\gamma')-\frac{\tau}{2}(\gamma'\gamma^{-1}\vert\gamma'\gamma^{-1}))\delta
\end{align} 
pour tout  $\phi\in L(\mathfrak g_0)^*, \tau\in \R$. On   notera  $\Ad^*(\gamma)$ au lieu de $\Ad^*(0,\gamma)$. 
Pour $a\in \mathfrak t$, en notant $\phi_a=(a\vert\cdot)$, alors, pour    $\gamma\in L(G_0)$, $\theta\in S^1$, $\Ad^*(\theta,\gamma)(\phi_a+\tau\Lambda_0)$ vaut
\begin{align}\label{orbite-a}
(\gamma a\gamma^{-1}-\tau\gamma'\gamma^{-1}\vert\cdot)+\tau\Lambda_0
-\frac{1}{2\tau} (\vert\vert \gamma a\gamma^{-1}-\tau\gamma'\gamma^{-1}\vert\vert^2-\vert\vert a\vert \vert^2 )\delta.
\end{align}
On remarque que cette quantit\'e ne d\'epend pas de $\theta$.  La coordonn\'ee le long de $\delta$ est une fonctionnelle d'\'energie. Notons par ailleurs que les actions coadjointes n'affectent pas la coordonn\'ee le long de $\Lambda_0$\footnote{Et que ce sont des actions !}.  
\paragraph{Orbite coadjointe.} On consid\`ere pour $\zeta\in \widehat L(\mathfrak g_0)^*$ l'orbite
$$\widehat{\mathcal O}_{\zeta}=\{\Ad^*(\theta,u)(\zeta): (\theta,u)\in S^1\ltimes L(G_0)\}.$$
Le champ fondamental $\underline{\theta d+X} $ associ\'e \`a $\theta d+X$ est donn\'e par
$$\langle \underline{\theta_1 d+X_1}_\xi,\theta_2 d+X_2\rangle =\langle \xi ,[\theta_1 d+X_1,\theta_2 d+X_2]\rangle,$$
$\xi \in\widehat{\mathcal O}_{\zeta},$ $\theta_2\in \R$, $X_2\in L(\mathfrak g_0)$. L'espace tangent en $\xi\in\widehat{ \mathcal O}_{\zeta}$ est
$$T_\xi \widehat{\mathcal O}_{\zeta} =\{\underline{\theta d+X}_\xi:\theta\in \R, X\in L(\mathfrak g_0)\},$$
et on munit $\widehat{\mathcal O}_{\zeta}$ d'une forme symplectique $\omega$ en posant pour $\xi \in \widehat{\mathcal O}_{\zeta}$,
$$\omega_\xi(\underline{\theta_1d+ X_1}_\xi,\underline{\theta_2 d+ X_2}_\xi)=\langle \xi,[\theta_2 d+X_2,\theta_1d + X_1]\rangle ,$$ 
$ \theta_1,\theta_2\in \R,$  $X_1,X_2\in L(\mathfrak g_0).$
\paragraph{Application moment.} L'action  de $S^1\ltimes L(\mathfrak g_0)$ sur $\widehat{\mathcal O}_{\zeta}$   est une action hamiltonienne, l'application identit\'e \'etant une application moment pour cette action. Ainsi l'action  du tore $S^1\ltimes T$ sur $\widehat{\mathcal O}_{\zeta}$  est hamiltonienne,    l'application $H$ de projection canonique sur $\R \Lambda_0\oplus \mathfrak t^*$ d\'efinie par
$$H(\xi)=\xi_{\vert_{\R d \oplus \mathfrak t}},\quad \xi \in \widehat{\mathcal O}_{\zeta},$$
\'etant une application moment pour cette action. 
\paragraph{Action infinit\'esimale.} En un point critique $\xi\in \widehat{\mathcal O}_{\zeta}$, l'action infinit\'esimale de $\theta d+X$ sur $T_\xi\widehat{ \mathcal O}_{\zeta}$, pour $\theta\in \R$, $X\in L(\mathfrak g_0)$ est  donn\'ee par 
$$\mathcal{L}_\xi (\theta d+X)(\underline{\kappa d+Y}_\xi)=\underline{\mbox{ad}(\theta d+X)(\kappa d+Y)}_\xi,$$
$\kappa\in \R$, $Y\in L(\mathfrak g_0)$.
 
\section{Une formule de localisation en dimension infinie}
Nous allons maintenant substituer aux orbites coadjointes de dimension finie pour lesquelles la formule de localisation de Duistermaat et de Heckman du th\'eor\`eme \ref{DH} est valable, les orbites coadjointes consid\'er\'ees \`a la section pr\'ec\'edente, et obtenir dans ce contexte, formellement au moins, une formule de localisation. 
Avec les notations de la section pr\'ec\'edente, pour $\theta\in \R$, $X\in \mathfrak{t}$, et $\zeta \in\widehat{L}( \mathfrak{g}_0)$, tels que  $M_0(\theta d+X)$ soit de cardinal fini, la formule du th\'eor\`eme \ref{DH} deviendrait
 \begin{align}\label{DH-infini} \int_{\widehat{\mathcal O}_{\zeta}}e^{  H^{(\theta,X)}} \omega^{\wedge \infty}\propto \sum_{\xi\in M_0(\theta d+X)}\frac{e^{  H^{(\theta,X)}(\xi)}}{\Pf(\mathcal L_\xi(\theta d+X))},
 \end{align}
 o\`u $H^{(\theta,X)}(\xi)=\xi(\theta d+X)$, $\xi \in \widehat{\mathcal O}_\zeta$.  Cette formule fait intervenir  plusieurs quantit\'es non d\'efinies. Dans le membre de gauche, la mesure contre laquelle on int\`egre, si elle existait,  devrait  \^etre une mesure $L( G_0)$-invariante  et  on sait qu'une telle mesure n'existe pas.    Dans celui de droite, le pfaffien  de l'op\'erateur lin\'eaire en dimension infini n'est pas d\'efini.  Il existe cependant un produit appel\'e produit z\^eta-regularis\'e qui permet de donner un sens \`a ce pfaffien. Il est utilis\'e dans un contexte proche par Wendt dans \cite{Wendt}. Nous allons donc essayer de donner un sens acceptable aux quantit\'es apparaissant dans (\ref{DH-infini}).

Afin de ne pas asphyxier le lecteur avec de nouvelles  d\'efinitions, on consid\`ere \`a partir de maintenant et jusqu'\`a la section \ref{subsectionorbite} que   $G_0$ est le groupe sp\'ecial unitaire $\SU(n)$. On choisit un tore maximal $T$ dans $\SU(n)$ \'egal \`a l'ensemble des matrices diagonales de $\SU(n)$.  On munit $\mathfrak{su}(n)$ du produit scalaire $(M\vert N)=\frac{1}{(2\pi)^2}\mbox{tr}(MN^*)$,  $M,N\in \mathfrak{su}(n)$.
\subsection{La somme discr\`ete}
On consid\`ere pour $\tau\in \R$ et $a\in \mathfrak t$, l'orbite $\widehat{\mathcal O}_{\tau \Lambda_0+\phi_a}$, o\`u $\phi_a=(.\vert a)$.
\paragraph{Points critiques.} Soit $\theta$ un r\'eel strictement positif et $b$ une matrice diagonale de $   \mathfrak{su}(n)$ dont les valeurs propres $2i\pi b_1,\dots,2i\pi b_n$ satisfont 
\begin{align}\label{cond1}  \left\{
    \begin{array}{l}
  \sum_{k=1}^{n}b_k=0,\\
 b_1>\dots>b_n>b_1-\theta.
    \end{array}
\right.
\end{align}
Posons pour $m\in \Z $, $k,l\in \{1,\dots,n\}$, avec $k< l$,
\begin{align*}
X^m_{k,l}(t)=e^{2im\pi t}E_{kl}-e^{-2im\pi t}E_{lk},\quad Y^m_{k,l}(t)=ie^{2im\pi t}E_{kl}+ie^{-2im\pi t}E_{lk}, 
\end{align*}
et pour $m\ge 1 $, $k,l\in \{1,\dots,n\}$, avec $k< l$,
\begin{align*} 
Z^{m,-}_{k,l}(t)=i\cos(2m\pi t)(E_{kk}-E_{ll}),\quad Z^{m,+}_{k,l}(t)=i\sin(2m\pi t)(E_{kk}-E_{ll}),
\end{align*}
o\`u $\{E_{ij}: i,j\in\{1,\dots,n\}\}$ est la base canonique de $\mathcal M_n(\C)$.
On a  pour $m\in \Z $, $k,l\in \{1,\dots,n\}$, avec $k< l$,
\begin{align}\label{valp1}  \left\{
    \begin{array}{lll}
\mbox{ad}(\theta d+b)(X^m_{k,l})&=&2\pi(m\theta+b_k-b_{l})Y^m_{k,l}\\
\mbox{ad}(\theta d+b)(Y^m_{k,l})&=&-2\pi(m\theta+b_k-b_{l})X^m_{k,l}
    \end{array}
\right.
\end{align}
et pour $m\ge 1$, $k,l\in \{1,\dots,n\}$, avec $k< l$,
 \begin{align}\label{valp2}  \left\{
    \begin{array}{lll} 
\mbox{ad}(\theta d+b)(Z^{m,-}_{k,l})&=&-2\pi m\theta Z^{m,+}_{k,l}\\
\mbox{ad}(\theta d+b)(Z^{m,+}_{k,l})&=&2\pi m\theta Z^{m,-}_{k,l}.
    \end{array}
\right.
\end{align}
 En \'evaluant $\underline{\theta d+b}_\xi$ en  les matrices  $X^m_{k,l},\, Y^m_{k,l}, \, Z^{m,-}_{k,l} \textrm{ et } Z^{m,+}_{k,l}$  et en utilisant les conditions (\ref{cond1}) on obtient  que \begin{center}  $\xi\in \widehat{\mathcal O}_{\tau \Lambda_0+\phi_a}$ est un point critique si et seulement si  $\xi\in \R\Lambda_0\oplus \R\delta\oplus \R\mathfrak t^*$.\end{center} 
Supposons maintenant que $\tau$ est choisi strictement positif et que $ a $ est une matrice diagonale de $\mathfrak{su}(n)$ dont les valeurs propres   $2i\pi a_1,\dots,2i\pi a_n$  satisfont 
\begin{align}\label{cond}  \left\{
    \begin{array}{l}
     \sum_{k=1}^{n}a_k=0,\\
       a_1>\dots>a_n>a_1-\tau.
    \end{array}
\right.
\end{align}
L'expression (\ref{orbite-a}) permet de montrer que l'ensemble des boucles $\gamma\in L(G_0)$ telles que $$\Ad^*(\gamma)(\tau\Lambda_0+\phi_a)\in \R\Lambda_0\oplus \R\delta\oplus \R\mathfrak t^*$$ est le sous-groupe $\widehat W_0$ de $L(G_0)$  d\'efini par 
$$\widehat W_0=\{\gamma_{\sigma,x} \in L(G_0):x\in \Lambda,\sigma\in \mathfrak{S}_n\},$$
o\`u $\gamma_{\sigma,x}(s)=\boldsymbol\sigma e^{sx} $,  $s\in S^1$, et o\`u  $\Lambda$ est le noyau de l'application exponentielle qui est ici le r\'eseau engendr\'e par 
$$ \{ 2i\pi(E_{kk}-E_{(k+1)(k+1)}): k\in \{1,\dots,n-1\}\}.$$
L'action de $\widehat W_0$ s'identifie \`a celle du groupe de Weyl affine $\widehat W$ d\'efinie au chapitre \ref{chap-Algebre-affine}. On a en effet, pour $\sigma\in \mathfrak S_n$ et $x\in \Lambda$
$$\Ad^*(\gamma_{\sigma,x})(\tau\Lambda_0+\phi_a)=(\boldsymbol\sigma(a-\tau  x)\boldsymbol\sigma^{-1}\vert\cdot)+(( x\vert\boldsymbol \sigma a \boldsymbol\sigma^{-1})-\frac{\tau}{2}(x\vert x))\delta.$$
\paragraph{Pfaffien de l'action infinit\'esimale.} Les identit\'es (\ref{valp1}) et  (\ref{valp2})  montrent que l'ensemble des valeurs propres de l'action inifnit\'esimale de $\theta d+b$ sur $T_\xi\widehat{\mathcal O}_{\tau\Lambda_0+\phi_a}$ est
$$\{\pm 2i\pi m: m\ge 1\} \cup\{2i\pi(\pm m\theta\pm (b_k-b_{l})): m\ge 1, k,l\in\{1,\dots,n\}, k<l\},$$
les valeurs propres d\'ependant de $\theta$ et de $b$ \'etant de multiplicit\'e $1$.  Le produit infini des valeurs propres n'est pas d\'efini mais la m\'ethode dite de r\'egularisation z\^eta permet de donner une  version r\'egularis\'ee du pfaffien \cite{Wendt}. On obtient    ici 
\begin{align*}
\Pf(\mathcal{L}_{\xi}(\theta d+b))= \det(\sigma) \prod_{1\le k<l\le n}\sin(\frac{\pi}{\theta}(b_k-b_l)),
\end{align*} 
pour $\xi=\gamma_{\sigma,x}.(\tau\Lambda_0+\phi_a)$, avec $\sigma\in \mathfrak S_n$ et $x\in \Lambda$.
\paragraph{Formule de localisation.}
Dans la formule (\ref{DH-infini}) le membre de droite devient pour $\zeta=\tau\Lambda_0+\phi_a$, $X=b$, 
\begin{align}\label{MDD-DH}
 \frac{1}{  \prod_{1\le k<l\le n}\sin(\frac{\pi}{\theta}(b_k-b_l))}\sum_{\sigma \in \mathfrak S_n, x\in \Lambda}\det(\sigma) e^{ \langle\Ad^*(\gamma_{\sigma,x})(\tau\Lambda_0+\phi_a),\, \theta d+b\rangle}.
 \end{align}
 Cette expression a un sens puisque la s\'erie est convergente. Comme dans (\ref{AC-SU}) la somme est \`a comparer au num\'erateur d'un caract\`ere de repr\'esentation tel qu'il apparait dans la formule des caract\`eres de Weyl (\ref{Weyl}). Cette fois, il s'agit d'un caract\`ere de repr\'esentation  d'alg\`ebre affine.  Il reste \`a donner un sens au membre de gauche de (\ref{DH-infini}).  
\subsection{L'int\'egrale orbitale}\label{subsectionorbite}
 Soit $\tau$ un r\'eel strictement postif  et $a$ un \'el\'ement de $\mathfrak t$. Remarquons d'abord que l'expression (\ref{orbite-a}) montre qu'une forme de $\widehat{\mathcal O}_{\tau\Lambda_0+\phi_a}$ est d\'etermin\'ee par sa projection canonique sur $\widetilde L(\mathfrak g_0)^*$ o\`u   $\widetilde L(\mathfrak g_0)$ est  l'alg\`ebre consid\'er\'ee dans l'introduction du m\'emoire telle que 
$$\widetilde L(\mathfrak g_0)=L(\mathfrak g_0)\oplus\R c.$$
   \'Etant donn\'e un \'el\'ement $\xi$ de $\widehat{\mathcal O}_{\tau\Lambda_0+\phi_a}$, dont la projection sur  $\widetilde L(\mathfrak g_0)^*$    est $$\tau\Lambda_0+\int_0^1(\cdot\vert \dot x_s )ds,$$ avec $\dot x\in L(\mathfrak g_0)$, on d\'etermine $\gamma\in L(G_0)$ telle que $\Ad^*(\gamma)(\tau\Lambda_0+\phi_a)=\xi$ en r\'esolvant 
 $$\tau \, X^{-1}dX=dx,$$
 avec $X(0)=I$, o\`u $I$ est la matrice identit\'e dans $G_0$  et en posant $\gamma(s)=X(s)^{-1}\gamma(0)e^{sa/\tau}$, $s\in S^1$, o\`u $\gamma(0)$ est un \'el\'ement de $G_0$   choisi tel que
 $$X(1)=\gamma(0)\exp(a/\tau)\gamma(0)^{-1}.$$
On \'etablit ainsi une correspondance entre $\widehat{\mathcal O}_{\tau\Lambda_0+\phi_a}$ et l'ensemble des chemins \`a valeurs dans $G_0$ dont la valeur finale est dans l'orbite de $\exp(a/\tau)$ pour l'action adjointe de $G_0$ sur lui-m\^eme.  Cette correspondance est classique. On la trouvera expos\'ee dans    \cite{Segal} par exemple.    

Revenons maintenant \`a l'int\'egrale orbitale de  (\ref{DH-infini}). Nous avons dit que la mesure contre laquelle on int\'egrait dans cette int\'egrale devait \^etre $L(G_0)$-invariante.  S'il existait une telle mesure, en notant  $\mu_a$ sa mesure image par la projection sur $\widetilde L(\mathfrak g_0)^*$ et $\widetilde{\mathcal O}_{\tau\Lambda_0+\phi_a}$ l'image de $\widehat{\mathcal O}_{\tau\Lambda_0+\phi_a}$ par cette projection, l'int\'egrale orbitale s'\'ecrirait, avec $\theta>0$,
\begin{align} \label{intmu}
\int_{\widetilde{\mathcal{O}}_{\tau\Lambda_0+\phi_a}}e^{  \phi(X)}e^{-\frac{ \theta}{2\tau}\vert\vert\phi\vert\vert^2}\mu_a(d\phi).
\end{align}
La mesure 
$$\nu_a(d\phi) =e^{-\frac{\theta}{2\tau}\vert \vert \phi\vert\vert^2}\mu_a(d\phi),$$
 serait quasi-invariante, c'est-\`a-dire v\'erifierait pour $\gamma\in L(G_0)$
$$\frac{d\Ad^*(\gamma)_*\nu_a}{d  \nu_a}=e^{-\theta\phi(\gamma'\gamma^{-1})-\frac{\tau\theta}{2 }\vert\vert\gamma^{-1} \gamma'\vert\vert^2},$$
o\`u $\Ad^*(\gamma)_*\nu_a$ serait la mesure image de $\nu_a$ par $\Ad^*(\gamma)$.  Si on n\'eglige les questions de r\'egularit\'e des trajectoires, on connait une telle mesure : c'est la loi de $$\tau\Lambda_0+ \int_0^1(.\vert \, db_s),$$ o\`u  $\sqrt{\theta/\tau}b$ est un brownien standard sur $\mathfrak g_0$, dont l'enroulement  sur $G_0$ \`a la vitesse $\tau$, est conditionn\'e \`a  avoir une valeur en $1$ dans l'orbite de $\exp(a/\tau)$ pour l'action de conjugaison de $G_0$ sur lui m\^eme. C'est ainsi que Frenkel construit effectivement l'int\'egrale orbitale, la mesure du brownien conditionn\'e par l'orbite du bout de son enroulement   jouant   le r\^ole dans le cadre affine de la mesure de probabilit\'e uniforme sur une $\Ad^*(G_0)$-orbite dans $\mathfrak g_0^*$ du cadre compact  qu'on peut elle aussi d\'ecrire  comme une mesure gaussienne sur $\mathfrak g_0^*$ conditionn\'ee \`a vivre dans  l'orbite. Nous   verrons cela plus pr\'ecis\'ement au chapitre suivant. Notons que  Michael Atiyah et Andrew Pressley d\'emontrent dans  \cite{AP} un  th\'eor\`eme de convexit\'e de type Kostant dans un contexte infini-dimensionnel. C'est dans leur travail je crois qu'appara\^it pour la premi\`ere fois l'\'energie comme application moment d'une action du cercle sur une orbite coadjointe. C'est cette action qui fait donc   surgir la mesure de Wiener comme mesure naturelle sur les orbites codajointes.
 
\chapter{Enroulement de draps browniens\\ \vspace{0.5cm}
\small{\it  o\`u l'on met  et les formes et le temps}}\label{chap-drap}

  \epigraph{Space is a swarming in the eyes, and Time a singing in the ears, says John Shade, a modern poet, as quoted by an invented philosopher ("Martin Gardiner" [sic]) in The Ambidextrous Universe.}{Vladimir Nabokov\footnotemark, \textit{Ada or Ardor: A Family Chronicle}}   
\footnotetext{Nabokov citant   Gardner citant   Nabokov, \c  ca boucle ici aussi.}
Un mouvement brownien r\'eel tu\'e en $0$ conditionn\'e au sens de Doob \`a rester positif a la m\^eme distribution qu'un processus de Bessel  de dimension trois. Consid\'erons plus g\'en\'eralement   l'action coadjointe d'un groupe de Lie compact  connexe semi-simple $G_0$ sur le dual $\mathfrak g_0^*$ de son alg\`ebre de Lie $\mathfrak g_0$, muni  d'un produit scalaire invariant pour cette action. On peut d\'efinir une application \og partie radiale \fg \, sur $\mathfrak g_0^*$ telle que le processus de la partie radiale d'un mouvement brownien sur cet espace ait la m\^eme distribution qu'un brownien sur le dual d'une sous-alg\`ebre de Cartan $\mathfrak t$ de $\mathfrak g_0$, conditionn\'e \`a rester dans une chambre de Weyl associ\'ee au syst\`eme de racines du groupe\footnotemark\footnotetext{Nous l'avons expliqu\'e dans l'introduction, le cas du bessel $3$ correspond \`a $G_0=\SU(2)$.}.   Dans \cite{defo3} nous consid\'erons l'action coadjointe d'un groupe de lacets \`a valeurs dans un groupe de Lie  compact sur le dual de l'extension centrale de son alg\`ebre de Lie et \'etablissons dans ce cadre une correspondance analogue. 

\begin{figure}
\begin{minipage}[c]{.46\linewidth}
     \begin{center}
             \includegraphics[width=6.1cm]{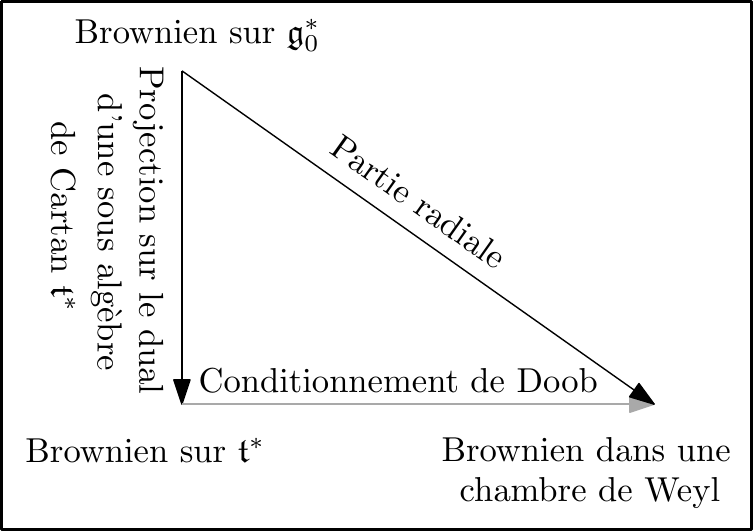}
         \end{center}
         \caption{Diagramme commutatif - Le cas compact}
         \label{D-C-1}
   \end{minipage} \hfill
   \begin{minipage}[c]{.46\linewidth}
    \begin{center}
            \includegraphics[width=6.1cm]{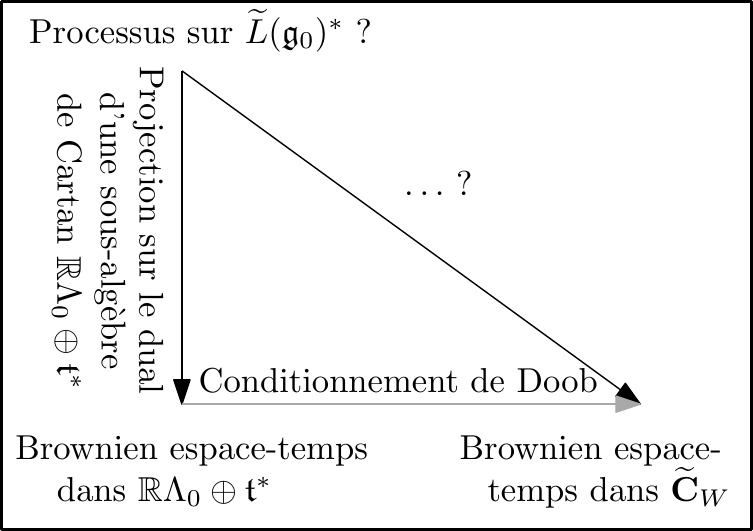}
        \end{center}
            \caption{Diagramme commutatif - Le cas affine }
            \label{D-C-2}
 \end{minipage}
 \end{figure}

 Lorsque nous avons commenc\'e \`a travailler sur le brownien et ses liens avec les repr\'esentations d'alg\`ebres affines, l'un des enjeux \'etait de compl\'eter dans ce contexte un diagramme commutatif analogue \`a celui  repr\'esent\'e en figure \ref{D-C-1}. Comme nous l'avons rappel\'e au chapitre \ref{chap-BETaff}, on construit une approximation discr\`ete du brownien conditionn\'e \`a rester dans une chambre de Weyl associ\'ee au syst\`eme de racines de $\mathfrak g_0$ en consid\'erant des produits tensoriels de repr\'esentations de $\mathfrak g_0$. Dans ce m\^eme chapitre, une construction analogue faisant intervenir des repr\'esentations de plus haut poids d'une alg\`ebre affine a fait apparaitre   un processus limite dans une chambre de Weyl associ\'ee au  syst\`eme de racines d'une telle alg\`ebre.  Cette construction a constitu\'e la premi\`ere \'etape dans l'\'elaboration d'un diagramme et nous a fourni des candidats pour compl\'eter une premi\`ere ligne, comme nous l'indiquons sur la figure \ref{D-C-2} avec les notations du chapitre \ref{chap-loc}. L'ensemble $\boldsymbol{\widetilde C}_W$ que nous d\'efinissons plus loin y est un sous-ensemble de $\R\Lambda_0\oplus\mathfrak t^*$  s'identifiant \`a la chambre de Weyl $\widetilde C_W$ introduite dans le chapitre \ref{chap-Algebre-affine}.  Il restait   \`a construire un processus dont la partie radiale pour l'action d'un  groupe serait le processus conditionn\'e \`a rester dans ce domaine.

C'est dans cette perspective que nous avons entam\'e la lecture du travail de Igor  Frenkel \cite{Frenkel}.  En effet,  on trouve  au c\oe ur du diagramme commutatif de la figure \ref{D-C-1} une formule de Harish-Chandra\footnote{Hans Duistermaat et Gert Heckman \cite{DH} le disent ainsi : "[Harish-Chandra] obtained [his] formula as a consequence of the computation of the radial part of the G-invariant differential operators on $\mathfrak{g}$. Conversely the formula for the radial part can be obtained from [his] formula  using the theory of Fourier integrals." On trouvera aussi   dans \cite{IZ}  une illustration de cette id\'ee pour le  groupe unitaire $\mbox{U}(n)$.},  \'equivalente dans le cas des groupes de Lie compacts \`a  une formule des caract\`eres de Kirillov.   Dans \cite{Frenkel} Igor Frenkel \'etablit une formule des caract\`eres de type Kirillov dans un  contexte infini-dimensionnel en consid\'erant une mesure gaussienne sur le dual d'une alg\`ebre de lacets \`a valeurs dans  une alg\`ebre de Lie simple compacte, essentiellement un mouvement brownien sur cette alg\`ebre. Sa formule est au c\oe ur de la correspondance que nous avons \'etablie.  Une part importante de notre travail a consist\'e \`a la comprendre. Une part importante du pr\'ec\'edent  chapitre \`a consist\'e \`a la pr\'esenter. Nous l'\'enon\c cons  plus rigoureusement maintenant. Nous indiquons ensuite dans la section \ref{drap} notre contribution.

Dans la suite,   $G_0$ est un groupe de Lie compact simple, connexe et simplement connexe, d'alg\`ebre de Lie $\mathfrak g_0$, $T$ est un tore maximal de $G_0$, d'alg\`ebre de Lie $\mathfrak t$. Nous supposons  comme habituellement et sans perte de g\'en\'eralit\'e que  $G_0$ est un groupe de matrices. On note $\Ad$ l'action adjointe de $G_0$ sur $\mathfrak g_0$, c'est-\`a-dire l'action par conjugaison de $G_0$ sur  $\mathfrak g_0$. On munit $\mathfrak g_0$ d'un produit scalaire $\Ad$-invariant not\'e $(\cdot\vert \cdot)$.  On note $\mathfrak g =\mathfrak g_0\oplus i\mathfrak g_0$ l'alg\`ebre complexifi\'ee de $\mathfrak g_0$ et $\mathfrak h$ le complexifi\'e de $\mathfrak t$.
 \section{Partie radiale et mouvement brownien sur $\mathfrak g_0$}
Nous redonnons quelques    d\'efinitions n\'ecessaires \`a  l'\'enonc\'e de la formule de Frenkel et de nos r\'esultats, afin que ceux-ci puissent \^etre expos\'es sans r\'ef\'erences trop fr\'equentes aux chapitres pr\'ec\'edents. 
\paragraph{Orbites coadjointes.} Consid\'erons l'alg\`ebre de lacets $L(\mathfrak g_0)$ de $S^1$ dans $ \mathfrak g_0,$ $S^1$ \'etant identifi\'e \`a $\R/\Z$, sans   pr\'eciser  le degr\'e de r\'egularit\'e des lacets. Sur ce point, nous renvoyons le lecteur \`a \cite{defo3}. On consid\`ere son extension centrale 
 $$\widetilde{L}(\mathfrak g_0)= {L}(\mathfrak g_0)\oplus \R c,$$   d\'efinie au chapitre \ref{chap-loc}, munie du crochet de Lie d\'efini par (\ref{Liehat}). L'int\'egrale orbitale de la formule (\ref{intmu})  du chapitre pr\'ec\'edent porte sur les formes   dans   $\widetilde{L}(\mathfrak g_0)^*$. Nous consid\'erons comme pr\'ec\'edemment le poids fondamental $\Lambda_0$ dans   $\widetilde L(\mathfrak g_0)^*$ d\'efini par 
$$\Lambda_0(L(\mathfrak g_0))=0, \quad \Lambda_0(c)=1,$$ et consid\'erons  l'action coadjointe $\Ad^*$ de $L(G_0)$ sur $ {L}(\mathfrak g_0)^*\oplus \R \Lambda_0$ d\'efinie par (\ref{coadjoint}) en prenant $\theta=0$ et en quotientant par $\R\delta$. On a  pour $\gamma\in L(G_0)$, $\phi\in L(\mathfrak g_0)^*$, $\tau\in \R$,
 \begin{align}\label{coadjointtilde} 
\Ad^*(\gamma)(\phi+\tau \Lambda_0)& =[\Ad_{  0}^*(\gamma)\phi-\tau (\gamma'\gamma^{-1}\vert \cdot)]+\tau\Lambda_0,
\end{align} 
o\`u $\Ad_{ 0}^*(\gamma)\phi(.)=\phi(\gamma^{-1}.\gamma)$. Pour un \'el\'ement  $\zeta$ de $\widetilde L(\mathfrak g_0)^*$ nous notons    $\widetilde{\mathcal O}_\zeta$ l'orbite coadjointe $\Ad^*(L(G_0))\{\zeta\}$ dans $\widetilde L(\mathfrak g_0)^*$ . Nous avons vu dans le chapitre pr\'ec\'edent qu'\'etant donn\'ee une forme $\xi$ de $ \widetilde{\mathcal{O}}_{\zeta}$,  on d\'etermine $\gamma\in L(G_0)$ tel que $\xi=\Ad^*(\gamma)(\zeta)$ en r\'esolvant une \'equation diff\'erentielle, pourvu que le niveau de $\zeta$, c'est-\`a-dire la coordonn\'ee le long de $\Lambda_0$ de $\zeta$,  soit non nulle.  Ainsi   un \'el\'ement $\xi$ s'\'ecrivant 
$$\xi=\tau\Lambda_0+\int_0^1(\cdot\vert\dot{x}_s)\, ds,$$
avec $\dot{x}\in L(\mathfrak g_0)$ et $\tau>0$, est dans    $\widehat{\mathcal O}_{\tau\Lambda_0+(a\vert \cdot)}$ pour $a\in \mathfrak t$ si et seulement si
la solution $\{X_s:s\ge 0\}$ issue de  $I$ de l'\'equation diff\'erentielle
$$\tau d X=Xdx $$ 
 est telle que  $X_1\in \Ad(G_0)\{\exp(a/\tau)\}$.  Les orbites dans $\tau\Lambda_0+L(\mathfrak g_0)^*$ pour l'action coadjointe de $L(G_0)$ sont donc en correspondance avec celles de  $G_0$ pour l'action  par conjugaison de $G_0$ sur lui-m\^eme. Pour param\'etrer celles-ci il est plus naturel de  travailler avec les   racines r\'eelles du groupe $G_0$\footnotemark\footnotetext{\`A ne pas confondre avec les racines r\'eelles de l'alg\`ebre affine... On nous tend des pi\`eges.} plut\^ot qu'avec ses racines infinit\'esimales (voir \cite{Brocker} pour ces notions).

\paragraph{Racines r\'eelles.} Nous l'avons vu, ce sont les racines  r\'eelles  du groupe $G_0$ plut\^ot que ses racines infinit\'esimales qui  apparaissent naturellement dans le cadre des orbites coadjointes d\'ecrites dans le chapitre \ref{chap-loc}.  Elles remplaceront  ici  les racines  infinit\'esimales consid\'er\'ees dans les chapitres \ref{chap-Algebre-affine} et \ref{chap-BETaff}. Nous les noterons   en caract\`ere gras, ainsi que les ensembles associ\'es.  Nous rappelons leur d\'efinition ci-dessous. L'introduction des racines r\'eelles ajoute une difficult\'e \`a l'expos\'e qui  peut para\^itre inutile. J'ai   cherch\'e une solution qui me permette de les \'eviter. Toutes poss\'edaient leur lot d'inconv\'enients et  j'ai  finalement choisi la pr\'esentation qui me semble la moins d\'esagr\'eable. Elle implique  l'utilisation de ces racines qui diff\`erent d'un facteur $2i\pi$ des racines infinit\'esimales\footnote{La solution ne me satisfait pas enti\`erement mais en r\'edigeant ce m\'emoire j'ai  trouv\'e quelque soutien dans les mots de T. Lancester et de S. J. Blundell  s'exprimant ainsi, dans Quantum Field Theory for the gifted amateur, \`a propos d'un f\^acheux facteur $2\pi$ : \guillemotleft We will try to formulate our equations so that every factor of $dk$ comes with a $(2\pi)$, hopefully eliminating one of the major causes of insanity in the subject, the annoying factors of $2\pi$.\guillemotright  \,.}.  Consid\'erons  donc  l'ensemble des racines r\'eelles de $\mathfrak g_0$
$$\boldsymbol\Phi=\{\boldsymbol\alpha\in \mathfrak h^*: \exists X\in \mathfrak{g}\setminus \{0\}, \, \forall H\in \mathfrak h, \, [H,X]=2i\pi \boldsymbol\alpha(H) X\}.$$ On suppose $\mathfrak{g}$ de rang $n$ et on choisit un ensemble de racines  r\'eelles simples $$\boldsymbol\Pi=\{\boldsymbol\alpha_k=\frac{1}{2i\pi}\alpha_k, \, k\in\{1,\dots,n\}\},$$ 
les racines $\{\alpha_k: k\in \{1,\dots,n\}\}$ \'etant les racines infinit\'esimales simples de $\mathfrak g$. On note 
$\boldsymbol\Phi_+$ l'ensemble des racines r\'eelles positives. L'ensemble des coracines simples r\'eelles est
$$\boldsymbol\Pi^\vee=\{\boldsymbol\alpha_k^\vee=2i\pi\alpha^\vee_k,  k\in\{1,\dots,n\}\},$$
o\`u les $\alpha_k^\vee$, $k\in \{1,\dots,n\}$, sont les coracines simples de $\mathfrak g$.  
Pour $\boldsymbol\alpha\in \boldsymbol\Pi$, la r\'eflexion $s_{ \boldsymbol\alpha^\vee}$ est definie sur $\mathfrak{t}$ par $$s_{\boldsymbol\alpha^\vee}(x)=x-\boldsymbol\alpha(x)\boldsymbol\alpha^\vee,\quad  \textrm{ pour } x \in \mathfrak{t}.$$ On  consid\`ere le groupe de Weyl affine \'etendu engendr\'e par les r\'eflexions  $s_{\boldsymbol\alpha^\vee}$ et les translations par $\boldsymbol\alpha^\vee$, $x\in\mathfrak t\mapsto x+ \boldsymbol\alpha^\vee,$ pour $\boldsymbol\alpha\in \boldsymbol\Pi$.  Un domaine  fondamental pour l'action de  ce groupe  sur $\mathfrak{t}$ est  
\begin{align*}
A=\{x\in \mathfrak{t} : \forall \boldsymbol\alpha\in \boldsymbol\Phi_+, \, \, 0\le \boldsymbol\alpha(x)\le  1\}
\end{align*} 
 Le groupe  $G_0$ \'etant  simplement connexe, le noyau de l'application exponentielle $\exp:\mathfrak g_0\to G_0$ est le r\'eseau $\boldsymbol Q^\vee$ engendr\'e par les coracines de $\boldsymbol\Pi^\vee$ et   l'ensemble des classes  de conjugaison  $G_0/\Ad(G_0)$ est en correspondance bijective  avec  le domaine fondamental $A$ \cite{Brocker}. Ainsi pour tout $u\in G_0$,  il existe un unique  $x\in A$ tel que    $u\in \Ad(G_0)\{\exp(x)\}$.  Pour  $\tau\in \R_+$, on d\'efinit  l'alc\^ove $A_\tau$ de niveau  $\tau$ par 
$$A_\tau=\{x\in \mathfrak{t} : \forall \boldsymbol\alpha\in \boldsymbol\Phi_+, \,\, 0\le \boldsymbol\alpha(x)\le \tau\},$$ i.e. $A_\tau=\tau A$.  

\paragraph{Partie radiale.} Soit $\tau$ un r\'eel strictement positif. Nous avons rappel\'e une correspondance bijective entre les orbites dans $\tau\Lambda_0+L(\mathfrak g_0)^*$ et celles dans $G_0$. On d\'efinit ainsi la partie radiale d'un \'el\'ement $\xi$ s'\'ecrivant 
$$\xi=\tau\Lambda_0+\int_0^1(\cdot \vert\dot{x}_s)\, ds,$$
avec $\dot{x}\in L(\mathfrak g_0)$, comme l'unique \'el\'ement $a\in A_\tau$ tel que  $$X_1\in \Ad(G_0)\{\exp(a/\tau)\}$$ o\`u $X$ est la solution $\{X_s:s\ge 0\}$ issue de  $I$ de l'\'equation diff\'erentielle
$$\tau d X=Xdx. $$
Si $x$ est un brownien \`a valeurs dans $\mathfrak g_0$  l'\'equation diff\'erentielle stochastique \begin{align}\label{eds}
\tau \, dX=X\circ dx,
\end{align} 
o\`u $\circ$ d\'esigne une int\'egrale de  Stratonovitch, a une unique solution issue de $I$. Une telle solution est un processus \`a valeurs dans $G_0$. C'est l'exponentielle stochastique de Stratonovitch de $\frac{x}{\tau}$.     On la note $\{\epsilon(\tau,x)_s, s\ge 0\}$. 
En rempla\c cant la trajectoire r\'eguli\`ere par un brownien et l'\'equation diff\'erentielle usuelle par  cette \'equation diff\'erentielle pour la d\'efinition de la partie radiale, on obtient  naturellement la d\'efinition suivante. 
 \begin{defn}  
Pour $\tau\in \R_+^*$, et $x=\{x_s: s\in[0,1]\}$ un mouvement brownien sur $\mathfrak g_0$, on   d\'efinit la partie radiale de $\tau\Lambda_0+\int_0^1(\cdot\vert dx_s)$ qu'on note $\rad(\tau\Lambda_0+\int_0^1(.\vert dx_s))$ par 
$$\rad(\tau\Lambda_0+\int_0^1(\cdot\vert dx_s))=\tau\Lambda_0+(\cdot\vert a),$$  
 o\`u $a$ est l'unique \'el\'ement de $ A_\tau$  tel que $$\epsilon(\tau,x)_1\in \Ad(G_0) \{\exp(a/\tau)\}.$$  \end{defn}
  
   \section{Une formule des caract\`eres de Kirillov--Frenkel}\label{section-Frenkel}
 On note $\boldsymbol\theta$ la plus grande racine r\'eelle, i.e. $\boldsymbol \theta=\theta/(2i\pi)$, o\`u $\theta$ est la plus haute racine de $\mathfrak g$, et on pose $\boldsymbol\alpha_0=\delta-\boldsymbol\theta$. On consid\`ere $\boldsymbol\theta^\vee=2i\pi\theta^\vee$ et on pose $\boldsymbol\alpha_0^\vee=c-\boldsymbol\theta^\vee$. On consid\`ere  $$\boldsymbol{\widehat{\mathfrak h}}=\mbox{Vect}_\C\{\boldsymbol\alpha^\vee_0,\boldsymbol\alpha_1^\vee,\dots,\boldsymbol\alpha^\vee_n,d\} \, \textrm{ et } \,\boldsymbol{ \widehat{\mathfrak h}}^*=\mbox{Vect}_\C\{\boldsymbol\alpha_0,\boldsymbol\alpha_1,\dots,\boldsymbol\alpha_n,\Lambda_0\}.$$ On rappelle que pour $i\in\{0,\dots,n\}$
 $$\boldsymbol\alpha_i(d)=\delta_{i0}, \quad \delta(\alpha_i^\vee)=0, \quad \Lambda_0(\boldsymbol\alpha_i^\vee)=\delta_{i0},\quad \Lambda_0(d)=0.$$ 
 On pose 
 $$\widehat{\boldsymbol\Pi }=\{\boldsymbol\alpha_i: i\in\{0,\dots,n\}\} \, \textrm{ et } \, \widehat{\boldsymbol\Pi }^\vee=\{\boldsymbol\alpha^\vee_i: i\in\{0,\dots,n\}\}.$$ Alors $(\widehat{\boldsymbol{\mathfrak h}},\widehat{\boldsymbol\Pi},\widehat{\boldsymbol\Pi}^\vee)$ est une r\'ealisation de matrice de Cartan g\'en\'eralis\'ee de type affine  et on peut   consid\'erer les objets associ\'es \`a ce triplet de m\^eme que  dans les  chapitres $1$ et $6$ de \cite{Kac}.   Nous les notons en caract\`ere  gras pour marquer la diff\'erence avec ceux associ\'es \`a  un triplet construit avec les racines infinit\'esimales comme dans le chapitre \ref{chap-Algebre-affine}.  Cette r\'ealisation permet de faire co\"incider l'action sur $\boldsymbol{\widehat{\mathfrak h}}^*$ du groupe de Weyl $\widehat{\boldsymbol W}$  d\'efini plus bas avec celle groupe de Weyl affine $\widehat W_0$ dans $L(G_0)$ sur $\boldsymbol{\widehat{\mathfrak h}}^*$.  L'isomorphisme  
\begin{align*}
\boldsymbol\nu:\, \,h \in \mathfrak{t}\mapsto (h\vert\cdot)\in \mathfrak{t}^*
\end{align*} identifie $\mathfrak{t}$ et $\mathfrak{t}^*$.   On note $(\cdot\vert\cdot)$ le produit scalaire sur $\mathfrak t^*$ induit par cet isomorphisme.  \`A partir de maintenant le produit scalaire sur $\mathfrak g_0$ est normalis\'e de telle sorte que $(\boldsymbol\theta\vert\boldsymbol\theta)=2$. Remarquons que  $\boldsymbol\nu(\boldsymbol\theta^\vee)=\boldsymbol\theta$ et $(\boldsymbol\theta^\vee\vert\boldsymbol\theta^\vee)=(\boldsymbol\theta\vert \boldsymbol\theta)=2$.   
On d\'efinit   le groupe de Weyl affine $\widehat{\boldsymbol W}$ comme le sous-groupe de  $\mbox{GL}(\boldsymbol{\widehat{\mathfrak{h}}}^*)$ engendr\'e par les r\'eflexions fondamentales  $s_{\boldsymbol\alpha}$, $\boldsymbol\alpha\in \widehat{\boldsymbol\Pi}$, d\'efinies par $$s_{\boldsymbol\alpha}(\beta)=\beta-  \beta(\boldsymbol\alpha^\vee) \boldsymbol\alpha,\ \quad \beta\in \widehat{\boldsymbol{\mathfrak{h}}}^*.$$   
La forme $(\cdot\vert \cdot)$ est $\widehat{\boldsymbol W}$-invariante.
Le groupe de Weyl affine $\widehat{\boldsymbol W}$  est le produit semi-direct $ {W}\ltimes  {\boldsymbol\Gamma}$, o\`u $  W$ est le groupe de Weyl associ\'e au syst\`eme de racines de $G_0$, et  $\boldsymbol\Gamma$ est le groupe de transformations  $t_{\gamma}$, $\gamma\in \boldsymbol\nu( \boldsymbol Q^\vee)$, d\'efinies par
$$t_\gamma(\lambda)=\lambda+\lambda(c)\gamma-\big[(\lambda\vert \gamma)+\frac{1}{2}(\gamma\vert\gamma)\lambda(c)\big]\delta, \quad \lambda\in \widehat{\boldsymbol{\mathfrak{h}}}^*.$$
  Un domaine fondamental pour l'action de $\widehat{\boldsymbol W}$ sur l'espace quotient $(\R_+\Lambda_0+ \mathfrak t^*+ \R\delta)/\R\delta$ est 
  $$\boldsymbol{  \widetilde C}_W=\{\lambda\in \R\Lambda_0\oplus \mathfrak t^*: \lambda(\boldsymbol\alpha^\vee)\ge 0, \boldsymbol\alpha\in\widehat{ \boldsymbol\Pi}\}. $$ On remarque que par d\'efinition de l'alc\^ove $A_\tau$ un \'el\'ement  $\lambda$ de $\R\Lambda_0+\mathfrak{t}^*$ est  dans $\boldsymbol{ \widetilde C}_W$ si et seulement si il existe  $\tau\in \R_+$ et $a\in A_\tau$ tels que $\lambda=\tau\Lambda_0+\phi_a$, o\`u $\phi_a=(a\vert \cdot)$. 
    Pour $\theta,\tau\in \R_+^*,$ $X\in \mathfrak t$, $a\in A_\tau$, on d\'efinit $\widehat{\boldsymbol\varphi}_{\theta d+X} (\tau\Lambda_0+\phi_a)$ par 
$$\widehat{\boldsymbol\varphi}_{\theta d+X} (\tau\Lambda_0+\phi_a)=\frac{1}{\pi(X/\theta)}\sum_{w\in \widehat{\boldsymbol W}}e^{\langle w(\tau\Lambda_0+\phi_a),\theta d+X\rangle},$$
o\`u $\boldsymbol\pi(X)=\prod_{i=1}^n \sin(\pi \boldsymbol\alpha_i(X))$.
Soit $\{x^{\tau/\theta}_s, s\ge 0\}$ un mouvement brownien sur $\mathfrak g_0$ tel que $\{\sqrt{\theta/\tau}x^{\tau/\theta}_s, s\ge 0\}$  soit un mouvement brownien  standard\footnote{Le $\theta$ ici n'a rien \`a voir avec la racine la plus haute !}. Le th\'eor\`eme suivant     donne la formule des caract\`eres de Kirillov--Frenkel.
\begin{theo}[I.B. Frenkel, \cite{Frenkel}]\label{theo-Frenkel}  Pour $X\in \mathfrak t$, $\theta,\tau\in \R_+^*$, $a\in A_t$,
$$\E\big(e^{(X\vert x_1^{\tau/\theta})}\vert \rad(\tau\Lambda_0+\int_0^1 (\cdot\vert dx_s^{\tau/\theta}))=\tau\Lambda_0+\phi_a\big)= \frac{\widehat{\boldsymbol\varphi}_{\theta d+X}(\tau\Lambda_0+\phi_a)}{\widehat{\boldsymbol\varphi}_{\theta d}(\tau\Lambda_0+\phi_a)}.$$
\end{theo} 
 La mesure image par la projection canonique sur $\R\Lambda_0\oplus \mathfrak t^*$ de la loi de $\tau\Lambda_0+\int_0^1(\cdot\vert x_s^{\tau/\theta})\, ds$ conditionnellement \`a $$\{\rad(\tau\Lambda_0+\int_0^1 (\cdot\vert dx_s^{\tau/\theta}))=\tau\Lambda_0+\phi_a\},$$ peut se comprendre comme une mesure de Duistermaat--Heckman associ\'ee \`a l'action hamiltonienne de $T$ sur l'orbite $\widetilde{\mathcal O}_{\tau\Lambda_0+\phi_a}$  munie de la mesure de Frenkel.  Le th\'eor\`eme en donne la transform\'ee de Fourier (en rempla\c cant $X$ par $iX$). Pour $G_0=\SU(n)$, on obtient, dans une version normalis\'ee, la m\^eme expression que celle donn\'ee  au chapitre \ref{chap-loc} dans la formule (\ref{MDD-DH}).
\section{Partie radiale d'un drap brownien}\label{drap}
\begin{figure}
     \begin{center}
            \includegraphics[scale=0.9]{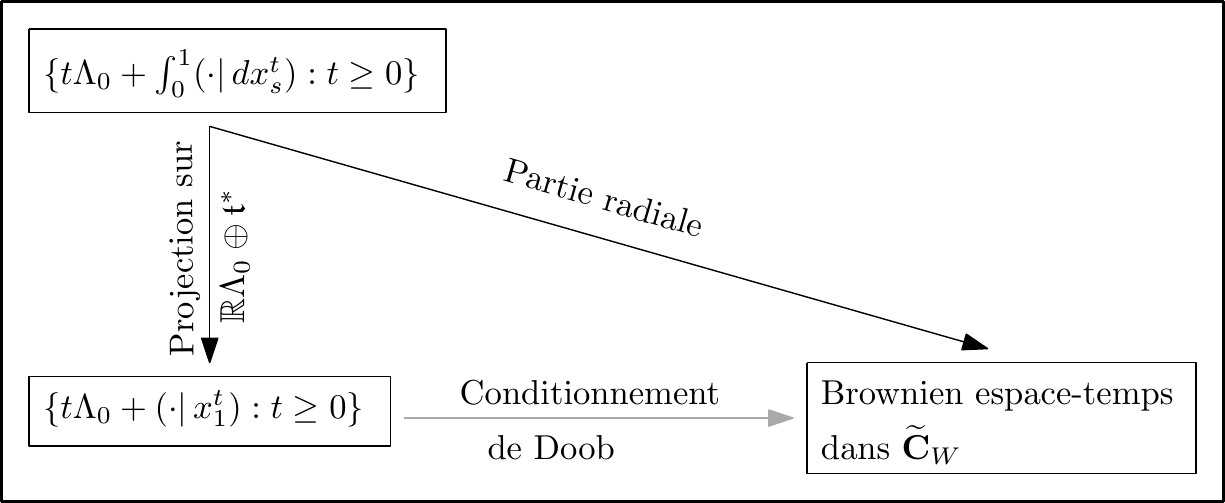}
            \caption{Diagramme commutatif  - Le cas affine}
            \label{DCCA} 
               \end{center}
 \end{figure}

Nous pouvons maintenant   compl\'eter le diagramme de la figure \ref{D-C-2}. Comme nous l'avons expliqu\'e dans l'introduction du m\'emoire, il s'agit de trouver un processus de L\'evy \`a valeurs dans $\widetilde{\mathcal L}(\mathfrak g_0)^*$ dont la projection sur $\R\Lambda_0\oplus \mathfrak t^*$ soit de m\^eme loi  que $\{t\Lambda_0+b_t:t\ge 0\}$ o\`u $\{b_t: t\ge 0\}$ est un brownien standard de $\mathfrak t^*$ et dont la loi en chaque temps fix\'e $t\ge 0$ soit la  m\^eme   que $\tau\Lambda_0+\int_0^1(\cdot\vert x_s^{\tau/\theta})\, ds$ pour des param\`etres $\tau$ et $\theta$ bien choisis. Un tel processus s'obtient \`a partir d'un drap brownien standard $$\{x_s^t : s,t\ge 0\}$$ \`a valeurs dans $\mathfrak g_0$.  Ce drap est un processus gaussien index\'e par $\R_+^2$  tel que pour tout $t\in \R_+^*$, $\{x_s^t/\sqrt{t}:s\ge 0\}$ est un brownien standard \`a valeurs dans $\mathfrak{g}_0$ et pour tout  $t,t'\in \R_+$, $$\{x_s^{t'+t}-x_s^{t}: s\ge 0\}$$ est un processus  ind\'ependant de $\sigma(x_s^r, r\le t, 0\le s)$ et de m\^eme loi que $\{x_s^{t'}: s\ge 0\}$. Il permet de construire un processus qui doit \^etre  pens\'e comme un processus de L\'evy 
\begin{align}\label{Levyaff} 
\{t\Lambda_0+\int_0^1(\cdot\vert \,dx_s^t): t\ge 0\}
\end{align}  \`a valeurs dans $\R\Lambda_0\oplus  {L}(\mathfrak g_0)^*$.  La projection canonique du processus (\ref{Levyaff}) sur $\R\Lambda_0\oplus \mathfrak t^*$  est $$\{t\Lambda_0+(\cdot\vert x_1^t): t\ge 0\},$$ o\`u la forme $(\cdot\vert x_1^t)$ est restreinte \`a  $\mathfrak t$. Si l'on identifie racines r\'eelles et racines infinit\'esimales, c'est le   brownien espace-temps   $\{B_t: t\ge 0\}$ du th\'eor\`eme \ref{TCL-cond}. 
\begin{prop} L'application $\widehat{\boldsymbol\varphi}_d$  est  une fonction harmonique  positive pour le processus $\{t\Lambda_0+(\cdot\vert x_1^t) : t\ge 0\}$ tu\'e sur les bords de $\widetilde{\boldsymbol{  C}}_W$.
\end{prop} 
 Si l'on identifie racines r\'eelles et racines infinit\'esimales,  le  transform\'e de Doob via $\widehat{\boldsymbol\varphi}_d$  du  processus $\{t\Lambda_0+(\cdot\vert x_1^t): t\ge 0\}$   tu\'e sur les bords de $\widetilde{\boldsymbol{  C}}_W$   est le   processus   $\{A_t: t\ge 0\}$ du th\'eor\`eme \ref{TCL-cond}.   Ainsi on a le th\'eor\`eme suivant. 
\begin{theo} [M. D. \cite{defo3}] Soit  $\{x_s^t : s,t\ge 0\}$ un drap brownien standard \`a valeurs dans $\mathfrak g_0$.  Alors
le processus $$\{t\Lambda_0+(\cdot\vert x_1^t): t\ge 0\}$$ est un brownien  espace-temps standard sur $\R\Lambda_0\oplus \mathfrak t^*$, la coordonn\'ee temporelle \'etant le long de $\Lambda_0$. 
Et  le processus $$\{\rad(t\Lambda_0+ \int_0^1(\cdot\vert dx_s^t)): t\ge 0\}$$ est  le  transform\'e de Doob via $\widehat{\boldsymbol\varphi}_d$  du  processus $$\{t\Lambda_0+(\cdot\vert x_1^t) : t\ge 0\}$$   tu\'e sur les bords de $\widetilde{\boldsymbol{  C}}_W$. 
\end{theo}
On peut ainsi compl\'eter  le diagramme commutatif  de la figure \ref{D-C-2}. On donne en figure \ref{DCCA} le diagramme commutatif compl\'et\'e. 
 
\chapter{Transformations de Pitman et brownien espace-temps\\ \vspace{0.5cm}
\small{\it  o\`u l'on file notre Queneau}}\label{chap-pit}
   \epigraph{couleurs de la corde

d\'ep\^ot de cette image

cristaux du temps

traces d'espace}{{\it Desiderata}, Raymond Queneau\footnotemark}
\footnotetext{{Dans L'instant fatal}, 1948.  Il semble  que Raymond Queneau, dont on connaissait  le go\^ut pour   l'analyse matricielle du langage, se soit  aussi  int\'eress\'e aux  alg\`ebres affines et aux cristaux de Kashiwara (plus de quarante ans avant l'invention de ces derniers, on aper\c coit d\'ej\`a la puissance de l'inversion temporelle). On peut en effet proposer la traduction math\'ematique suivante pour les quatre vers cit\'es : $0$ ou $1$/application moment/cristaux de niveau $t$/coordonn\'ee spatiale.}
 
Le th\'eor\`eme de Pitman donne une repr\'esentation du brownien conditionn\'e au sens de Doob \`a rester positif  par une fonctionnelle du mouvement brownien r\'eel. Les mod\`eles de chemins de Littelmann sont des mod\`eles combinatoires pour les repr\'esentations d'alg\`ebres de Lie. La fonctionnelle du th\'eor\`eme de Pitman y  joue un r\^ole important, ce qui offre une voie possible  pour sa preuve (au th\'eor\`eme) :   montrer d'abord une version \`a temps discret impliquant les marches et les processus de Markov introduits dans la section \ref{section-MA-su(2)}, puis  la version \`a temps continu  par une application du th\'eor\`eme central limite. 

Le brownien espace-temps conditionn\'e $\{(t,a_t): t\ge 0\}$ introduit \`a  la section \ref{section-sl2-affine} est la limite d'une suite de  processus de Markov construits \`a partir de produits tensoriels de repr\'esentations d'une alg\`ebre affine de type $A_1^{(1)}$. La th\'eorie de Littelmann est valable pour une telle alg\`ebre et il est naturel de se demander s'il existe une repr\'esentation \`a la Pitman du brownien espace-temps conditionn\'e. Nous avons d\'emontr\'e dans \cite{boubou-defo} qu'une telle repr\'esentation existait et obtenu comme cons\'equence une repr\'esentation de m\^eme type  du brownien conditionn\'e \`a vivre dans un intervalle.    Ph. Biane,  Ph. Bougerol et N. O'Connell ont \'etabli dans \cite{bbobis} une repr\'esentation de type Pitman du brownien conditionn\'e \`a vivre   dans un domaine fondamental pour l'action d'un groupe de Coxeter fini. Le th\'eor\`eme de repr\'esentation  que nous pr\'esentons est  une  variation nouvelle sur le th\'eor\`eme originel faisant intervenir un groupe de Coxeter  infini : le groupe de Weyl associ\'e au syst\`eme de racines de l'alg\`ebre affine $\widehat{\mathcal L}(\mathfrak{sl}_2(\C))$.  La preuve  repose sur   une approximation de ce groupe par des groupes di\'edraux, qui sont des groupes de Coxeter de cardinal fini, pour lesquels nous disposons des r\'esultats de \cite{bbobis}.

Nous pr\'esentons dans la section \ref{Pitmangen} le th\'eor\`eme de Pitman et un exemple de  g\'en\'eralisation pour un groupe  de coxeter de cardinal fini. Nos th\'eor\`emes de repr\'esentation sont \'enonc\'es dans  la section \ref{pitint}.  Ils s'obtiennent en consid\'erant des transformations  du brownien espace-temps $\{(t,b_t): t\ge 0\}$ introduit dans la section \ref{section-sl2-affine}. Nous expliquons dans la section \ref{sectionpitaff}    comment nos r\'esultats  s'envisagent \`a travers le prisme des mod\`eles de chemins de Littelmann.  La section \ref{pit-DH} porte sur  la mesure de Duistermaat--Heckman  associ\'ee \`a l'action hamiltonienne  d'un tore maximal de $\SU(2)$    sur une orbite coadjointe de  $\widetilde L(\mathfrak{su}(2))^*$ munie d'une mesure de Frenkel. Nous y indiquons comment   cette mesure appara\^it lorsque l'on consid\`ere la loi jointe du brownien espace-temps $\{(t,b_t): t\ge 0\}$ et de son image par nos transformations. Enfin nous avons avons \'evoqu\'e dans l'introduction du m\'emoire la pr\'esence d'un terme dit correcteur dans les th\'eor\`emes de repr\'esentation. Nous faisons dans la tr\`es courte section \ref{correc} quelques remarques sur cette correction.
\section{Le th\'eor\`eme de Pitman et ses g\'en\'eralisations}\label{Pitmangen}
Nous pr\'esentons le th\'eor\`eme de Pitman ainsi qu'un th\'eor\`eme de repr\'esentation du brownien dans un c\^one d'angle $\pi/4$ \'etabli dans \cite{bbo}.    Nous esp\'erons ainsi faire apparaitre dans la section suivante notre th\'eor\`eme de repr\'esentation comme une extension naturelle du th\'eoreme de  repr\'esentation de Pitman originel, valable pour un processus dans un domaine fondamental associ\'e au groupe de Coxeter de cardinal infini qu'est le groupe de Weyl associ\'e au syst\`eme de racines de l'alg\`ebre affine $\widehat{\mathcal L}(\mathfrak{sl}_2(\C))$.
 
\subsection{Le th\'eor\`eme de Pitman}
La transformation de Pitman $\mathcal P$  op\`ere sur les chemins \`a valeurs r\'eelles, c'est-\`a-dire les  fonctions continues $f:\R_+\to \R$ telles que $f(0)=0$. Elle est  d\'efinie par
$$\mathcal Pf(t)=f(t)-2\inf_{0\le s\le t}f(s), \, \, t\ge 0.$$ 
Le th\'eor\`eme de Pitman s'\'enonce  ainsi :   si $\{b(t),t\ge 0\}$ est un mouvement brownien r\'eel standard, alors $\{\mathcal Pb(t),t\ge 0\}$ est un processus de Bessel de dimension trois, c'est-\`a-dire un mouvement brownien r\'eel tu\'e en $0$ conditionn\'e au sens de Doob \`a rester positif. Nous l'avons vu dans les chapitres ant\'erieurs, ses densit\'es de transition sont 
$$q_t(x,y)=\frac{h(y)}{h(x)}p_t^0(x,y), \, x,y\in \R^*_+, \, t\ge 0,$$
o\`u  $h$ est la fonction $h(x)=x,$ $x\in \R_+$, et $p_t^0$ est la densit\'e de transition du brownien tu\'e en $0$, qu'un principe de r\'eflexion usuel permet d'\'ecrire 
$$p_t^0(x,y)=p_t(x,y)-p_t(-x,y),\, x,y\in \R^*_+, \, t\ge 0,$$
avec $p_t$ le noyau de la chaleur sur $\R$.  
\begin{figure}[!t]
\centering
\includegraphics[scale=0.8]{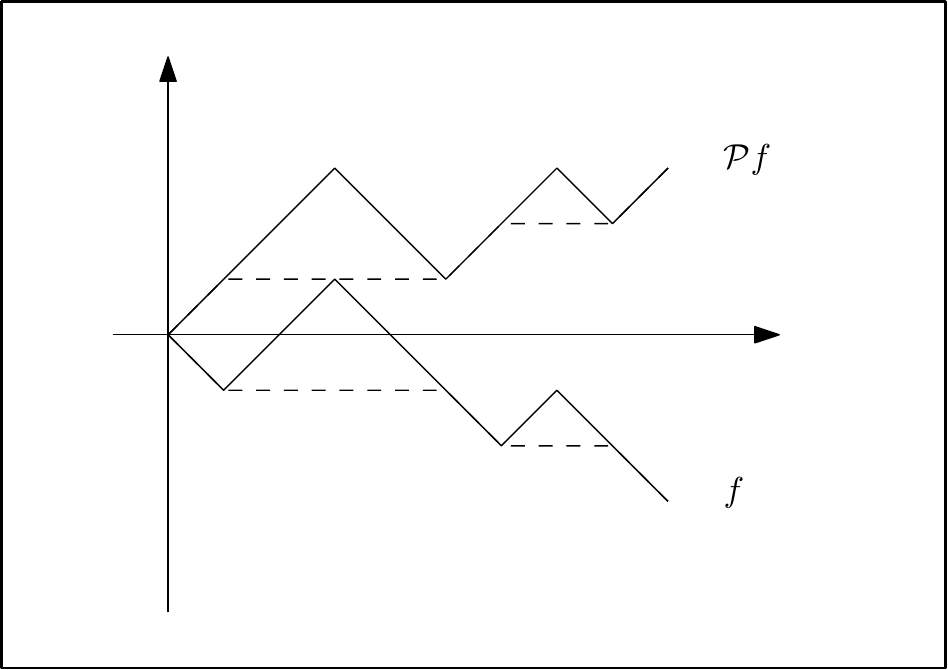}
\caption{Transformation de Pitman} 
\label{Pitman}
\end{figure}
Faisons d\`es \`a pr\'esent quelques remarques simples qui s'av\'ereront, je l'esp\`ere, \'eclairantes. 
\begin{enumerate}
\item Le Bessel de dimension trois vit dans $\R_+$.
\item $\R_+$ est un domaine fondamental pour le groupe de transformations de $\R$   engendr\'e par la sym\'etrie par rapport \`a $0$. 
\item La densit\'e de transition   du brownien tu\'e en $0$ s'exprime comme une somme  altern\'ee  portant sur ce groupe de transformations.  
\item La transform\'ee $\mathcal Pf$ s'obtient en appliquant \`a certaines portions de la courbe repr\'esentative de $f$ une sym\'etrie par rapport \`a l'axe des abscisses (voir la figure \ref{Pitman}).
\item $\mathcal P$ envoie tout chemin sur un chemin \`a valeurs positives. 
\item Les chemins \`a valeurs positives sont laiss\'es invariants par la transformation $\mathcal P$. 
\end{enumerate}
\subsection{Brownien dans un c\^one d'angle $\pi/4$}
\begin{figure}[!t]
\centering
\includegraphics[scale=0.8]{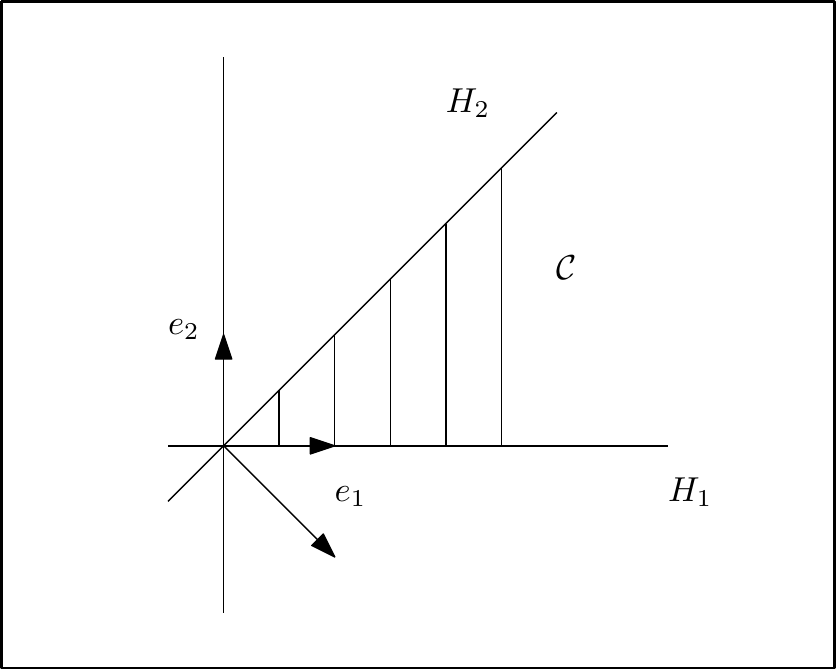}
\caption{Chambre de Weyl de type $BC_2$} 
\label{C_2}
 \end{figure}

   Nous d\'ecrivons ici le cas o\`u le groupe de Coxeter est le groupe  de Weyl associ\'e \`a un syst\`eme de racines de type $B_2$ ou $C_2$, c'est-\`a-dire le groupe de transformations de $\R^2$ engendr\'e par les sym\'etries orthogonales $s_1$ et $s_2$  par rapport aux hyperplans $H_1$ et $H_2$ repr\'esent\'es en figure \ref{C_2}. On note $W(BC_2)$ ce groupe. Les sym\'etries $s_1$ et $s_2$ sont  d\'efinies par 
 $$s_1(x)=x-2(e_2\vert x)e_2, \textrm{ et } s_2(x)=x-(e_1-e_2\vert x)(e_1-e_2), \quad \, x\in \R^2,$$
o\`u $(\cdot \vert \cdot)$ est le produit scalaire usuel sur $\R^2$ et $\{e_1,e_2\}$ la base usuelle de $\R^2$. L'adh\'erence du c\^one  ouvert  $\mathcal C$ d'angle $\pi/4$ 
$$\mathcal C=\{(x_1,x_2)\in \R^2: 0<x_2<x_1\},$$
repr\'esent\'e en figure \ref{C_2}, est un domaine fondamental pour l'action de $W(BC_2)$ sur $\R^2$ : pour tout $x\in \R^2$, il existe un unique vecteur dans l'intersection de $W(BC_2)\cdot\{x\}$ et de     l'adh\'erence de $ \mathcal C$. 

\paragraph{Brownien dans le c\^one.} Consid\'erons l'application  $h$ d\'efinie sur $\R^2$ par $$h(x)=x_1x_2(x_1^2-x_2^2), \, x=(x_1, x_2)\in \R^2.$$ La fonction $h$ est  strictement  positive sur $\mathcal C$ et nulle sur la fronti\`ere de $\mathcal C $.  Elle est harmonique pour le brownien sur $\R^2$ tu\'e sur les bords de $\mathcal C $ et le transform\'e de Doob de  ce processus a pour densit\'es de transition 
$$q_t(x,y)=\frac{h(y)}{h(x)}p^{ \pi/4}_t(x,y), \, \,x,y\in \mathcal C, \, t>0, $$
o\`u $p^{\pi/4}_t$ est la densit\'e de transition du brownien du plan tu\'e sur les bords de $\mathcal C $.  Comme pour le cas du brownien tu\'e en $0$, un principe de r\'eflexion montre que 
$$p^{\pi/4}_t(x,y)=\sum_{w\in W(BC_2)}\det(w)p_t(x,wy),$$
o\`u $p_t$ est le noyau de la chaleur sur $\R^2$.

\begin{figure}[!t]
\centering
\includegraphics[scale=0.8]{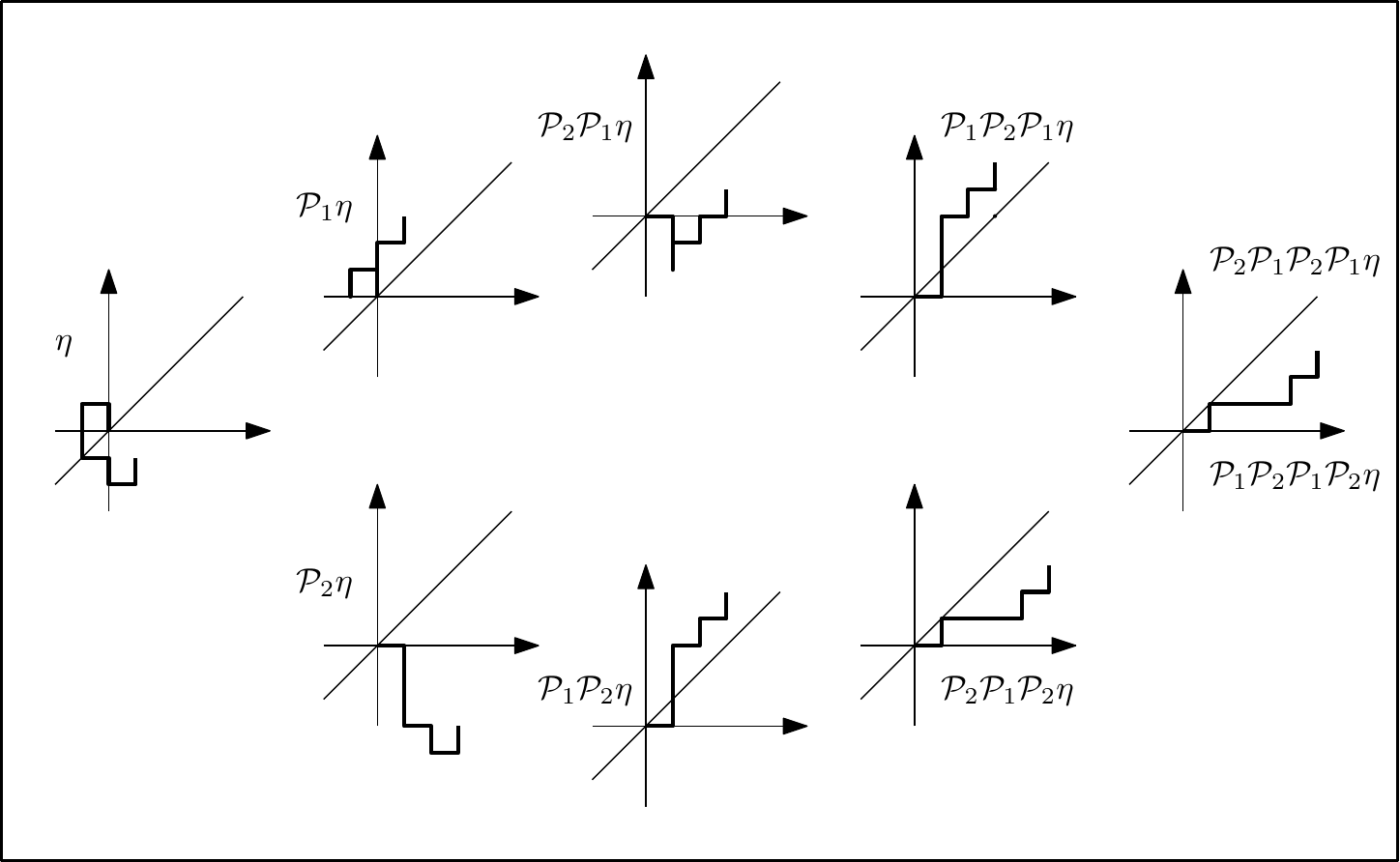}   
\caption{Transformations de Pitman $\mathcal P_1$ et $\mathcal P_2$} 
\label{Pitmans}
\end{figure}

\paragraph{Les transformations de Pitman.} Nous avons vu que le th\'eor\`eme de Pitman   faisait  intervenir une transformation de chemins obtenue \`a partir de la sym\'etrie par rapport \`a $0$. Consid\'erons les transformations de Pitman $\mathcal P_1$ et $\mathcal P_2$ respectivement  associ\'ees aux sym\'etries $s_1$ et $s_2$. Elles agissent sur les chemins \`a valeurs dans $\R^2$, c'est-\`a-dire les applications continues $\eta$ \`a valeurs dans $\R^2$, d\'efinies sur $\R_+$, et telles que  $\eta(0)=0$, de la fa\c con suivante 
\begin{align*}
\mathcal P_1\eta(t)&=\eta(t)-2\inf_{s\le t}(e_2\vert \eta(s))e_2\\
\mathcal P_2\eta(t)&=\eta(t)-\inf_{s\le t}(e_1-e_2\vert \eta(s))(e_1-e_2), \quad t\ge 0,
\end{align*}
La figure \ref{Pitmans} montre les transformations successives d'un chemin $\eta$ trac\'e dans $\R^2$, par les applications $\mathcal P_1$ et $\mathcal P_2$.
Nous remarquons qu'en appliquant  ces transformations successivement et alternativement au chemin $\eta$ on envoie celui-ci en quatre transformations au plus sur un chemin \`a valeurs dans l'adh\'erence du  c\^one  $\mathcal C $. Cette remarque est en fait une propri\'et\'e dont la d\'emonstration non triviale utilise le fait que $W(BC_2)$ poss\`ede un unique plus long \'el\'ement $$w_0=s_1s_2s_1s_2=s_2s_1s_2s_1=-\mbox{Id},$$
qui est de longueur quatre. Ainsi tout vecteur de $\R^2$ est envoy\'e dans l'adh\'erence de $\mathcal C $ par application de   quatre r\'eflexions  $s_1$ et $s_2$ au plus.  Et quatre applications successives des transformations $\mathcal P_1$ et $\mathcal P_2$ envoient tout chemin dans l'adh\'erence de $\mathcal C $. 

\paragraph{Un th\'eor\`eme de repr\'esentation de type  Pitman.} Le th\'eor\`eme de Pitman dans ce cadre s'\'enonce ainsi :  si $\{b(t), t\ge 0\}$ est un brownien standard du plan  alors 
$$\{\mathcal P_1\mathcal P_2\mathcal P_1\mathcal P_2 b(t), t\ge 0\}$$
est le brownien conditionn\'e \`a rester dans $\mathcal C $ d\'efini plus haut.

\section{Transformations de Pitman et brownien espace-temps}\label{pitint} \subsection{Th\'eor\`eme de repr\'esentation}\label{affine}
\begin{figure}[!t]
\centering
\includegraphics[scale=0.8]{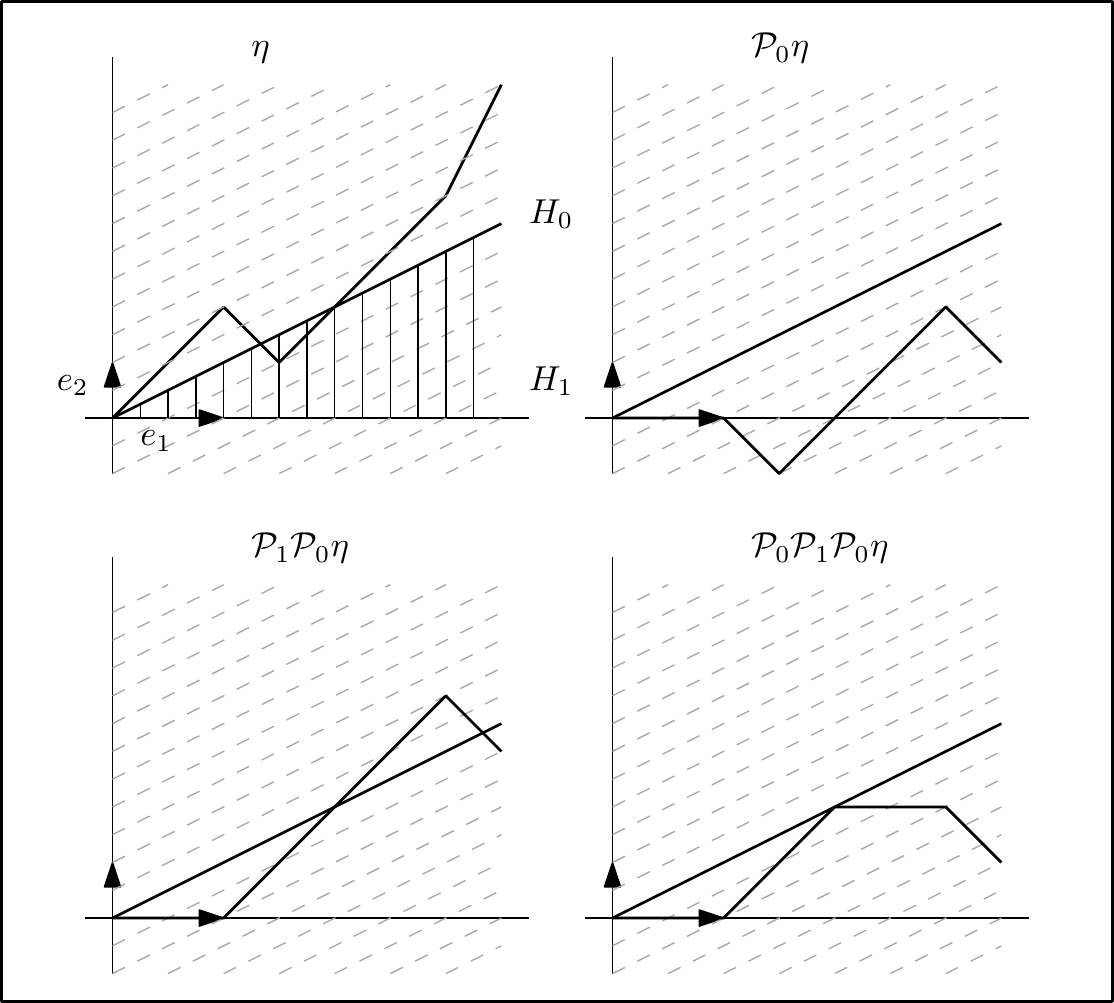}   
\caption{Transformations de Pitman $\mathcal P_0$ et $\mathcal P_1$} 
\label{PitmanAff}
\end{figure}

Consid\'erons un brownien espace-temps r\'eel standard $\{(t,b_t): t\ge 0\}$ et le processus espace-temps $\{(t,a_t): t\ge 0\}$ consid\'er\'e en (\ref{cond-sl2}).  C'est le brownien espace-temps tu\'e sur la fronti\`ere du c\^one 
$$\widetilde C_W=\{(t,x): 0\le x\le t\},$$
conditionn\'e au sens de Doob \`a rester \`a l'int\'erieur du c\^one pass\'e le temps initial. Sur la figure \ref{PitmanAff}, $\widetilde C_W$ est la partie hachur\'ee d\'elimit\'ee par les hyperplans $H_0$ et $H_1$.  Le c\^one $\widetilde C_W$ est un domaine fondamental  pour l'action sur $\R_+\times \R$ d'un groupe $\widetilde W$ engendr\'e par deux sym\'etries $s_0$ et $s_1$ d\'efinies par 
\begin{align*}
s_0(x)&=x+2(e_1-e_2\vert x)e_2=(x_1,-x_2+2x_1),\\
s_1(x)&=x-2(e_2\vert x)e_2,
\end{align*}
$x\in \R^2$. La premi\`ere est la sym\'etrie  par rapport \`a l'hyperplan $H_0$ repr\'esent\'e en  figure \ref{PitmanAff}, parall\`element \`a $e_2$. La seconde est la  sym\'etrie orthogonale par rapport \`a l'hyperplan $H_1$.  Remarquons que les deux sym\'etries laissent invariante la   coordonn\'ee le long de $e_1$.  
\paragraph{Les transformations de Pitman.} Consid\'erons les transformations de Pitman $\mathcal P_{s_0}$ et $\mathcal P_{s_1}$ correspondant aux sym\'etries $s_0$ et $s_1$. Elles  agissent sur un chemin (espace-temps) $\eta(t)=(t,f(t)), t\in \R_+$, o\`u $f(t)\in \R$, de la fa\c con suivante.
\begin{align*}
\mathcal P_{s_0}\eta(t)&=(t, f(t)+2\inf_{s\le t} (s-f(s)))\\
\mathcal P_{s_1}\eta(t)&=(t, f(t)-2\inf_{s\le t} f(s)), \quad t\ge 0.
\end{align*} 
La figure \ref{PitmanAff}  montre les transformations successives d'un chemin $\eta$ par $\mathcal P_{s_0}$ et $\mathcal P_{s_1}$. 
\paragraph{Applications successives de transformations de Pitman.} Le groupe $\widetilde W$ est de cardinal infini et n'a pas de plus long \'el\'ement. On peut montrer que si $f$ est $C^1$ par morceaux sur $[0,T]$, alors il existe $k$ tel que pour $t\in[0,T]$
$$\mathcal P_{s_k}\dots \mathcal P_{s_1}\mathcal P_{s_0}\eta(t)\in \widetilde{ C}_W  ,$$
o\`u $\mathcal P_{s_{2n}}=\mathcal P_{s_{0}}$ et $\mathcal P_{s_{2n+1}}=\mathcal P_{s_{1}}$ . Dans ce cas, pour $n\ge k$, 
$$\mathcal P_{s_n}\dots \mathcal P_{s_1}\mathcal P_{s_0}\eta(t)=\mathcal P_{s_k}\dots \mathcal P_{s_1}\mathcal P_{s_0}\eta(t), \,\, t\in[0,T].$$
\paragraph{Transformations de Pitman et trajectoires browniennes.} Il est naturel de se demander ce qu'il advient  quand on remplace $\eta$ par un brownien espace-temps $\{B(t)=(t,b_t), t\ge 0\}$. En fait, dans ce cas, la limite de $\mathcal P_{s_n}\dots\mathcal P_{s_0}B(t)$  quand $n$ tend vers l'infini n'existe pas.  On peut montrer que 
\begin{align}\label{deux}
\lim_{n\to\infty}\vert\vert \mathcal P_{s_{n+1}}\dots \mathcal P_{s_0} B(t)-\mathcal P_{s_{n}}\dots \mathcal P_{s_0} B(t)\vert\vert=2,
\end{align}
d\`es que $t>0$. Cependant, en introduisant les transformations $\mathcal L_{s_0}$ et $\mathcal L_{s_1}$ d\'efinies par 
\begin{align*}
\mathcal L_{s_0}\eta(t)= (t, f(t)+\inf_{s\le t} (s-f(s))),\quad 
\mathcal L_{s_1}\eta(t)= (t, f(t)-\inf_{s\le t} f(s)),
\end{align*}
pour $\eta(t)=(t,f(t)), t\ge 0$, 
on a le r\'esultat de convergence suivant, en posant $\mathcal L_{s_{2n}} =\mathcal L_{s_{0}}$ et  $\mathcal L_{s_{2n+1}} =\mathcal L_{s_{1}}$.  
\begin{theo}[Philippe Bougerol,  M.D, \cite{boubou-defo}] \label{Pitaff} Les suites de  processus
$$\{\mathcal L_{s_{n+1}}\mathcal P_{s_n}\dots\mathcal P_{s_0}B(t), t\ge 0\} \textrm{ et } \{\mathcal L_{s_{n+1}}\mathcal P_{s_n}\dots\mathcal P_{s_1}B(t), t\ge 0\}, \quad n\ge 0,$$ convergent en loi vers le processus conditionn\'e  $\{(t,a_t), t\ge 0\}$.
\end{theo}
 La preuve de ce th\'eor\`eme repose sur une approximation du brownien espace-temps par un brownien du plan, et du processus conditionn\'e $\{(t,a_t):t\ge 0\}$ par un brownien du plan conditionn\'e \`a rester dans un c\^one de $\R^2$ d'angle ${\pi}/{m}$,   pour lequel nous disposons des r\'esultats de Biane et al.\ \cite{bbo}. Enfin une inversion du temps permet d'obtenir une repr\'esentation  de type Pitman du brownien dans l'intervalle. Elle est \'enonc\'ee dans le th\'eor\`eme \ref{pit01}.  

\subsection{Le brownien dans l'intervalle}
Le brownien conditionn\'e au sens de Doob \`a rester dans $(0,1)$ est un processus de Markov issu de $0$,  et \`a valeurs dans l'intervalle $(0,1)$, pass\'e le temps initial. Ses densit\'es  de transition sont
$$q_t(x,y)=\frac{\sin(\pi y)}{\sin(\pi x)}e^{\pi^2 t/2}p_t^{0,1}(x,y), \, x,y\in (0,1),$$
o\`u $p_t^{0,1}$ est la densit\'e de transition  du brownien tu\'e en $0$ et en $1$. Un principe de r\'eflexion permet d'\'ecrire   $p_t^{0,1}$  comme une somme altern\'ee de noyau de la chaleur $p_t$  sur $\R$. On a  ainsi
\begin{align}\label{somalt01}
p_t^{0,1}(x,y)=\sum_{k\in \Z}p_t(x,y+2k)-p_t(x,-y+2k), \, x,y\in (0,1).
\end{align}
On remarque que la somme porte cette fois sur le groupe engendr\'e par la sym\'etrie par rapport \`a $0$ et les translations de $+2$ ou $-2$, qui est aussi le groupe engendr\'e par la sym\'etrie par rapport \`a $0$ et la sym\'etrie par rapport \`a $1$. Par ailleurs, l'intervalle $[0,1]$ est un domaine fondamental pour ce groupe de transformations.  Ici, le groupe pr\'esente une diff\'erence essentielle avec les groupes de transformations  pr\'ec\'edents : ses transformations ne sont pas lin\'eaires.  Travailler dans la chambre affine plut\^ot que dans l'intervalle aura donc \'et\'e une mani\`ere de lin\'eariser le probl\`eme.   Pour d\'eduire du th\'eor\`eme de repr\'esentation \ref{Pitaff} un th\'eor\`eme de repr\'esentation du brownien dans l'intervalle,     on remarque que pour une trajectoire issue de z\'ero $\{f(t), t\ge 0\}$ \`a valeurs r\'eelles, on a l'\'equivalence 
$$(\forall t>0, \, 0\le f(t)\le 1)\Leftrightarrow (\forall t>0, \, 0\le tf(1/t)\le t).$$
La loi du brownien est invariante par la transformation $$\{f(t),t\ge  0\}\mapsto \{tf(1/t), t\ge 0\},$$ et on peut  montrer   la proposition suivante. 
\begin{prop} 
Si $\{Z_t,t\ge 0\}$ est un Brownien issu de $0$, conditionn\'e \`a rester dans $(0,1)$, alors $\{(t,tZ_{1/t}), t\ge 0\}$ a la m\^eme loi que $\{(t,a_t):t\ge 0\}$. \end{prop} 
Nous pouvons alors \'enoncer le th\'eor\`eme de repr\'esentation du brownien dans l'intervalle de \cite{boubou-defo}.  Il implique les transformations de chemins $\mathcal L_1$ et $\mathcal P_1$ d\'efinies par  
\begin{align*}\mathcal L_1\varphi(t)=\varphi(t)-\inf_{0\leq s\le t}\varphi(s),\quad 
\PP_1\varphi(t)=\varphi(t)-2\inf_{0\leq s\le t}\varphi(s), \quad t\ge 0,\end{align*}
pour $\varphi:\R_+\to \R,$ tel que $\varphi(0)=0$. Ce sont  les transformations de L\'evy et de Pitman usuelles. Nous introduisons les transformations $\mathcal L_0$ et $\mathcal P_0$ de $\varphi$ en posant 
\begin{align*}
\mathcal L_0\varphi(t)&=\varphi(t)+\inf_{0\leq s\le t}(s-\varphi(s)),\quad 
\PP_0\varphi(t)=\varphi(t)+2\inf_{0\leq s\le t}(s-\varphi(s)), \quad t\ge 0. \end{align*}
Alors, si $\{b_t: t\ge 0\}$ est un brownien r\'eel standard, en notant pour $n\ge 0$, $\mathcal P_{2n}=\mathcal P_0,\mathcal L_{2n}=\mathcal L_0$ et $\mathcal P_{2n+1}=\mathcal P_1,\mathcal L_{2n+1}=\mathcal L_1$, on a le th\'eor\`eme suivant.
\begin{theo}[Philippe Bougerol, M.D, \cite{boubou-defo}]\label{pit01} Les suites de processus 
$$\{ t\mathcal L_{n+1} \PP_{n}\cdots \PP_1 \PP_0 b(1/t),\,  t\ge 0 \} \textrm{ et } \{ t\mathcal L_{n+1} \PP_{n}\cdots \PP_2\PP_1 b(1/t),\,  t\ge 0 \} $$
convergent en loi, quand $n$ tend vers l'infini, vers le mouvement brownien issu de $0$, conditionn\'e \`a rester dans l'intervalle $(0,1)$. 
\end{theo}

  \section{Mod\`eles de chemins de Littelmann}\label{sectionpitaff}
Les mod\`eles de chemins de Littelmann sont des mod\`eles combinatoires d\'evelopp\'es dans le cadre de la th\'eorie des repr\'esentations. Le th\'eor\`eme \ref{Pitaff} dont d\'ecoule le th\'eor\`eme \ref{pit01} fait intervenir un brownien conditionn\'e \`a rester dans un domaine fondamental pour l'action du groupe de Weyl associ\'e au syst\`eme de racines de l'alg\`ebre affine $\widehat{\mathcal L}(\mathfrak{sl}_2(\C))$. Ce brownien  est apparu dans le chapitre \ref{chap-BETaff} comme limite d'un processus de Markov sur les plus hauts poids   de repr\'esentations de $\widehat{\mathcal L}(\mathfrak{sl}_2(\C))$. Les transformations $\mathcal P_0$ et $\mathcal P_1$ jouent un r\^ole important dans le contexte des mod\`eles de chemins associ\'es \`a ces repr\'esentations. Pour autant,  notre d\'emonstration du th\'eor\`eme \ref{Pitaff} n'utilise aucunement  les mod\`eles de Littelmann.  Il nous semble cependant int\'eressant pour plusieurs raisons d'expliquer comment ce th\'eor\`eme peut s'envisager dans ce cadre. D'abord, cela nous permettra de l'int\'egrer   \`a un \'edifice que nous esp\'erons coh\'erent dans lequel sont explor\'es les liens entre propri\'et\'es du brownien et th\'eorie des repr\'esentations. Ensuite, ces mod\`eles \'etant valables pour une alg\`ebre de Kac--Moody quelconque, ils ouvrent   la voie \`a un th\'eor\`eme de repr\'esentation du   brownien dans un domaine fondamental associ\'e \`a une  alg\`ebre de Kac--Moody affine autre que $\widehat{\mathcal L}(\mathfrak{sl}_2(\C))$. 
 Enfin, l'approche par les mod\`eles de chemins \'etait la n\^otre   quand nous avons commenc\'e \`a r\'efl\'echir \`a la question d'une repr\'esentation de type Pitman dans un cadre affine. Malgr\'e nos efforts, certains points techniques n'ont pu \^etre surmont\'es. Pourtant si ce point de vue  ne nous a pas permis d'aboutir au th\'eor\`eme de repr\'esentation final, il a largement contribu\'e,  par l'intuition qu'il nous a donn\'e, \`a son \'elaboration.  

 \subsection{Le cas de $\mathfrak{sl}_2(\C)$} Nous indiquons  d'abord ce qu'est le mod\`ele des chemins de Littelmann  \cite{littel} pour l'alg\`ebre de Lie $\mathfrak{sl}_2(\C)$  et rappelons son lien avec le th\'eor\`eme de repr\'esentation de Pitman, tel qu'il a \'et\'e mis au jour dans  \cite{bbo} et approfondi dans  \cite{bbobis}.
  \paragraph{Module de Littelmann.} Consid\'erons $V=\R\alpha_1$, o\`u $\alpha_1$ est la racine positive de $\mathfrak{sl}_2(\C)$. Un chemin $\pi$ \`a valeurs dans $V$ est une application continue $\pi:\R_+\to V$ telle que $\pi(0)=0$. Il peut s'\'ecrire, pour $s\ge 0$, $\pi(s)=f(s)\alpha_1$ o\`u $f(s)\in \R$. Un chemin dominant est un chemin \`a valeurs dans une chambre de Weyl,  qui est ici $\R_+\alpha_1$. Ainsi $\pi$ est dominant si $f$ est \`a valeurs positives. On fixe $t>0$. Un chemin int\'egral sur $[0,t]$ est un chemin lin\'eaire par morceaux tel que $2f(t)$ et $2\min_{s\le t}f(s)$ soient dans $\Z$. Pour un module irr\'eductible de plus haut poids  $\omega$ de $\mathfrak{sl}_2(\C)$, si on choisit un   chemin  dominant $\pi$ int\'egral sur $[0,t]$ tel que $\pi(t)=\omega$, alors l'ensemble $\mathbb B_{\pi}$ des chemins int\'egraux   $\eta$  sur $[0,t]$ tels que $\mathcal P_{\alpha_1}\eta$ vaille $\pi$, avec $$\mathcal P_{\alpha_1}\eta(s)=(\varphi(s)-2\inf_{0\le u\le s}\varphi(u))\alpha_1,\quad s\in[0,t],$$
 o\`u $\eta(s)=\varphi(s)\alpha_1$, est un module de Littelmann associ\'e au module de plus haut poids $\omega$ (voir P. Littelmann \cite{littel}). On reconna\^it dans $\mathcal P_{\alpha_1}$ la transformation de Pitman. 
 \paragraph{Th\'eor\`eme de Pitman.} Cette approche fournit une preuve du th\'eor\`eme de repr\'esentation de Pitman. Consid\'erons en effet un plus haut poids $\omega$ de $\mathfrak{sl}_2(\C)$, c'est-\`a-dire un \'el\'ement de $\N\alpha_1/2$,  et un chemin dominant $ \pi$ int\'egral  sur $[0,1]$ tel que $\pi(1)=\omega$. On consid\`ere une suite $(\eta_i)_{i\ge 1}$ de chemins al\'eatoires ind\'ependants identiquement distribu\'es selon  une mesure de probabilit\'e uniforme sur le module de Littelmann $\mathbb B_\pi$ et on pose pour $t\in [k-1,k]$, $k\ge 1$,
 $$\eta(t)=\eta_1*\dots *\eta_k(t),$$ 
o\`u $*$ est la concat\'enation usuelle de chemins. Notons que lorsque $\omega=\alpha_1/2$, $\{\eta(k): k\ge 0\}$ est en identifiant $n$ \`a $n\alpha_1/2$ la marche simple sur $\Z$. La th\'eorie de Littelmann assure que le processus
$$\{\mathcal P_{\alpha_1}\eta(k), k\ge 0\}$$ est un processus de Markov \`a valeurs dans l'ensemble des poids dominants de $\mathfrak{sl}_2(\C)$ dont le noyau de transition, si l'on  identifie $n$ \`a $n\alpha_1/2$, est donn\'e par (\ref{MarkovSU2}), avec $q=1$. Par la m\^eme identification, c'est la marche simple conditionn\'ee \`a rester positive lorsque $\omega=\alpha_1/2$. Le processus $\{\frac{1}{\sqrt{n}}\eta(nt), t\ge 0\}$ converge vers un brownien de $V$ quand $n$ tend vers l'infini et le processus $\{\frac{1}{\sqrt{n}}\mathcal{P}_{\alpha_1}\eta(nt), t\ge 0\}$  converge vers un processus de  Bessel de dimension trois. La transformation de Pitman commutant avec la multiplication par un scalaire, on obtient ainsi le th\'eor\`eme de repr\'esentation de Pitman.
 
 \subsection{Le cas de $\widehat{\mathcal L}(\mathfrak{sl}_2(\C))$} Expliquons maintenant ce que donne cette approche lorsque l'on remplace $\mathfrak{sl}_2(\C)$ par $\widehat{\mathcal L}(\mathfrak{sl}_2(\C))$, c'est-\`a-dire  $\mathfrak{sl}_2(\C)$ par  une alg\`ebre affine de type $A_1^{(1)}$.
Une r\'ealisation d'une sous-alg\`ebre de Cartan $\widehat{\mathfrak h}_\R$ est ici donn\'ee par $\R^3$. On identifie $\widehat{\mathfrak h}_\R^*$ \`a $\R^3$ et on note 
\begin{align*}
c=(1,0,0), \, \alpha_1^\vee=(0,1,0), \, d=(0,0,1)\\
\Lambda_0=(1,0,0), \, \alpha_1=(0,2,0), \, \delta=(0,0,1).
\end{align*} 
 On pose $ \alpha^\vee_0=c- \alpha^\vee_1$ et $\alpha_0=\delta-\alpha_1$.
Le groupe de Weyl est engendr\'e par les r\'eflexions  $s_{\alpha_0}$ et $s_{\alpha_1}$ d\'efinies sur $\widehat{\mathfrak{h}}_\R^*$ par  
$$s_{\alpha_i}(v)=v-  \langle v,\alpha^\vee_i\rangle\alpha_i,$$
$v \in \widehat{\mathfrak{h}}_\R^*$ et $i\in\{0,1\}$. Ce sont les r\'eflexions par rapport aux murs de la chambre de Weyl 
$$\widehat C_W=\{t\Lambda_0+x\alpha_1/2+y \delta, (t,x)\in \widetilde C_W, y \in \R\}= \widetilde C_W\times \R.$$ L'action du groupe de Weyl sur l'espace quotient\'e par $\R\delta$ s'identifie \`a celle du groupe de transformation sur $\R^2$ consid\'er\'e dans la section \ref{affine}.
 \paragraph{Module de Littelmann.} Dans la th\'eorie de Littelmann, un chemin est maintenant une application continue et affine par morceaux  $\eta:\R^+\to \widehat{\mathfrak{h}}_\R^*$ telle que     $\eta(0)=0$.   On d\'efinit les transformations de chemins  $\mathcal P_{\alpha_i}$, $i\in\{0,1\}$, par $$\mathcal P_{\alpha_i}\eta(t)=\eta(t)- \inf_{0 \leq s\le t} \langle \eta(s),\alpha^\vee_i \rangle\alpha_i.$$ 
 Remarquons que dans l'espace quotient\'e  $\widehat{\mathfrak h}_\R^*/\R\delta$ identifi\'e \`a $\R^2$ on a  
 $\mathcal P_{\alpha_i}\eta(t)=\mathcal P_{s_i}\eta(t),$ pour $\eta(s)=(s,f(s))$, $s\ge 0$.
Un chemin dominant est un chemin \`a valeurs dans    $ \widehat C_W$ et un chemin  int\'egral   sur $[0,t]$ est un chemin $\eta$ tel que 
$$\langle \eta(t),\alpha^\vee_i \rangle,\, \min_{0\le s\le t}\langle \eta(s),\alpha^\vee_i \rangle\in \Z,\, i=0,1.$$ Pour  un r\'eel $t >0$ et un chemin dominant $\pi$ int\'egral  sur  $[0,t]$, on d\'efinit  le module de Littelmann $\mathbb B_\pi$ engendr\'e par $\pi$ comme \'etant   l'ensemble des chemins  int\'egraux $\eta $ sur $[0,t]$ tels qu'il existe   $n\in \N$ tel que  $$\mathcal P_{\alpha_n}\mathcal P_{\alpha_{n-1}}\cdots \mathcal P_{\alpha_1}\mathcal P_{\alpha_0}\eta(s) = \pi(s)$$ pour $0 \leq s \leq t,$ o\`u   $\alpha_{2k}=\alpha_{0}$ et $\alpha_{2k+1}=\alpha_1$ pour $k\in \N$. Ce module donne une description  du cristal de Kashiwara de plus haut poids  $\pi(t)$ qui est un objet combinatoire   param\'etrant les poids de la repr\'esentation de plus haut poids $\pi(t)$ de $\widehat{\mathcal{L}}(\mathfrak{sl}_2(\C))$ (voir M. Kashiwara \cite{kash95}). Pour un chemin int\'egral   $\eta$ sur $[0,t]$, il existe   $k\in \N$ tel que  pour tout $n\ge k$, 
$$\mathcal P_{\alpha_n}\cdots\mathcal P_{\alpha_0}\eta(s)=\mathcal P_{\alpha_k}\cdots\mathcal P_{\alpha_0}\eta(s),$$ 
pour  $0 \leq s \leq t$. On note  alors
$$\widehat{\mathcal P}\eta(s)=\mathcal P_{\alpha_k}\cdots\mathcal P_{\alpha_0}\eta(s).$$
Le chemin $\widehat{\mathcal P}\eta$ est dominant.

\paragraph{Marches al\'eatoires et cha\^ines de Markov sur le r\'eseau des poids.} Nous avons dit que la preuve du th\'eor\`eme  \ref{pit01} reposait sur une approximation  du brownien espace-temps conditionn\'e $\{(t,a_t): t\ge 0\}$ par un brownien du plan conditionn\'e \`a rester dans un c\^one de $\R^2$ d'angle ${\pi}/{m}$. Nous avons vu dans le chapitre \ref{chap-BETaff} une autre approximation de ce processus espace-temps. Indiquons son lien avec le mod\`ele des chemins.  Consid\'erons  un chemin dominant $ \pi$ int\'egral  sur $[0,1]$ tel que $\pi(1)=\Lambda_0$ et  pour tout $n\ge 0$,  une suite $(\eta_i)_{i\ge 1}$ de chemins al\'eatoires ind\'ependants, identiquement distribu\'es, de loi
$$\P(\eta_1=\zeta)=\mu_{\Lambda_0}(\zeta(1)), \quad \zeta\in \mathbb B_\pi,$$
o\`u  $\mu_{\Lambda_0}$ (qui d\'epend de $n$) est  donn\'ee dans la section \ref{section-sl2-affine}. On construit comme pr\'ec\'edemment un processus al\'eatoire $\{\eta(t), t\ge 0\}$ par concat\'enation de ces chemins. La th\'eorie de Littelmann montre que $$\{\widehat{\mathcal P}\eta(k): k\ge 0\}$$ est le processus de Markov de noyau de transition $Q_{\Lambda_0}$ d\'efini dans   la section \ref{section-sl2-affine}. Nous pouvons  d\'eduire des r\'esultats de cette section que la suite de  processus $$\{\frac{1}{n}\eta([nt]) : t\ge 0\}$$ converge  dans l'espace quotient $\widehat{\mathfrak{h}}^*_\R/\R\delta$ vers un brownien espace-temps, la coordonn\'ee temporelle \'etant le long de $\Lambda_0$ et la coordonn\'ee spatiale le long de $\alpha_1/2$, et que la suite $$\{\frac{1}{n}\widehat{\mathcal P}\eta([nt]) : t\ge 0\}$$  converge dans l'espace quotient vers $\{(t,a_t):t\ge 0\}$.  La convergence (\ref{deux}) montre que dans l'espace quotient $\widehat{\mathfrak{h}}^*_\R/\R\delta$ l'interversion des limites dans 
$$\lim_{n\to\infty}\lim_{k\to\infty}\mathcal{P}_{\alpha_{k+1}}\dots\mathcal P_{\alpha_0}(\frac{1}{n}\eta([nt])$$
n'est pas valable. Il serait cependant int\'eressant d'\'etablir une preuve du th\'eor\`eme \ref{Pitaff} empruntant la voie des chemins de Littelmann. Cela serait en effet une premi\`ere \'etape importante pour l'obtention d'un th\'eor\`eme de Pitman valable pour toute alg\`ebre de Kac--Moody affine.
\section{Cristaux de Kashiwara et mesure de  Duistermaat--Heckman}\label{pit-DH}
Consid\'erons le groupe $G_0$ et son alg\`ebre de Lie $\mathfrak g_0$ consid\'er\'es au chapitre pr\'ec\'edent. Nous avons d\'efini dans un cadre affine la mesure de Duistermaat--Heckman comme la mesure image par la projection sur $\R\Lambda_0\oplus \mathfrak t^*$ d'une  mesure de Frenkel sur une orbite coadjointe de $\widetilde{L}(\mathfrak g_0)^*$. Dans un cadre   compact   la mesure de Duistermaat--Heckman associ\'ee \`a l'action d'un tore maximal de $G_0$ sur une orbite coadjointe de $\mathfrak g_0^*$ peut   s'obtenir par approximation semi-classique des poids d'une repr\'esentation irr\'eductible de $\mathfrak g_0$ dont le plus haut poids est grand. Les cristaux de Kashiwara sont des objets combinatoires qui param\`etrent les poids des repr\'esentations irr\'eductibles d'une alg\`ebre de Lie et on peut ainsi d\'ecrire la mesure de Duistermaat--Heckman comme la mesure image par une fonction \og poids \fg \, d'une mesure d\'efinie sur une version continue des cristaux de Kashiwara.  Dans un cadre affine, comme observ\'e dans \cite{boubou-defo}, la mesure de Duistermaat--Heckman apparaissant aux chapitres \ref{chap-loc} et \ref{chap-drap} s'obtient de m\^eme par approximation semi-classique et se d\'ecrit de m\^eme comme une mesure image par une fonction \og poids\fg \, d'une certaine mesure  sur un analogue continu d'un cristal de Kashiwara. Pour le type $A_1^{(1)}$ notre mod\`ele brownien fournit une construction possible de la mesure \`a consid\'erer sur un tel cristal. Nous expliquons comment dans cette section. 
\paragraph{Coordonn\'ees en corde.} Rappelons  dans le contexte de la section pr\'ec\'edente, ce que sont les coordonn\'ees en corde d'un chemin int\'egral. On   associe \`a un chemin  $\eta$ int\'egral sur $[0,t]$, une suite d'entiers $(x_k)_{k\ge 0}$ nulle \`a partir d'un certain rang, d\'efinie par 
$$\mathcal P_{\alpha_n}\cdots\mathcal P_{\alpha_0}\eta(t)=\eta(t)+\sum_{k=0}^nx_k\alpha_k.$$ Les $x_k, k\ge 0,$ s'appellent les coordonn\'ees en corde de $\eta$ sur $[0,t]$.   
\paragraph{Cristaux.} Pour  $\lambda$ un poids dominant, i.e. $\lambda$ un \'el\'ement de $\widehat C_W$ \`a coordonn\'ees enti\`eres, on pose
$$B(\lambda)=\{x\in B(\infty) : x_k\le  \langle \lambda-\sum_{i=k+1}^\infty x_i\alpha_i, \alpha^\vee_k\rangle \textrm{ pour tout }  k\ge 0\},$$
o\`u
$$B(\infty)=\{x=(x_k)\in  \N^{(\N)} : \frac{x_k}{k}\ge \frac{x_{k+1}}{k+1}\ge 0,  \textrm{ pour tout } k\ge 1\}. $$
Les ensembles $B(\infty)$ et $B(\lambda)$ sont des param\'etrisations des cristaux de Kashiwara, respectivement associ\'es au module de Verma et au module de plus haut poids $\lambda$. Si $\pi$ est un chemin dominant int\'egral sur $[0,t]$ tel que $\pi(t)=\lambda$, alors $\mathbb B_\pi$ et $B(\lambda)$ sont isomorphes en tant que cristal. Les coordonn\'ees en corde d'un chemin int\'egral quelconque sont dans $B(\infty)$ et les coordonn\'ees en corde d'un chemin de $B_\pi$ sont dans $B(\lambda)$. 

\paragraph{Coordonn\'ees en corde d'un chemin brownien.} On consid\`ere le processus  $\{B(t)=(t,b_t): t\ge 0\}$ o\`u $\{b_t :t\ge 0\}$ est un brownien standard r\'eel standard. On d\'efinit, pour $t\ge 0$, une suite $\xi(t)$ de variables al\'eatoires positives  $(\xi_k(t))_{k\ge 0}$ par
  $$\mathcal P_{\alpha_n}\dots\mathcal P_{\alpha_0}B(t)=B(t)+\sum_{k=0}^n \xi_k(t)\alpha_k, \quad n\ge 0.$$
  Nous nous pla\c cons d\'esormais dans l'espace quotient $\widehat{\mathfrak{h}}^*_\R/\R\delta$   identifi\'e \`a $\R^2$.  Les racines et coracines s'\'ecrivent dans    l'espace quotient 
$$\alpha_0=(0,-2), \,\alpha_1=(0,2), \,\alpha^\vee_0=(1,-1), \, \alpha^\vee_1=(0,1).$$
Dans cet espace
  \begin{align*}
  \mathcal P_{\alpha_n}\dots\mathcal P_{\alpha_0}B(t)=\mathcal P_{s_n}\dots\mathcal P_{s_0}B(t),
  \end{align*}
  et 
  $$\mathcal L_{s_{n}}\mathcal P_{s_{n-1}}\dots\mathcal P_{s_0}B(t)=B(t)+\sum_{k=0}^{n-1} \xi_k(t)\alpha_k+\frac{1}{2}\xi_{n}(t)\alpha_{n}.$$
 \paragraph{Cristaux continus.} Pour d\'ecrire la loi des $\xi_k$, $k\ge 0$, on introduit des analogues continus des cristaux de Kashiwara. Pour une suite  $x=(x_k)\in \R_+^{\N}$ on pose, quand la limite existe
\begin{align}\label{sigma}\sigma(x)=\lim_{n \to +\infty} \sum_{k=0}^{n-1}x_k\alpha_k+\frac{1}{2}x_n\alpha_n.\end{align} 
On d\'efinit pour $\lambda\in \widetilde{ C}_{W} $, 
\begin{align*}&\Gamma(\infty)=\{x=(x_k)\in  \R^\N : \frac{x_k}{k}\ge \frac{x_{k+1}}{k+1}\ge 0,  \textrm{ pour tout } k\ge 1,\,  x_0\ge 0, \,  \sigma(x) \textrm{ existe}\},\\
&\Gamma(\lambda)=\{x\in \Gamma(\infty) : x_k\le   \langle\lambda-\sigma(x)+\sum_{i=0}^kx_i\alpha_i,\alpha^\vee_k\rangle , \textrm{ pour tout }  k\ge 0\}.
\end{align*}   
  \paragraph{Mesure de Duistermaat--Heckman.} Les convergences dans le th\'eor\`eme \ref{Pitaff} sont en fait presque s\^ures et si on pose  pour $t>0$,
  $$\Lambda(t)=\lim_{n\to\infty}\mathcal L_{s_{n+1}}\mathcal P_{s_n}\dots\mathcal P_{s_0}B(t).$$
 La loi  de $\xi(t)$ conditionnellement \` a $\{\Lambda(t)=\lambda\}$ fournit une mesure sur $\Gamma(\lambda)$.  Si l'on identifie les racines r\'eelles de $G_0$ \`a ses racines infinit\'esimales, l'image de cette mesure  par la fonction \og poids \fg \, $ \lambda-\sigma$ est une mesure de Duistermaat--Heckman. En rempla\c cant la fonction  \og poids \fg\, par le brownien, on obtient le th\'eor\`eme suivant. 
 \begin{theo}[Philippe Bougerol,  M.D, \cite{boubou-defo}]  La loi de $B(t)$ conditionnellement \`a $\Lambda(t)=\lambda$ est la mesure de Duistermaat--Heckman associ\'ee \`a l'action d'un tore maximal de $\SU(2)$ sur une orbite de $\widetilde{\mathcal O}_\lambda$ de $\widetilde{L}(\mathfrak{su}(2))^*$ munie d'une mesure de Frenkel.
 \end{theo}
 \section{Remarque sur le terme correcteur}\label{correc}
 L'introduction dans la section pr\'ec\'edente des cristaux de Kashiwara  et de leur version continue nous permet de faire quelques nouvelles remarques \`a propos du terme correcteur apparaissant dans le th\'eor\`eme de repr\'esentation \ref{Pitaff}. Nous avons dit qu'il provenait d'un manque de r\'egularit\'e des trajectoires browniennes. L'examen  des cristaux nous permet de nous forger une intuition plus alg\'ebrique de l'origine de cette correction.  On remarque en effet que  $x\in \Gamma(\infty)$ est dans $\Gamma(\lambda)$ si et seulement si pour tout $k\geq 0$,
\begin{align}\label{condG} \langle \sigma(x)-\sum_{i=0}^{k-1}x_i\alpha_i-\frac{1}{2}x_k\alpha_k,\alpha^\vee_k\rangle\leq \langle\lambda,\alpha^\vee_k\rangle .\end{align}
Ainsi $\Gamma(\lambda)$ est d\'efini de telle sorte que les conditions d'appartenance fassent appara\^itre un coefficient $1/2$, tout comme   la transformation de L\'evy. Cela sera d\'ecisif pour les preuves. Par ailleurs un cristal de plus haut poids $B(\lambda)$ pour $\lambda$ dominant   montre la m\^eme correction et si on suppose qu'un th\'eor\`eme de repr\'esentation du m\^eme type que celui que nous avons montr\'e existe pour une alg\`ebre de Kac--Moody affine quelconque, cette observation permet d'augurer la correction qu'il faudra effectuer. 
  \chapter{Produit de fusion et produit de convolution\\ \vspace{0.5cm}
\small{\it  o\`u  l'on fusionne  et convole}} \label{chap-fusion}
 Ce chapitre est un peu \`a part dans l'ensemble du m\'emoire dans la mesure o\`u les r\'esultats qui y sont \'enonc\'es n'interviennent pas dans  les diagrammes pr\'esent\'es dans l'introduction qui avaient jusqu'\`a pr\'esent, du moins je l'esp\`ere, donn\'e \`a l'ensemble son architecture. Il n'est cependant pas  sans rapport avec les questions que nous avons   abord\'ees pour l'instant.  
 Dans le cadre semi-simple  nous avons rappel\'e deux   approximations  de mesures standard obtenues \`a partir de   repr\'esentations. La premi\`ere implique les poids d'une repr\'esentation irr\'eductible dont le  plus haut poids est grand. Une mesure sur ces poids fournit une approximation de la mesure de Duistermaat--Heckman associ\'ee \`a une orbite coadjointe. La seconde implique de grandes puissances tensorielles d'une repr\'esentation fix\'ee. Leur  d\'ecomposition en composantes irr\'eductibles  fournit un processus de Markov qui est une approximation d'un processus \`a temps continu dans une chambre de Weyl. Nous pr\'esentons dans ce chapitre les   r\'esultats de \cite{defo0}. Ils rel\`event, dans le cadre de la fusion, des deux types de proc\'ed\'e mis en \oe uvre dans ces approximations.  Les orbites dans un groupe compact pour l'action par conjugaison du groupe sur lui-m\^eme remplacent ici les orbites coadjointes, le produit  de fusion tel qu'il est d\'efini dans \cite{Kac} remplace le produit tensoriel et des processus de Markov dans des alc\^oves, les processus dans des chambres de Weyl.  Nous d\'efinissons le produit de fusion et les coefficients de fusion dans la premi\`ere section, et  donnons dans la deuxi\`eme, notamment  \`a travers des exemples simples, une interpr\'etation probabiliste de ces coefficients.  Nous y voyons qu'il jouent pour une large classe de marches al\'eatoires dans des alc\^oves le m\^eme r\^ole que les coefficients de Littelwood-Richardson et leurs g\'en\'eralisations pour les marches dans une chambre de Weyl associ\'ee \`a une alg\`ebre de Lie semi-simple. Nous exposons dans la troisi\`eme section leur lien avec   la convolution sur un groupe  de Lie compact et   d\'ecrivons l'hypergroupe de la fusion comme une approximation de celui des classes de conjugaison sur ce groupe. Ce dernier aspect  est \`a rapprocher du travail de Eckhard Meinrenken   \cite{Meinrenken} portant sur des applications moments \`a valeurs dans un groupe compact.

Nous consid\'erons comme habituellement    un groupe de Lie $G_0$     simple, compact, connexe et simplement connexe, qu'on suppose sans perte de g\'en\'eralit\'e \^etre un groupe de matrices. Son alg\`ebre de Lie est not\'ee $\mathfrak g_0$. On consid\`ere un tore    maximal $T$ de $G_0$, d'alg\`ebre de Lie $\mathfrak t$.  On note $\Ad$ l'action adjointe de $G_0$ sur $\mathfrak g_0$, c'est-\`a-dire l'action par conjugaison de $G_0$ sur  $\mathfrak g_0$. On munit $\mathfrak g_0$ d'un produit scalaire $\Ad$-invariant not\'e $(\cdot\vert \cdot)$ et  on note $\mathfrak g =\mathfrak g_0\oplus i\mathfrak g_0$ l'alg\`ebre de Lie complexifi\'ee de $\mathfrak g_0$ et $\mathfrak h$ le complexifi\'e de $\mathfrak t$.    
 
\section{Produit de fusion et coefficients de fusion}   Les d\'efinitions suivantes se trouvent     dans le chapitre $13$ de \cite{Kac}.   Nous reprenons les notations de la section \ref{section-roots-coroots-A} de  notre chapitre \ref{chap-Algebre-affine}, que nous rappelons le plus souvent afin de faciliter la lecture. Nous munissons comme dans ce chapitre l'alg\`ebre de Cartan $\mathfrak h$ de la forme bilin\'eaire invariante standard normalis\'ee   $(\cdot\vert\cdot)$.   On rappelle que $\mathfrak h_\R=i\mathfrak t$. Pour tout  $y \in \mathfrak{h}_\R^*$, on \'ecrit $t_y$ pour la  translation d\'efinie sur  $\mathfrak{h}_\R^*$ par $$t_y(x)=x+y, \quad x\in \mathfrak{h}_\R^*.$$ Pour $k\in \N^*$, on consid\`ere  le groupe $W_k$ engendr\'e par  le groupe de Weyl $W$ associ\'e \`a $\mathfrak g$  et les translations $t_{(k+h^\vee)\theta}$, o\`u $\theta$ est la racine la plus haute de $\mathfrak g$ et $h^\vee$ le nombre de Coxeter dual qui vaut $1+\rho(\theta^\vee)$, $\rho$ \'etant la demi-somme des racines positives et $\theta^\vee$ la coracine la plus haute. En fait  $W_k$ est le produit semi-direct $W\ltimes T_{(k+h^\vee)M}$, o\`u $M=\nu(Q^\vee)$ et $T_{(k+h^\vee)M}=\{t_{(k+h^\vee)x}: x\in M\}$, $Q^\vee$ \'etant le r\'eseau des coracines et $\nu$ l'isomorphisme d\'efini par (\ref{iso}). Ainsi pour $w\in W_k$,  on d\'efinit  $\det(w)$ comme le  determinant de la composante lin\'eaire de $w$.   On rappelle que   $$\Pi^\vee=\{\alpha_i^\vee, i=1,\dots,n\}$$
est l'ensemble des coracines simples de $\mathfrak g$ suppos\'ee de rang $n$. Un domaine fondamental pour l'action de  $W_k$ sur $\mathfrak{h}_\R^*$ est 
$$A_k=\{\lambda\in \mathfrak{h}_\R^*:  \lambda(\theta^\vee)\le k+h^\vee \textrm{ et } 0\le \lambda(\alpha_i^\vee), i=1\dots,n\}.$$ 
On consid\`ere le r\'eseau des poids $P$ de $\mathfrak g$ et $P_+$ celui de ses poids dominants.
On introduit le sous-ensemble $P^k_+$ de $P_+$ d\'efini par $$P^k_+=\{\lambda\in P_+ : \lambda(\theta^\vee)\le k\},$$
qu'on appelle  alc\^ove de niveau  $k$. 
\begin{rema}\label{alc-aff}
Remarquons que, dans le cadre affine, $P_+^k$ s'identifie \`a l'ensemble des poids de $\widehat{P}_+$ dans $\R\Lambda_0\oplus \mathfrak{h}_\R^*$ dont la coordonn\'ee le long du poids fondamental $\Lambda_0$ est $k$. 
\end{rema}
D\'efinissons les caract\`eres discr\'etis\'es  $\Upsilon_\lambda$ de niveau $k$ pour $\lambda\in P$. Ils sont d\'efinis sur  l'alc\^ove discr\`ete $P_+^k$   ce qui justifie l'expression choisie pour leur nom. Pour  $\lambda\in P_+$,  on a
$$\Upsilon_\lambda(\sigma)=\mbox{ch}_{ \lambda }\big(-2i\pi\nu^{-1}(\frac{\sigma+\rho}{k+h^\vee})\big), \quad \sigma\in P_+^k,$$
 $\mbox{ch}_\lambda$ \'etant le caract\`ere de la repr\'esentation de plus haut poids $\lambda$ de $\mathfrak g$. Pour $\lambda\in P$, $\Upsilon_\lambda$ est d\'efini de m\^eme en utilisant pour $\mbox{ch}_\lambda$  l'expression donn\'ee par  la formule des caract\`eres de Weyl. On a  pour tout  $\lambda\in P$ et  $w\in W_k$
\begin{align}\label{alt}
\Upsilon_{w(\lambda+\rho)-\rho}=\det(w)\Upsilon_\lambda,
\end{align}
ce qui implique en particulier  que $\Upsilon_\lambda=0$ lorsque $(\lambda+\rho)$ est sur un mur $$\{x\in \mathfrak{h}_\R^*: x(\alpha^\vee)=0\}$$ pour un $\alpha^\vee\in \Pi^\vee$, ou sur le  mur $$\{x\in  \mathfrak{h}_\R^* : x(\theta^\vee)=k+h^\vee\}.$$   
Notons que pour $\lambda\in P_+^k$ on a 
\begin{align}\label{chi0}
\Upsilon_\lambda(0)=\frac{\prod_{\alpha\in \Phi_+}\sin(\pi(\lambda+\rho\vert \alpha)/(k+h^\vee))}{\prod_{\alpha\in \Phi_+}\sin(\pi( \rho\vert \alpha)/(k+h^\vee))},
\end{align} 
o\`u $\Phi_+$ est l'ensemble des racines positives. 
En particulier $\Upsilon_\lambda(0)$ est positif pour tout $\lambda\in P_+^k$.
  Les coefficients de fusion  $N^\beta_{\lambda,\gamma}$ de niveau $k$, pour $\lambda,\gamma,\beta\in P_+^k$, sont les constantes de structure de l'hypergroupe  des caract\`eres discr\'etis\'es de niveau $k$.  Ce sont les uniques entiers positifs tels  que
\begin{align}  
\forall \sigma\in P_+^k, \quad \Upsilon_{ \lambda }(\sigma)\Upsilon_{ \gamma }(\sigma)=\sum_{\beta\in P_+^k}N_{\lambda,\gamma}^\beta \Upsilon_{ \beta }(\sigma),\label{FP}.
\end{align}    
Ni l'unicit\'e, ni la positivit\'e des coefficients de fusion ainsi d\'efinis ne s'obtiennent imm\'ediatement. L'unicit\'e est d\'emontr\'ee dans \cite{Kac}. On trouvera une preuve de  la positivit\'e dans le chapitre 16 de \cite{Difrancesco}. Notons que cette d\'efinition des coefficients de fusion donn\'ee par Victor G. Kac dans \cite{Kac} co\"incide avec celle intervenant dans le cadre des mod\`eles de Wess-Zumino-Witten.
\paragraph{Lorsque  $G_0=\SU(2)$.}  Donnons l'exemple de $\mathfrak{sl}_2(\C)$ qu'on trouve dans \cite{Kac}.  Choisissons un entier $k$. Les plus hauts poids de $\mathfrak{sl}_2(\C)$ sont en correspondance avec les entiers naturels et l'alc\^ove de niveau $k$ avec $\{0,\dots,k\}$. Pour $n\in\{0,\dots,k\}$, le caract\`ere discr\'etis\'e $\Upsilon_n$ de niveau $k$  est ici donn\'e  par 
\begin{align}\label{cardSU(2)}
\Upsilon_n(m)=\sin[\frac{\pi(n+1)(m+1)}{k+2}]/\sin[\frac{\pi(m+1)}{k+2}],
\end{align} $m\in\{0,\dots,k\}$.
De plus pour $i,j\in \{0,\dots,k\}$, on a 
\begin{align}\label{FusionSU(2)}
\Upsilon_i\Upsilon_j=\sum_{s=0}^kN_{i,j}^s\, \Upsilon_s,
\end{align}
 o\`u
 $$N_{i,j}^s= \left\{
    \begin{array}{ll}
        1 & \mbox{si } \vert i-j\vert\le s\le \min(i+j,2k- i-j),\textrm{ et }  i+j+s\in 2\Z\\
        0 & \mbox{sinon.}
    \end{array}
\right.$$

\section{Marches dans une alc\^ove}
La remarque \ref{alc-aff} montre que le r\'eseau des poids dominants de l'alg\`ebre de Kac--Moody affine  $\widetilde{\mathcal L}(\mathfrak g)$ s'identifie \`a  la r\'eunion  des alc\^oves. Dans le chapitre \ref{chap-BETaff}, nous avons construit, en consid\'erant des produits tensoriels de repr\'esentations d'une alg\`ebre affine, des marches al\'eatoires \`a valeurs dans un tel r\'eseau   dont  le niveau augmente \`a chaque pas. Dans \cite{grabiner} David  Grabiner consid\`ere  des marches al\'eatoires dans des alc\^oves de niveau fix\'e afin d'\'etablir pour ces derni\`eres des formules de type Karlin--McGregor. Ce sont des marches au plus proche voisin pour lesquelles il \'etablit des formules donnant le nombre de trajectoires possibles assujetties \`a rester dans une alc\^ove,  la valeur initiale, la valeur finale et le nombre de pas \'etant fix\'es.    Il pose explicitement la question de leur lien \'eventuel avec la th\'eorie des repr\'esentations. Suivant une id\'ee de Ph. Bougerol, j'ai montr\'e dans    \cite{defo0} qu'un tel lien existait en effet. Il implique non plus le produit tensoriel de repr\'esentations mais le produit de fusion. Nous expliquons comment dans cette partie.\subsection{Principes g\'en\'eraux}
\paragraph{Coefficients de fusion et marche dans une alc\^ove.} Fixons un niveau $k\ge 1$ et  choisissons un poids $\gamma$ dans l'alc\^ove $P_+^k$. On associe aux coefficients de fusion  un noyau sous-markovien $K$ sur $P_+^k$ d\'efini par 
\begin{align}\label{killedM} 
K(\lambda,\beta)=\frac{1}{\dim V(\gamma)}N_{\lambda,\gamma}^\beta ,\quad \lambda,\beta\in P_+^k,
\end{align}
o\`u $V(\gamma)$ est une repr\'esentation de plus haut poids $\gamma$ de $\mathfrak g$.    Les r\`egles de fusion  donn\'ees par (\ref{FP})  permettent de d\'eterminer les fonction propres de $K$. Pour tout $\sigma\in P_+^k$, la fonction d\'efinie sur $P_+^k$ par 
$$x\in P_+^{k}\mapsto \Upsilon_x(\sigma)$$ est une fonction propre de $K$ associ\'ee \`a la valeur propre $\Upsilon_\gamma(\sigma)/\dim V(\gamma)$. On diagonalise ainsi $K$. On obtient la fonction propre de Perron-Frobenius pour $\sigma=0$ et on d\'efinit   une transformation de Doob $Q$ de $K$ en posant \begin{align}\label{qalcov} 
Q(\lambda,\beta) =N_{\lambda,\gamma}^\beta\frac{\Upsilon_\beta(0)}{\Upsilon_\lambda(0)\Upsilon_\gamma(0)}, \quad \lambda,\beta\in P_+^k.
\end{align}
On dit qu'un poids entier dominant $\gamma$ est minuscule lorsque l'orbite de $\gamma$ sous l'action du groupe de Weyl est  l'ensemble des poids de la repr\'esentation de plus haut poids $\gamma$. Pour un tel poids $\gamma$,  on a pour $\beta,\lambda\in P_+^k$,
  $$N_{\lambda,\gamma}^\beta= \left\{
    \begin{array}{ll}
        1 & \mbox{si } \beta-\lambda\in  W\cdot\{\gamma\} \textrm{ et } \beta\in P_+^k\\
        0 & \mbox{sinon.}
    \end{array}
\right.$$
Dans le cas minuscule donc, $K$ est   le noyau restreint \`a  $P_+^k$  d'une marche al\'eatoire \`a valeurs dans $P$ dont les pas sont uniform\'ement distribu\'es sur $W\cdot\{\gamma\}$.  En consid\'erant diff\'erentes alg\`ebres de Lie et diff\'erentes repr\'esentations minuscules on  obtient   les marches dans les alc\^oves consid\'er\'ees par David Grabiner et retrouve ainsi d'une mani\`ere unifi\'ee les formule de type Karlin--McGregor qu'il a \'etablies.  En effet, comme nous le voyons dans le paragraphe suivant,  dans la perspective de la fusion, ces formules   peuvent \^etre vues comme des formules du m\^eme  type valables pour les coefficients de fusion. 
 
\paragraph{Coefficients de fusion et nombre de trajectoires possibles dans une alc\^ove.}   D\'efinissons en effet  pour $p\ge 1$, et $\lambda,\gamma,\beta\in P_+^k$, les  coeffcients $N_{\lambda,\gamma,p}^\beta$ comme les uniques coefficients  positifs tels que pour tout $\sigma\in P_+^k$
$$\Upsilon_\lambda(\sigma)\Upsilon_\gamma(\sigma)^p=\sum_{\beta\in P_+^k}N_{\lambda,\gamma,p}^\beta\Upsilon_\beta(\sigma).$$ Si pour $\beta \in P$ on note  $m_{V(\gamma)^{\otimes p}}(\beta)$ la multiplicit\'e du poids $\beta $ dans le module $V(\gamma)^{\otimes p}$ de $\mathfrak g$, alors  on a l'identit\'e
\begin{align}\label{BKF}
N_{\lambda, \gamma,p}^ \beta=\sum_{w\in W_k} \det(w) m_{V(\gamma)^{\otimes p}}(w(\beta +\rho)-(\lambda+\rho)),
\end{align}
qui est l'analogue dans le cadre de la fusion de l'identit\'e (\ref{BrauerKlimyk}). On trouvera dans \cite{Walton} une preuve de cette formule pour $p=1$. La preuve est la m\^eme pour $p$ quelconque. Lorsque $\gamma$ est minuscule on montre facilement (voir proposition 5.2 de \cite{defo0}) que $N_{\lambda, \gamma,p}^ \beta$ est le nombre de trajectoires possibles d'une marche allant, en restant dans l'alc\^ove $P_+^k$,   de $\lambda$ \`a $\beta$ en $p$ pas de $W.\{\gamma\}$. Le coefficient $m_{V(\gamma)^{\otimes p}}(\beta)$ est le nombre de trajectoires possibles d'une marche allant de $\lambda$ \`a $\beta$ en $p$ pas de $W.\{\gamma\}$ mais pouvant sortir de l'alc\^ove. Ainsi dans le cas minuscule,  la formule (\ref{BKF}) est    une formule de type Karlin--McGregor pour la marche tu\'ee de noyau $K$.
 
\paragraph{Coefficients de fusion et asymptotique du nombre de trajectoires.} La  mesure $\mu$ d\'efinie sur $P_+^k$ par $$\mu(z)=\chi_z(0)^2,\quad z\in P_+^k,$$ est $Q$-invariante (proposition 5.6 de \cite{defo0})  et les r\'esultats standard sur les cha\^ines de Markov permettent d'obtenir une approximation de $N_{\lambda,\gamma,p}^\beta$ lorsque $p$ est grand  
proportionnelle \`a  
$$ \Upsilon_\gamma(0)^p\Upsilon_\lambda(0)\Upsilon_{\beta}(0),$$
le coefficient de proportionnalit\'e ne d\'ependant pas de $\gamma,\lambda,\beta$, et moyennant des conditions sur $\lambda-\beta$ et $p$ d\'ependant de la p\'eriodicit\'e de la cha\^ine $Q$ (proposition 5.7 de \cite{defo0}).
En particulier dans le cas o\`u  $\gamma$ est  minuscule,  on obtient une approximation pour le nombre de trajectoires possibles de la marche tu\'ee de noyau $K$ d\'efini par (\ref{killedM}), les points initiaux et finaux et le nombre de pas \'etant fix\'es. Nous retrouvons ainsi, toujours de mani\`ere unifi\'ee, les approximations obtenues par Christian Krattenthaler dans  \cite{Krattenthaler} pour une large classe de marches al\'eatoires. On pourra se reporter \`a la section $6$ de \cite{defo0} pour   des exemples d\'etaill\'es. Nous en d\'eveloppons deux dans la section suivante.

 \subsection{ Marches al\'eatoires dans une   alc\^ove de type $A$}  
\paragraph{Un premier exemple.} Fixons $d\ge 1$. Ici nous donnons l'exemple d'une marche al\'eatoire construite \`a partir d'une repr\'esentation minuscule de l'alg\`ebre de Lie $\mathfrak{sl}_d(\C)$ des matrices de taille $d\times d$ de trace nulle \`a coefficients dans $\C$.  Consid\'erons $\{e_1,\dots,e_d\}$ la base usuelle de $\R^d$. Une marche al\'eatoire sur $\Z^d$ dont les pas sont distribu\'es selon la mesure uniforme sur $\{e_1,\dots,e_d\}$ se d\'ecompose en une somme impliquant une marche al\'eatoire d\'eterministe et une marche al\'eatoire $\{X(k):k\ge 0\}$ dont les pas sont distribu\'es uniform\'ement sur $$\{e_1-\frac{1}{d}e,\dots,e_d-\frac{1}{d}e\},$$ o\`u $e=\sum_{i=1}^de_i$. Ainsi  l'ensemble des pas de  $\{X(k):k\ge 0\}$ est l'ensemble des poids de la repr\'esentation standard de $\mathfrak{sl}_d(\C)$. L'ensemble des racines de $\mathfrak{sl}_d(\C)$ est  $\Phi=\{e_i-e_j,i\ne j\}$,  l'ensemble des poids dominant est $$P_+=\{\lambda\in \R^d: \sum_{i=1}^d\lambda_i=0,\, \lambda_{i}-\lambda_{i+1}\in \N, \, i\in\{1,\dots,d-1\}\}.$$ 
 La plus haute coracine est $\theta^\vee=e_1-e_d,$  et pour $k\ge 1$, l'alc\^ove de niveau $k$ est $$P_+^k=\{\lambda\in P_+: \lambda_1-\lambda_d\le k\}.$$
 La demi-somme des racines positives vaut    $\rho=\frac{1}{2}\sum_{i=1}^d(d-2i+1)e_i$ et le nombre de Coxeter dual $h^\vee$ vaut $d$. La repr\'esentation standard de $\mathfrak{sl}_d(\C)$ a pour plus haut poids $e_1-\frac{1}{d}e$.   Pour $\lambda\in P_+^k$ et  $\gamma= e_1-\frac{1}{d}e $ dans (\ref{FP}) les coefficients de fusion valent 
  $$N_{\lambda,\gamma}^\beta= \left\{
    \begin{array}{ll}
        1 & \mbox{si } \beta-\lambda\in\{e_i-e/d: i\in\{1,\dots,d\}\} \textrm{ et }  \beta\in P_+^k \\
        0 & \mbox{sinon.}
    \end{array}
\right.$$
On consid\`ere $K$ le noyau de transition d\'efini par (\ref{killedM}) dans notre situation. C'est le noyau de la marche $\{X(k):k\ge 0\}$ tu\'ee en dehors de $P_+^k$, et le noyau markovien $Q$ de (\ref{qalcov}) vaut ici
$$Q(\lambda,\beta)=\frac{d\, \Upsilon_\beta(0)}{\Upsilon_\lambda(0)\Upsilon_\gamma(0)}K(\lambda,\beta), \quad \lambda,\beta\in P_+^k,$$
o\`u  
\begin{align}\label{chiA}
\Upsilon_x(0)=\prod_{1\le i<j\le d}\frac{\sin(\pi\frac{x_i-x_j+j-i}{k+d})}{ \sin(\pi\frac{j-i)}{k+d})}, \quad x\in P_+^k.\end{align}
   Moyennant des conditions sur $n\in \N$, $x,y\in P_+^k$, d\'ependant de la p\'eriodicit\'e de la cha\^ine, le nombre de trajectoires possibles d'une marche dans $P_+^k$ dont les pas sont dans  $\{e_1-\frac{1}{d}e,\dots,e_d-\frac{1}{d}e\}$, le point initial \'etant $x$, le point final  $y$, et nombre de pas $n$ \'etant fix\'e, est \'equivalent, d'apr\`es ce qu'on a vu,  pour $n$  grand (\`a une constante multiplicative pr\`es qui ne d\'epend pas de $(x,y)$) \`a 
 $$ \prod_{i=2}^{d} \frac{(\sin(\pi\frac{i}{d+k}))^n}{(\sin(\pi\frac{i-1}{d+k}))^n}  \prod_{1\le i<j\le d}\sin(\pi\frac{x_i-x_j+j-i}{k+d})\sin(\pi\frac{y_i-y_j+j-i}{k+d}).$$ 
 \paragraph{Variation sur le premier exemple.}
 Au lieu de la marche al\'eatoire pr\'ec\'edente, on peut consid\'erer une marche sur $P$ dont les pas sont  uniform\'ement distribu\'es sur $$\{\pm(e_i-e/d): i\in\{1,\dots,d\}\},$$  consid\'erer une telle marche tu\'ee en dehors de l'alc\^ove $P_+^k$ et poser pour celle-ci les m\^emes questions que pr\'ec\'edemment. La m\^eme m\'ethode peut \^etre mise en \oe uvre moyennant quelques ajustements. Il faut cette fois consid\'erer au lieu de la d\'ecompostion (\ref{FP}) avec $\gamma= e_1-e/d$, la d\'ecomposition  de  
 $$(\Upsilon_{e_1-e/d}+\Upsilon_{-(e_d-e/d)})\Upsilon_\lambda,$$ pour   $\lambda\in P_+^k$,  en somme de caract\`eres discr\'etis\'es, o\`u $\Upsilon_{-(e_d-e/d)}$ est le caract\`ere discr\'etis\'e associ\'e \`a la repr\'esentation de plus haut  poids  $ -(e_d-e/d)$ de $\mathfrak{sl}_d(\C)$, c'est-\`a-dire la repr\'esentation standard duale. Si on note $K$ le noyau de la marche tu\'ee en dehors de $P_+^k$ on peut exprimer cette  d\'ecomposition en fonction de $K$. On a en effet pour   $\lambda\in P_+^k$
  \begin{align}\label{fusionpm} 
  \frac{1}{2d}(\Upsilon_{e_1-e/d}+\Upsilon_{-(e_d-e/d)})\Upsilon_\lambda=\sum_{\beta \in P_+^k}K(\lambda,\beta)\Upsilon_\beta.
  \end{align}
 Comme dans l'exemple pr\'ec\'edent, les caract\`eres discr\'etis\'es donnent   les fonctions propres de l'op\'erateur $K$.  Les formules exactes ou asymptotiques donnant le nombre possible de trajectoires  de la marche dans une alc\^ove s'obtiennent de m\^eme. 
 \paragraph{Particules et ruine de joueurs.} Ce dernier exemple est li\'e \`a deux mod\`eles probabilistes  largement \'etudi\'es: un mod\`ele de particules auto-\'evitantes \'evoluant sur les sommets d'un polygone  et un mod\`ele de ruine de joueurs en dualit\'e avec le premier.    \begin{figure}
\begin{minipage}[c]{.46\linewidth}
     \begin{center}
             \includegraphics[width=6.1cm]{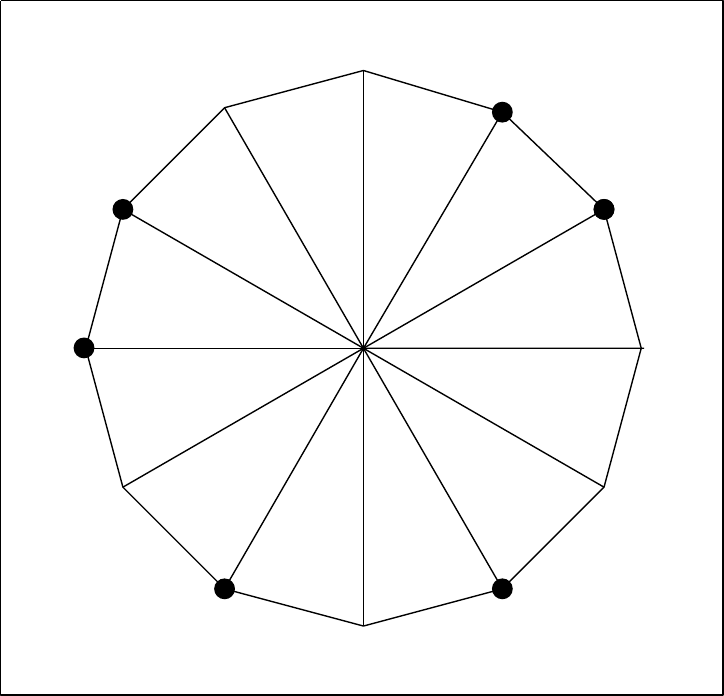}
         \end{center}
         \caption{Particules}
         \label{polyN1}
   \end{minipage} \hfill
   \begin{minipage}[c]{.46\linewidth}
    \begin{center}
            \includegraphics[width=6.1cm]{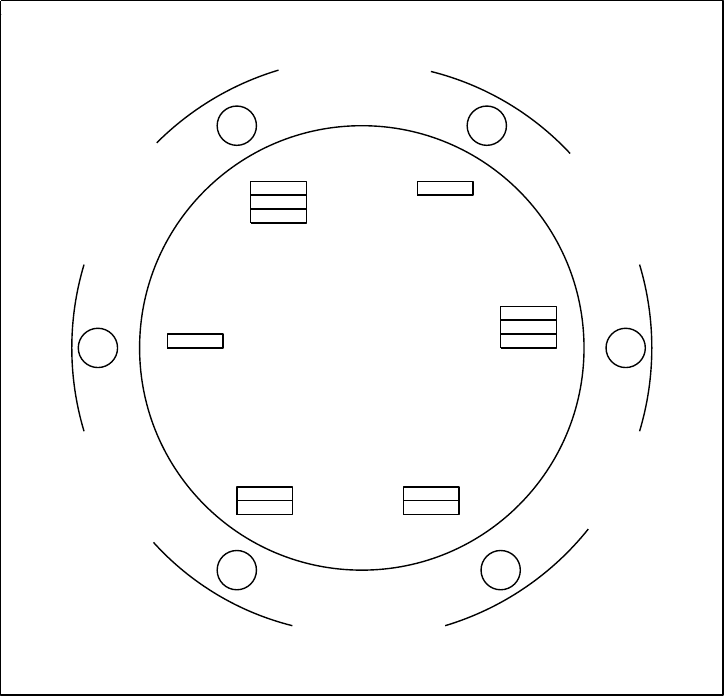}
        \end{center}
            \caption{Joueurs}
            \label{polyN0}
 \end{minipage}
 \end{figure}En effet, la d\'ecomposition (\ref{fusionpm}) permet de d\'eterminer les fonctions propres du noyau d'une marche  sym\'etrique simple sur le cercle discret $(\Z/N\Z)^d$ tu\'ee lorsque deux coordonn\'ees se rencontrent. Cette d\'ecomposition est ainsi li\'ee \`a un mod\`ele de particules assujetties \`a rester ordonn\'ees comme indiqu\'e \`a la figure \ref{polyN1}. Ce mod\`ele  est  par exemple  \'etudi\'e par Wolfgang K\"oning  et   Neil O'Connell  dans \cite{CK}  dans un tout autre contexte. Sur la figure, six particules \'evoluent sur un dod\'ecagone. Le mod\`ele de ruine de joueurs \`a $d$ joueurs  tel qu'il est d\'ecrit dans \cite{Diaconis} s'\'etudie de m\^eme dans le contexte de la fusion.  On suppose que $d$ joueurs   sont assis autour d'une table ronde. Chacun  poss\`ede un certain nombre de pi\`eces au temps initial et  le total des pi\`eces vaut $N$. \`A chaque instant entier, l'un d'entre eux est choisi selon une mesure de probabilit\'e uniforme, qui choisit lui-m\^eme de fa\c con \'equiprobable l'un de ses voisins. Les deux joueurs jouent ensuite \`a un jeu de pile ou face avec une pi\`ece \'equilibr\'ee et le vainqueur  re\c coit une pi\`ece du perdant. Le jeu continue tant qu'il reste \`a chacun de l'argent. Les deux mod\`eles sont en dualit\'e comme indiqu\'e aux figures \ref{polyN1} et \ref{polyN0}, les \'ecarts entre les particules devenant les sommes poss\'ed\'ees par chaque joueur.  Sur la figure \ref{polyN0}, six joueurs sont autour d'une table sur laquelle chacun a empil\'e ses pi\`eces devant lui. Le mod\`ele de la ruine de joueurs se comprend donc de m\^eme que celui des particules dans le contexte de la fusion.

\section{Un probl\`eme de Horn compact}

Le probl\`eme de Horn qui porte sur le spectre de la somme de deux matrices hermitiennes est li\'e, comme l'explique Klyachko dans \cite{Klyachko} et comme nous le rappelons plus bas, \`a la d\'ecomposition en composantes irr\'eductibles de produits tensoriels  de repr\'esentations du groupe unitaire.
Nous avons montr\'e dans \cite{defo0} que l'hypergroupe de la fusion fournissait quant \`a lui une approximation de celui des classes de conjugaison d'un groupe de Lie compact. Notre travail \'etablit donc un pont entre le produit de fusion d'une part et un probl\`eme de d\'ecomposition en classes de conjugaison d'autre part, qu'on peut qualifier de probl\`eme de Horn multiplicatif. Notre r\'esultat r\'esout une conjecture de \cite{Shaffaf}.

\paragraph{Probl\`eme de Horn.}
Le probl\`eme de Horn s'\'enonce ainsi : \'etant donn\'ees deux matrices hermitiennes de m\^eme taille dont seuls les spectres sont connus, que peut-on dire du spectre de leur somme ? Ce probl\`eme est li\'e via l'approximation semi-classique \`a une question de th\'eorie des repr\'esentations du groupe unitaire. Consid\'erons en effet deux matrices hermitiennes $A$ et $B$ de taille $d\times d$. Notons $a_1,\dots, a_d$  et $b_1, \dots, b_d$    les valeurs propres ordonn\'ees (compt\'ees avec  multiplicit\'e) de $A$ et de $B$   telles que $$a_1\ge\dots\ge a_d\quad \textrm{ et }\quad b_1\ge \dots\ge b_d.$$ Consid\'erons la d\'ecomposition
$$V\big([na]\big)\otimes V\big({[nb]}\big)=\bigoplus_{x\in E_n(A,B)}V(x)^{\oplus m_x},$$
o\`u $[na]=([na_1],\dots,[na_d])$, $[nb]=([nb_1],\dots,[nb_d])$ et o\`u  pour  $x\in \N^d$,  $V(x)$ est la repr\'esentation de plus haut poids $x$ de $\mbox{U}(d)$. Alors si $u$ est une variable al\'eatoire distribu\'ee selon la mesure de Haar sur $\mbox{U}(d)$, la mesure $$\sum_{x\in E_n(A,B)} \frac{m_x\dim V(x)}{\dim V\big([na]\big)\dim V\big([nb]\big)}\delta_{x/n},$$ converge  vers la loi des valeurs propres de $A+uBu^*$ quand $n$ tend vers l'infini.   Nous avons montr\'e dans \cite{defo0} un r\'esultat analogue dans lequel  le produit tensoriel est remplac\'e par le produit de fusion et la somme de matrices hermitiennes par un produit de matrices dans un groupe de Lie tel que $G_0$.

\paragraph{Probl\`eme de Horn compact.} 
On consid\`ere $$A=\{x\in \mathfrak h_\R^*:  x(\theta^\vee)\le  1, \, 0\le x(\alpha_i^\vee), \,i=1,\dots,n\}.$$
Le groupe $G_0$ \'etant simplement connexe, l'ensemble des classes de conjugaison dans $G_0$ est en correspondance bijective avec $A$. Choisissons $\xi$ et $\gamma$ deux \'el\'ements de  $A$ et $(\xi_{n})_{n\ge 1}$ et $(\gamma_{n})_{n\ge 1}$ deux suites \`a valeurs dans $P_{+}$   telles que  les suites $(\frac{1}{n}\xi_n)$ et $(\frac{1}{n}\gamma_n)$  convergent  respectivement   vers $\xi$ et $\gamma$ quand $n$ tend vers l'infini et telles que pour tout  $k\in \N^*$,
$\xi_{k}\in P_{+}^k$, $\gamma_{k}\in P_{+}^k$. On d\'efinit   une suite   $(\mu_{k})_{k\ge 1}$ de mesures de probabilit\'e  sur $A$ en posant 
\begin{align}\label{muk}
\mu_{k}=\sum_{\beta\in P^k_+}N^\beta_{\xi_k,\gamma_k}\frac{\Upsilon_\beta(0)}{\Upsilon_{\xi_k}(0)\Upsilon_{\gamma_k}(0)}\, \delta_{\frac{\beta+\rho}{k+h^\vee}}.
\end{align}
 Alors  la suite $( \mu_{k})_{k\ge 1}$   converge \'etroitement vers une mesure  $\mu_{\xi,\gamma}$ pour laquelle on a le th\'eor\`eme suivant.

\begin{theo}[M.D, \cite{defo0}] \label{theoOM}  Si    $u$ est une variable al\'eatoire distribu\'ee selon la mesure de Haar sur   $G_0$ alors  $$\exp(2i\pi\nu^{-1}(\xi))u\exp(2i\pi\nu^{-1}(\gamma))u^{*}$$ a la m\^eme loi que $$u\exp(2i\pi\nu^{-1}(\beta))u^{*},$$ o\`u  $\beta$ est distribu\'ee selon la loi  $\mu_{\xi,\gamma}$.
\end{theo}

Ce th\'eor\`eme r\'esout la conjecture de \cite{Shaffaf}. 
 Nous en donnons  une illustration  lorsque le groupe $G_0$ est le  groupe sp\'ecial unitaire  $\SU(2)$. On consid\`ere un tore maximal de $\SU(2)$
$$T=\{T_x= \left(
\begin{array}{cc}
e^{ 2i\pi x} &  0    \\
0  & e^{-2i\pi x}      
\end{array}
\right) : x\in[0,1]\}.$$
Choisissons un entier $k$. Les plus hauts poids de $\SU(2)$ sont les entiers naturels, l'alc\^ove de niveau $k$ est $\{0,\dots,k\}$ et pour $n\in\{0,\dots,k\}$, nous avons rappel\'e que le caract\`ere discr\'etis\'e $\Upsilon_n$ de niveau $k$  est dans ce cas donn\'e  par 
$$\Upsilon_n(m)=\sin[\frac{\pi(n+1)(m+1)}{k+2}]/\sin[\frac{\pi(m+1)}{k+2}],$$ $m\in\{0,\dots,k\}$.
Pour $i,j\in \{0,\dots,k\}$, on a 
\begin{align}\label{FusionSU(2)}
\Upsilon_i\Upsilon_j=\sum_sN_{i,j}^s\, \Upsilon_s,
\end{align}
 o\`u
 $$N_{ij}^s= \left\{
    \begin{array}{ll}
        1 & \mbox{si } \vert i-j\vert\le s\le \min(i+j,2k- i-j),\textrm{ et }  i+j+s\in 2\Z\\
        0 & \mbox{sinon.}
    \end{array}
\right.$$
On choisit  $a$ et $b$ dans l'intervalle $]0,1[$. La mesure $\mu_k$ d\'efinie pour $\xi_k=[ka]$ et $\gamma_k=[kb]$   vaut   ici
$$\mu_k=\sum_{s=0}^kN_{[ka],[kb]}^s \frac{\sin[\pi(s+1)/(k+2)]\sin[\pi/(k+2)]}{\sin[\pi([ka]+1)/(k+2)]\sin[\pi([kb]+1)/(k+2)]}\delta_{\frac{s+1}{k+2}}$$
Par ailleurs pour tout  $X$ dans $\SU(2)$ il existe un unique $x\in [0,1]$ tel que  $X=uT_{x/2}u^{-1}$ pour un  $u$ dans  $\SU(2)$. On appelle $x$ la partie radiale   de $X$. On montre par un calcul \'el\'ementaire que  si   $u$ est choisi selon la mesure de Haar sur $\SU(2)$ alors la partie radiale de $uT_{a/2}u^{-1}T_{b/2}$, pour $a,b\in[0,1]$, a une densit\'e d\'efinie sur $\R$ par 
$$\frac{1}{2} \frac{\pi \sin(\pi x)}{\sin(\pi a)\sin(\pi b)}1_{[r,s]}(x), \quad x\in \R,$$
o\`u  $r=\min(\vert a-b\vert,\min(a+b,2-(a+b)))$, $s= \max(\vert a-b\vert,\min(a+b,2-(a+b)))$.
C'est aussi ce que dit le th\'eor\`eme \ref{theoOM}. Cet exemple est \`a comparer \`a l'exemple de  $\SU(2)$ de \cite{Shaffaf} (voir aussi \cite{Dooley}).   
  \chapter{Bruit blanc  quantique \`a temps discret \\ \vspace{0.5cm}
\small{\it  o\`u l'on fait des projets d'avenir}} \label{chap-non-com}
 
Les diff\'erentes questions que nous avons abord\'ees dans ce m\'emoire t\'emoignent toutes d'un lien entre certaines mesures sur les orbites coadjointes dans le dual d'une alg\`ebre de Lie  et les repr\'esentations de cette alg\`ebre. Il me semble que dans le cas des alg\`ebres de Lie semi-simples complexes c'est pour un probabiliste l'approche non commutative d\'evelopp\'ee par Philippe Biane   \cite{biane1} qui r\'ev\`ele de la mani\`ere la plus lumineuse la nature de ce lien. En effet, pour ces alg\`ebres et dans un contexte non commutatif, les mesures sur les orbites     apparaissent comme des limites en loi de variables al\'eatoires non commutatives et ainsi certaines de leurs propri\'et\'es comme l'\'echo de ph\'enom\`enes li\'es par nature \`a la th\'eorie des repr\'esentations.  Je  m\`ene aujourd'hui une r\'eflexion - pr\'ecisons que c'est une r\'eflexion collective - visant \`a d\'evelopper une analyse semblable pour des alg\`ebres affines.   Il s'agit principalement de donner un fondement alg\'ebrique \`a la formule des caract\`eres de Kirillov--Frenkel discut\'ee au chapitre \ref{chap-loc} et de mieux comprendre, dans ce contexte alg\'ebrique, le r\^ole de l'enroulement du brownien. 
J'explique bri\`evement dans  ce chapitre le cadre dans lequel nous - un nous exclusif donc ici - avons inscrit notre r\'eflexion.

Avant cela  je fais quelques rappels sur l'approche non commutative,  reprenant le formalisme de Beno\^it Collins et Piotr \'Sniady  \cite{CS}. Ce formalisme fait apparaitre naturellement des variables al\'eatoires \`a valeurs dans des espaces duaux et    me semble ainsi s'accorder  tout \`a fait \`a la philosophie de la m\'ethode des orbites de Kirillov. 

Notons que le point de vue non commutatif est celui que Fran\c cois Chapon et moi-m\^eme  avions adopt\'e  dans \cite{chapondefo} donnant ainsi une interpr\'etation   alg\'ebrique aux r\'esultats de \cite{ANM}. Nous ne parlerons cependant pas de ce travail ici. 

\section{Variables al\'eatoires non commutatives}
Consid\'erons l'ensemble $\mathcal L^{-\infty}(\Omega)$ des variables al\'eatoires sur $\Omega$ \`a valeurs r\'eelles admettant des moments de tous les ordres. Alors l'esp\'erance $\E$ d\'efinit une forme lin\'eaire positive sur $\mathcal L^{-\infty}(\Omega)$ v\'erifiant $\E(1)=1$ et l'ensemble $\mathcal L^{-\infty}(\Omega)$ est ainsi une $\R$-alg\`ebre unitaire involutive munie  d'une forme lin\'eaire   normalis\'ee et positive. On a plus g\'en\'eralement la d\'efinition suivante. 
\begin{defn}  Un espace de probabilit\'e non commutatif est un couple  $(\mathcal U,\varphi)$, o\`u   $\mathcal U$ est une $\R$-alg\`ebre unitaire  munie d'une involution $*$ et  $\varphi$ est un \'etat, c'est-\`a-dire une forme lin\'eaire  sur $\mathcal U$,  normalis\'ee, positive et traciale, i.e. v\'erifiant $\varphi(1)=1$,   $\varphi(xx^*)>0$ pour tout \'el\'ement $x$ non nul de $\mathcal U$ et $\varphi(xy)=\varphi(yx)$ pour tout $x,y$ dans  $U$.  Un \'el\'ement de $\mathcal U$ est une variable al\'eatoire non commutative.
\end{defn}

 \paragraph{Vecteur al\'eatoire non commutatif.} Soit $V$ un $\R$-espace vectoriel de dimension finie et $v:\Omega\to V$ un vecteur al\'eatoire usuel. On d\'efinit le moment d'ordre $k\in \N$ de $v$ en posant 
 $$m_k(v)=\E(v\otimes\dots\otimes v)\in V^{\otimes k}.$$
Autrement dit si $\mathcal B=\{e_1,\dots,e_n\}$ est une base de $V$ et 
$v=\sum_{i=1}^na_ie_i,$
avec $a_i\in \mathcal L^{-\infty}(\Omega)$, alors 
$$m_k(v)=\sum_{i_1,\dots,i_k} \E(a_{i_1}\dots a_{i_k})e_{i_1}\otimes\dots\otimes e_{i_k}.$$
Ainsi la donn\'ee des moments   de $v$ est la donn\'ee des moments joints   de ses coordonn\'ees dans la base $\mathcal B$. Le vecteur al\'eatoire $v$ peut \^etre vu comme un \'el\'ement de $\mathcal L^{-\infty}(\Omega)\otimes V$. Pour un espace de probabilit\'e non commutatif $(\mathcal U,\varphi)$, on appelle vecteur al\'eatoire non commutatif \`a valeurs dans $V$ un \'el\'ement $v\in \mathcal U\otimes V$. On d\'efinit le moment d'ordre $k$ d'un \'el\'ement
$v=\sum_ia_i\otimes e_i$ de $\mathcal U\otimes V$ par 
$$m_k(v)=\sum_{i_1,\dots,i_k} \varphi(a_{i_1}\dots a_{i_k})e_{i_1}\otimes\dots\otimes e_{i_k}.$$  
Ainsi si $v$ est pens\'e comme une collection $\{a_i:i\in\{1,\dots,n\}\}$ de variables al\'eatoires non commutatives, alors les moments de $v$ fournissent  les moments joints non commutatifs de ces variables al\'eatoires. 
\begin{defn} Soit $V$ un espace  vectoriel sur $\R$, $ (\mathcal A,\varphi) $ un espace de  probabilit\'e non commutatif et $\{(\mathcal A_n,\varphi_n) : n\ge 0\}$  une famille d'espaces de  probabilit\'e non commutatifs. 
\begin{enumerate} 
\item On dit que deux vecteurs al\'eatoires $v_1\in \mathcal A_1\otimes V $, $v_2\in \mathcal A_2\otimes V$ ont m\^eme loi s'ils ont les m\^emes moments. 
\item On dit qu'une suite  de vecteurs al\'eatoires $\{v_n:n\ge 0\}$ telle que $v_n\in \mathcal A_n\otimes V$ converge en loi vers $v\in \mathcal A\otimes V$ si les moments de $v_n$ convergent vers les moments de $v$. 
\end{enumerate}
\end{defn}
\begin{defn} Soit $V$ un espace  vectoriel sur $\R$  et $\{(\mathcal A_t,\varphi_t) : t\in T\}$  une famille d'espaces de  probabilit\'e non commutatifs, avec $T=\N$, ou $T=\R_+$, et $\mathcal A_s\subset \mathcal A_t$, pour $s\le t$.  Un processus non commutatif est une famille de vecteurs non commutatifs $\{v_t : t\in T\}$  avec $v_t\in \mathcal A_t\otimes V$.   
\end{defn}

On dit qu'une famille $\{v_t^n:t\in T\}$, $n\ge 0$, de processus non commutatifs converge en loi vers un processus non commutatif $\{v_t: t\in T\}$ si pour    $t_1,\dots,t_k\in T$, le vecteur non commutatif $(v_{t_1}^n,\dots, v_{t_k}^n)$ converge en loi vers  $(v_{t_1},\dots, v_{t_k})$ quand $n$ tend vers l'infini.

\section{Repr\'esentations et mesures d'orbites}
Un exemple classique de vecteur al\'eatoire non commutatif  est celui d'un vecteur al\'eatoire  obtenu \`a partir d'une repr\'esentation d'alg\`ebre de Lie simple complexe $\mathfrak{g}$. On pourra le trouver d\'ecrit dans \cite{CS}. Un tel  vecteur al\'eatoire est \`a valeurs dans  le  dual $\mathfrak{g}_0^*$ de la partie compacte $\mathfrak g_0$ de   $\mathfrak{g}$.    L'avantage du cadre non commutatif est qu'il fournit des approximations de mesures classiques sur $\mathfrak{g}_0^*$ qui partagent avec elles une propri\'et\'e d'invariance en loi pour l'action coadjointe sur $\mathfrak{g}_0^*$ du groupe $G_0$ dont $\mathfrak{g}_0$ est l'alg\`ebre de Lie.  

Il y a deux grands types  d'approximation. D'abord, une repr\'esentation irr\'eductible de grande dimension fournit une approximation de la mesure uniforme  sur une orbite coadjointe dans $\mathfrak{g}_0^*$. Dans cette perspective, la formule des caract\`eres de Kirillov appara\^it comme une d\'eg\'en\'erescence commutative d'une formule analogue impliquant le caract\`ere de la repr\'esentation et valable dans un contexte non commutatif.   D'autre part, une puissance tensorielle de repr\'esentations de degr\'e \'elev\'e fournit  une approximation de la mesure gaussienne sur $\mathfrak{g}_0^*$ et une famille de  produits tensoriels de repr\'esentations une  approximation d'un brownien sur $\mathfrak{g}_0^*$.  Nous avons rappel\'e au chapitre \ref{chap-drap} que le processus de la partie radiale d'un mouvement brownien sur $\mathfrak{g}_0$  a la m\^eme distribution qu'un brownien sur le dual d'une sous-alg\`ebre de Cartan, conditionn\'e \`a rester dans une chambre de Weyl associ\'ee au syst\`eme de racines de $\mathfrak{g}_0^*$. Dans un cadre non commutatif,   en \'evaluant les moments de l'approximation   du brownien  sur $\mathfrak{g}^*_0$ sur le centre de l'alg\`ebre enveloppante de $\mathfrak{g}$, on obtient une version non commutative de cette identit\'e en loi, qui appara\^it donc elle  aussi comme l'\'echo   d'un ph\'enom\`ene discret li\'e \`a la th\'eorie des repr\'esentations.   
 \section{Repr\'esentations d'alg\`ebres affines et mesure de Wiener}
  Igor  Frenkel a \'etabli pour les alg\`ebres affines une formule des caract\`eres de type Kirillov.  Comme dans le cas semi-simple on peut associer \`a certaines repr\'esentations   d'une telle alg\`ebre, une famille de variables al\'eatoires. Il s'agit cependant maintenant d'une famille infinie. Dans le travail de Frenkel,  c'est la mesure de Wiener sur la forme compacte $\mathfrak{g}_0$ qui fournit   une mesure  sur une orbite coadjointe dans $\widetilde L(\mathfrak g_0)^*$. On peut construire    \`a partir d'une repr\'esentation de $\widehat{\mathcal{L}}(\mathfrak{g})$ une famille de variables al\'eatoires  non commutatives  qui partage avec la mesure de Wiener une propri\'et\'e de quasi-invariance en loi   rempla\c cant la propri\'et\'e d'invariance en loi valable dans le cas semi-simple.   Nous indiquons dans cette partie comment sont construites ces variables al\'eatoires non commutatives. Il est probable qu'elles  fournissent une approximation de la mesure de Frenkel sur une orbite coadjointe. Nous consid\'erons ensuite le cas d'une puissance tensorielle de repr\'esentations. 
 Dans le contexte d'une alg\`ebre semi-simple complexe, une puissance tensorielle de repr\'esentations fournit  une approximation de la gaussienne sur la forme compacte de cette alg\`ebre. Une famille de produits tensoriels de repr\'esentations fournit quant \`a elle une approximation du brownien sur la forme compacte. Dans le cas des alg\`ebres affines, la m\^eme construction est possible. Le bruit blanc sur la forme compacte $\mathfrak{g}_0$ remplace la gaussienne, et une famille de bruits blancs associ\'ee \`a un drap brownien sur $\mathfrak{g}_0$ remplace le brownien.

 \subsection{Repr\'esentations d'alg\`ebres affines et variables al\'eatoires non commutatives} 
Consid\'erons l'alg\`ebre affine $\widehat{\mathcal{L}}(\mathfrak g)$ associ\'ee \`a $\mathfrak{g}$ d\'efinie dans le chapitre \ref{chap-Algebre-affine}. Sa forme forme compacte $\widehat{\mathcal{L}}(\mathfrak{g})_0$   s'\'ecrit 
$$\widehat{\mathcal{L}}(\mathfrak{g})_0=\mathcal{L}(\mathfrak{g})_0\oplus \R d\oplus \R c,$$ o\`u  $\mathcal{L}(\mathfrak{g})_0$ s'identifie \`a un sous-ensemble de boucles \`a valeurs dans $\mathfrak{g}_0$.   L'alg\`ebre  $\mathfrak{g}_0$ est munie d'un produit scalaire $\Ad(G_0)$-invariant et nous   munissons $\mathcal{L}(\mathfrak{g})_0$ d'un produit scalaire $(\cdot\vert \cdot)$ d\'efini par 
 $$(x,y)=\int_0^1(x_s\vert y_s)\, ds,$$ $x,y\in \mathcal{L}(\mathfrak{g})_0$ et choisissons   une base orthonorm\'ee $\{e_k: k\in \Z\}$ de $\mathcal{L}(\mathfrak{g})_0$. Consid\'erons   une repr\'esentation $(\rho,V)$  de $\widehat{\mathcal{L}}(\mathfrak{g})$ de la cat\'egorie $\mathcal O_{int}$\footnote{On suppose que $V$ est un produit tensoriel de repr\'esentations irr\'eductibles de plus haut poids de sorte que l'\'etat  $\tau^r_\rho$ sera bien d\'efini.}. Elle fournit une famille de variables al\'eatoires  non commutatives.   Consid\'erons en effet  l'alg\`ebre enveloppante $\mathcal U(\widehat{\mathcal{L}}(\mathfrak{g})_0)$ de $\widehat{\mathcal{L}}(\mathfrak{g})_0$ et  l'alg\`ebre enveloppante  vue dans cette repr\'esentation   $$\mathcal U_\rho(\widehat{\mathcal{L}}(\mathfrak{g})_0)=\{\rho(x): x\in \mathcal U(\widehat{\mathcal{L}}(\mathfrak{g})_0)\},$$
 munie de l'involution $\rho(x)^*=-\rho(x)$.  Pour un r\'eel strictement positif $r$, nous munissons $\mathcal U_\rho(\widehat{\mathcal{L}}(\mathfrak{g})_0)$    d'un \'etat $\tau^r_\rho$ en posant
 $$\tau^r_\rho(\rho(x))=\frac{\mbox{Tr}(e^{r\rho(d)}\rho(x))}{\mbox{Tr}(e^{r\rho(d)})}, \, \, x\in \mathcal U(\widehat{\mathcal{L}}(\mathfrak{g})_0).$$
La famille
 $\{\rho(e_k): k\in \Z\}$ forme alors une famille de variables al\'eatoires non commutatives sur $(\mathcal U_\rho(\widehat{\mathcal{L}}(\mathfrak{g})_0),\tau^r_\rho)$, qu'on peut \'ecrire sous la forme d'une somme formelle $$ \rho=\sum_{k\in\Z}\rho(e_k)e^*_k+\rho(c)\Lambda_0.$$
 Remarquons que la repr\'esentation $\rho$ est ici restreinte \`a l'alg\`ebre affine sans d\'erivation. Cette d\'erivation est cependant essentielle pour que la trace d\'efinissant l'\'etat existe.  La loi des variables al\'eatoires   $\{\rho(e_k): k\in \Z\}$  est d\'etermin\'ee par les moments d'ordre $k$
 $$m_k(\rho)=\sum_{i_1,\dots,i_k}\tau^r_\rho(e_{i_1}\dots e_{i_k})e^*_{i_1}\otimes \dots\otimes e_{i_k}^*$$ 
On exclut volontairement   la variable $\rho(c)$ qui commute avec toutes les autres\footnote{Nous ne consid\'ererons que des repr\'esentations dont les poids sont tous de m\^eme niveau et $\rho(c)$ sera donc toujours d\'eterministe.}.  La loi du  vecteur construit ici satisfait une propri\'et\'e de quasi-invariance pour l'action coadjointe d'un sous-groupe de lacets de $L(G_0)$.    Consid\'erons en effet le groupe de Lie complexe simplement connexe  $G$ associ\'e \`a $\mathfrak g$. Nous notons $L(G)$ l'ensemble des boucles \`a valeurs dans $G$ et $\widehat{L}(\mathfrak g)$ l'alg\`ebre de Lie complexifi\'ee de $\widehat{L}(\mathfrak g_0)$. On peut d\'efinir\footnote{On  trouvera   les d\'etails dans \cite{Frenkel} .} une action adjointe $\Ad$    de $L(G)$ sur $\widehat{ {L}}(\mathfrak{g})$ qui co\"incide avec celle de $L(G_0)$ sur $\widehat{ {L}}(\mathfrak{g}_0)$ d\'efinie au chapitre \ref{chap-loc}.  
Si $z$ est un \'el\'ement d'un sous-espace radiciel de  $\mathcal{L}(\mathfrak{g})$ autre que la sous-alg\`ebre de Cartan, alors $\ad(z)$ et 
 $\rho(z)$ agissent de mani\`ere  localement nilpotente  et l'action adjointe satisfait
 $$\rho(\Ad(\exp(z))(x))=\exp(\rho(z))\rho(x)\exp(-\rho(z)).$$
 Ainsi si $\gamma$ est un \'el\'ement de $L(G_0)$ s'\'ecrivant 
 $$\gamma=\exp(z_1)\dots\exp(z_n),$$ o\`u  $z_1,\dots,z_n$ sont   des \'el\'ements de sous-espaces radiciels de  $\mathcal{L}(\mathfrak{g})$  alors  pour    $k\in\N$ et $x_1,\dots,x_k\in \mathcal{L}(\mathfrak{g})_0$  on a 
 \begin{align*} 
\Ad^*(\gamma)m_k(\rho)(x_1\otimes\dots\otimes  x_k)&=\tau^r_\rho(\rho(\Ad(\gamma^{-1}) (x_1))\dots\rho(\Ad(\gamma^{-1})( x_k)))\\ 
& =\tau^r_\rho(e^{-\frac{r}{2}(\gamma'\gamma^{-1},\gamma'\gamma^{-1})\rho(c)-r\rho(\gamma'\gamma^{-1})}\rho(x_1)\dots\rho(x_k)).
\end{align*}
Finalement   si $\rho$ est une repr\'esentation de plus haut poids de niveau $n$ alors
 \begin{align*} 
\Ad^*(\gamma)m_k(\rho)(x_1\otimes\dots\otimes  x_k) =\tau^r_\rho(e^{-\frac{rn}{2}(\gamma'\gamma^{-1},\gamma'\gamma^{-1})-\rho(\gamma'\gamma^{-1})}\rho(x_1)\dots\rho(x_k)).
\end{align*} 
Cette derni\`ere identit\'e devra \^etre compar\'ee \`a l'identit\'e (\ref{QIC}).
 \subsection{Produit tensoriel de repr\'esentations d'alg\`ebres affines et bruit blanc non commutatif}
 Commen\c cons par rappeler que la d\'eriv\'ee d'un brownien $\{x_t:t\ge 0\}$  \`a valeurs sur $\mathfrak{g}_0$, si elle n'existe pas, s'exprime usuellement comme une somme formelle 
  \begin{align} \label{BB} dx_s=\sum_{k\in \Z}a_k e^*_k,
  \end{align}
  la famille $\{a_k: k\in \Z\}$ \'etant une famille de gaussiennes centr\'ees ind\'ependantes et identiquement distribu\'ees.  Cette somme formelle est le bruit blanc. Si elle  n'est pas convergente  pour une norme $L^2$, elle      a cependant un sens dans le dual d'un  certain espace fonctionnel.   Nous avons pu montrer que  pour une repr\'esentation $\rho$ de $\widehat{\mathcal{L}}(\mathfrak{g})$ particuli\`ere de niveau $1$,   la suite de variables al\'eatoires 
 $$\{\frac{1}{n}\rho^{\otimes n}(e_k) : k\in \Z\}, \, n\ge 1,$$ 
 les variables al\'eatoires $\rho^{\otimes n}(e_k)$, $k\in\Z$, \'etant d\'efinies sur l'espace de probabilit\'e non commutatif $(\mathcal U_{\rho^{\otimes n}}(\widehat{\mathcal{L}}(\mathfrak{g})_0),\tau^{r}_{\rho^{\otimes n}})$, avec $r=1/n$, converge\footnote{\`a un facteur $i$ pr\`es qu'on peut omettre moyennant  quelques identifications.} quand $n$ tend vers l'infini vers la famille de variables al\'eatoires $\{a_k: k\in \Z\}$. La repr\'esentation $\rho^{\otimes n}$ doit donc \^etre pens\'ee comme une approximation non commutative du bruit blanc sur $\mathfrak{g}_0$ donn\'e en (\ref{BB}), ou plut\^ot, si on veut garder l'information de niveau, une approximation de la forme  
  $$\Phi=\int_0^1(\cdot\vert dx_s)+\Lambda_0.$$ 
Cette approximation partage avec $\Phi$ une propri\'et\'e de quasi-invariance en loi.  En effet les moments de $\Phi$    s'\'ecrivent
$$m_k(\Phi)=\sum \E(\Phi(e_{i_1})\dots\Phi(e_{i_k}))e_{i_1}^*\otimes\dots\otimes e^*_{i_k},\, k\in\N,$$    et pour $\gamma\in L(G_0)$, $\Ad^*(\gamma)m_k(\Phi)(x_1\otimes\dots\otimes  x_k)$ vaut 
\begin{align}\label{QIC} \E(e^{-\frac{1}{2}\int_0^1(\gamma'_s\gamma_s^{-1}\vert\gamma'_s\gamma_s^{-1})\, ds-\int_0^1(\gamma'_s\gamma_s^{-1}\vert \, dx_s)}\Phi(x_1)\dots\Phi(x_k)),
\end{align}
ce qui  traduit le caract\`ere quasi-invariant de la mesure de Wiener.    Quant \`a la famille de processus 
 $$\{\frac{1}{n}\rho^{\otimes [nt]}: t\ge 0\},\, n\ge 1,$$
 elle fournit une approximation du processus  $$\{t\Lambda_0+\int_0^1(\cdot\vert \,dx_s^t): t\ge 0\},$$ consid\'er\'e au chapitre \ref{chap-drap}, o\`u $\{x_s^t: s,t\ge 0\}$ est un drap brownien  standard sur $\mathfrak{g}_0$.

\end{document}